\documentclass[11pt]{amsart}
\usepackage{amssymb}
\usepackage{enumerate}
\usepackage{amsmath,stmaryrd,amssymb,amsfonts,amsthm,
	latexsym, amscd, epsfig, epic,color,bbm}
\usepackage{mathtools, graphics, epstopdf} 
\usepackage{extarrows}
\usepackage{enumitem}
\usepackage{enumerate}
\usepackage{mathrsfs}
\usepackage{eucal}
\usepackage{eufrak}
\usepackage{xspace}
\usepackage{a4wide}
\usepackage{colordvi}
\usepackage{comment}
\usepackage{hyperref,eucal}
\usepackage{tikz}\usetikzlibrary{matrix,arrows}
\usepackage{imakeidx}
\makeindex

\newtheorem{Theorem}{Theorem}[section]
\newtheorem{theorem}[Theorem]{Theorem}
\newtheorem{corollary}[Theorem]{Corollary}
\newtheorem{lemma}[Theorem]{Lemma}
\newtheorem{definition}[Theorem]{Definition}

\newtheorem{proposition}[Theorem]{Proposition}

\newtheorem{example}[Theorem]{Example}

\theoremstyle{plain}

\def\E{\CMcal{E}}

\def\bt{\mathbf{t}}
\def\bx{\mathbf{x}}

\newcommand{\Z}{\mathbb{Z}}

\newcommand{\fP}{\mathfrak{P}}
\newcommand{\LTer}{\mathrm{LTer}}
\newcommand{\Ter}{\mathrm{Ter}}

\newcommand{\For}{\mathrm{For}}
\newcommand{\Trim}{\mathrm{Trim}}
\newcommand{\NCPerm}{\mathrm{NCPerm}}
\newcommand{\fS}{\mathfrak{S}}
\newcommand{\id}{\mathrm{id}}

\theoremstyle{remark}

\numberwithin{equation}{section}
\numberwithin{figure}{section}

\newcommand{\FPExampleActive}[1]{%
	\draw[red] (#1) circle (10pt);
	\filldraw[red] (#1) circle (5pt);
}
\newcommand{\FPExampleForest}[3]{%
	\begin{tikzpicture}[scale=0.28,baseline=(current bounding box.center)]
		\path[use as bounding box] (-1,-2.6) rectangle (9,3.7);
		\draw[black,thin] (-0.8,0) -- (7,0);
		\ifcase#1\relax
		\or % F: (1 vee 2), 3, 4, ...
		\draw[blue,thick] (0,0) -- (1,1) -- (2,0);
		\filldraw[black] (1,1) circle (5pt);
		\or % T: ((1 vee 2) vee 3), 4, ...
		\draw[blue,thick] (0,0) -- (2,2) -- (4,0);
		\draw[blue,thick] (2,0) -- (1,1);
		\filldraw[black] (1,1) circle (5pt);
		\filldraw[black] (2,2) circle (5pt);
		\or % U: (1 vee (2 vee 3)), 4, ...
		\draw[blue,thick] (0,0) -- (2,2) -- (4,0);
		\draw[blue,thick] (2,0) -- (3,1);
		\filldraw[black] (2,2) circle (5pt);
		\filldraw[black] (3,1) circle (5pt);
		\or % C_3: 1 vee ((2 vee 3) vee 4)
		\draw[blue,thick] (0,0) -- (3,3) -- (6,0);
		\draw[blue,thick] (2,0) -- (4,2);
		\draw[blue,thick] (4,0) -- (3,1);
		\filldraw[black] (3,3) circle (5pt);
		\filldraw[black] (4,2) circle (5pt);
		\filldraw[black] (3,1) circle (5pt);
		\or % C_4: ((1 vee 2) vee 3) vee 4
		\draw[blue,thick] (0,0) -- (3,3) -- (6,0);
		\draw[blue,thick] (2,0) -- (1,1);
		\draw[blue,thick] (4,0) -- (2,2);
		\filldraw[black] (3,3) circle (5pt);
		\filldraw[black] (2,2) circle (5pt);
		\filldraw[black] (1,1) circle (5pt);
		\or % C_5: (1 vee (2 vee 3)) vee 4
		\draw[blue,thick] (0,0) -- (3,3) -- (6,0);
		\draw[blue,thick] (2,0) -- (3,1);
		\draw[blue,thick] (4,0) -- (2,2);
		\filldraw[black] (3,3) circle (5pt);
		\filldraw[black] (2,2) circle (5pt);
		\filldraw[black] (3,1) circle (5pt);
		\fi
		\foreach \x/\j in {0/1,2/2,4/3,6/4}{%
			\filldraw[black] (\x,0) circle (5pt);
			\node at (\x,-1) {$\j$};
		}
		\node at (8,0) {$\cdots$};
		#2
		\node at (3,-2.3) {$#3$};
	\end{tikzpicture}%
}
\newcommand{\FPNonfreeForest}[2]{%
	\begin{tikzpicture}[scale=0.24,baseline=(current bounding box.center)]
		\path[use as bounding box] (-1,-2.7) rectangle (11,4.6);
		\draw[black,thin] (-0.8,0) -- (9,0);
		\ifcase#1\relax
		\or % F
		\draw[blue,thick] (1,1) -- (0,0);
		\draw[blue,thick] (1,1) -- (2,0);
		\filldraw[black] (1,1) circle (5pt);
		\draw[blue,thick] (5,1) -- (4,0);
		\draw[blue,thick] (5,1) -- (6,0);
		\filldraw[black] (5,1) circle (5pt);
		\or % D_1
		\draw[blue,thick] (4,4) -- (3,3);
		\draw[blue,thick] (3,3) -- (1,1);
		\draw[blue,thick] (1,1) -- (0,0);
		\draw[blue,thick] (1,1) -- (2,0);
		\filldraw[black] (1,1) circle (5pt);
		\draw[blue,thick] (3,3) -- (5,1);
		\draw[blue,thick] (5,1) -- (4,0);
		\draw[blue,thick] (5,1) -- (6,0);
		\filldraw[black] (5,1) circle (5pt);
		\filldraw[black] (3,3) circle (5pt);
		\draw[blue,thick] (4,4) -- (8,0);
		\filldraw[black] (4,4) circle (5pt);
		\or % D_2
		\draw[blue,thick] (4,4) -- (3,3);
		\draw[blue,thick] (3,3) -- (0,0);
		\draw[blue,thick] (3,3) -- (4,2);
		\draw[blue,thick] (4,2) -- (2,0);
		\draw[blue,thick] (4,2) -- (5,1);
		\draw[blue,thick] (5,1) -- (4,0);
		\draw[blue,thick] (5,1) -- (6,0);
		\filldraw[black] (5,1) circle (5pt);
		\filldraw[black] (4,2) circle (5pt);
		\filldraw[black] (3,3) circle (5pt);
		\draw[blue,thick] (4,4) -- (8,0);
		\filldraw[black] (4,4) circle (5pt);
		\or % D_3
		\draw[blue,thick] (4,4) -- (3,3);
		\draw[blue,thick] (3,3) -- (0,0);
		\draw[blue,thick] (3,3) -- (4,2);
		\draw[blue,thick] (4,2) -- (3,1);
		\draw[blue,thick] (3,1) -- (2,0);
		\draw[blue,thick] (3,1) -- (4,0);
		\filldraw[black] (3,1) circle (5pt);
		\draw[blue,thick] (4,2) -- (6,0);
		\filldraw[black] (4,2) circle (5pt);
		\filldraw[black] (3,3) circle (5pt);
		\draw[blue,thick] (4,4) -- (8,0);
		\filldraw[black] (4,4) circle (5pt);
		\or % D_4
		\draw[blue,thick] (4,4) -- (3,3);
		\draw[blue,thick] (3,3) -- (2,2);
		\draw[blue,thick] (2,2) -- (1,1);
		\draw[blue,thick] (1,1) -- (0,0);
		\draw[blue,thick] (1,1) -- (2,0);
		\filldraw[black] (1,1) circle (5pt);
		\draw[blue,thick] (2,2) -- (4,0);
		\filldraw[black] (2,2) circle (5pt);
		\draw[blue,thick] (3,3) -- (6,0);
		\filldraw[black] (3,3) circle (5pt);
		\draw[blue,thick] (4,4) -- (8,0);
		\filldraw[black] (4,4) circle (5pt);
		\or % D_5
		\draw[blue,thick] (4,4) -- (3,3);
		\draw[blue,thick] (3,3) -- (2,2);
		\draw[blue,thick] (2,2) -- (0,0);
		\draw[blue,thick] (2,2) -- (3,1);
		\draw[blue,thick] (3,1) -- (2,0);
		\draw[blue,thick] (3,1) -- (4,0);
		\filldraw[black] (3,1) circle (5pt);
		\filldraw[black] (2,2) circle (5pt);
		\draw[blue,thick] (3,3) -- (6,0);
		\filldraw[black] (3,3) circle (5pt);
		\draw[blue,thick] (4,4) -- (8,0);
		\filldraw[black] (4,4) circle (5pt);
		\or % D_6
		\draw[blue,thick] (4,4) -- (0,0);
		\draw[blue,thick] (4,4) -- (5,3);
		\draw[blue,thick] (5,3) -- (4,2);
		\draw[blue,thick] (4,2) -- (2,0);
		\draw[blue,thick] (4,2) -- (5,1);
		\draw[blue,thick] (5,1) -- (4,0);
		\draw[blue,thick] (5,1) -- (6,0);
		\filldraw[black] (5,1) circle (5pt);
		\filldraw[black] (4,2) circle (5pt);
		\draw[blue,thick] (5,3) -- (8,0);
		\filldraw[black] (5,3) circle (5pt);
		\filldraw[black] (4,4) circle (5pt);
		\or % D_7
		\draw[blue,thick] (4,4) -- (0,0);
		\draw[blue,thick] (4,4) -- (5,3);
		\draw[blue,thick] (5,3) -- (4,2);
		\draw[blue,thick] (4,2) -- (3,1);
		\draw[blue,thick] (3,1) -- (2,0);
		\draw[blue,thick] (3,1) -- (4,0);
		\filldraw[black] (3,1) circle (5pt);
		\draw[blue,thick] (4,2) -- (6,0);
		\filldraw[black] (4,2) circle (5pt);
		\draw[blue,thick] (5,3) -- (8,0);
		\filldraw[black] (5,3) circle (5pt);
		\filldraw[black] (4,4) circle (5pt);
		\fi
		\foreach \x/\j in {0/1,2/2,4/3,6/4,8/5}{%
			\filldraw[black] (\x,0) circle (5pt);
			\node at (\x,-1) {$\j$};
		}
		\node at (10,0) {$\cdots$};
		\node at (4,-2.3) {$#2$};
	\end{tikzpicture}%
}

\title[Monk and Pieri rules for forest polynomials]{Monk and Pieri rules for forest polynomials}

\author{Haojun Bai}
\address{Center for Combinatorics, Nankai University, LPMC, Tianjin 300071, P.R. China
		}
\email{haojunbai@mail.nankai.edu.cn}

\author{Shu Xiao Li}
\thanks{Shu Xiao Li was supported in part by the Fundamental Research Funds for the Central Universities.}
\address{
	Research Center for Mathematics and Interdisciplinary Sciences,
	Shandong University \& Frontiers Science Center for Nonlinear Expectations, Ministry of Education,
	Qingdao, Shandong 266237, P.R. China}
\email{lishuxiao@sdu.edu.cn}

\author[D.G.L.~Wang]{David G.L. Wang}
\thanks{David G.L. Wang was supported by the General Program of the
	National Natural Science Foundation of China (Grant No.~12171034).}
\address[David G.L. Wang]{School of Mathematics and Statistics
	\& Beijing Key Laboratory on MCAACI,
	Beijing Institute of Technology;
	MIIT Key Laboratory of Mathematical Theory and Computation in Information Security,
	Beijing 102400,
	P. R. China}
\email{glw@bit.edu.cn}

\subjclass[2010]{05E05}

\keywords{}

\begin{document}
	\begin{abstract}
		This article gives Pieri rules for forest polynomials and Monk rule for double forest polynomials. We construct binary relations $\lessdot_i$ and $\prec_i$ based on the local Tamari moves on indexed forests. These relations are explicit and constructive. For ordinary forest polynomials, the $h$- and $e$-Pieri rules are given by all chains of a fixed length along the cover relations, with certain restriction on the left or right endpoints. For double forest polynomials, our Monk rule shows every off-diagonal degree-one structure coefficient is always a single positive root, whose weight is also controlled by Tamari moves. Additionally, we introduce an insertion algorithm on positive injective words that yields a Monk rule for the corresponding slide polynomials. Regrouping this rule along Sylvester classes gives an independent Monk rule for ordinary forest polynomials.
	\end{abstract}
	\maketitle
	
\section{Introduction}

The multiplication of distinguished polynomial bases is a meeting point of algebraic geometry, equivariant topology, and enumerative combinatorics. In classical type-$A$ Schubert calculus, the Schubert polynomials of Lascoux and Sch\"utzenberger~\cite{LS} represent cohomology classes of flag varieties and make their intersection products accessible to explicit computation. Monk's formula~\cite{Monk} describes multiplication by the Schubert polynomial $\fS_{s_i}(\bx)=x_1+\cdots+x_i$ through suitable Bruhat covers, each occurring with coefficient one. The Pieri formulas for complete homogeneous and elementary symmetric polynomials~\cite{Sottile} extend this principle to higher degrees. The relationship between structure constants and chains in Bruhat order, developed further by Bergeron and Sottile~\cite{BS}, suggests a general objective: to express polynomial multiplication by local combinatorial operations whose iterations enumerate the coefficients.

Equivariant cohomology adds more data to this picture. The localization framework of Goresky, Kottwitz, and MacPherson~\cite{GKM} reconstructs equivariant cohomology from fixed points and one-dimensional torus orbits under suitable hypotheses. In the flag-variety setting, Billey's restriction formulas~\cite{Billey} connect localization values with Bruhat order, while Graham's positivity theorem~\cite{Graham} places equivariant Schubert structure constants among nonnegative integer combinations of monomials in simple roots.

The quasisymmetric side has its origins in the enumerative theory of labelled posets in Gessel's work on $P$-partitions and quasisymmetric functions~\cite{Gessel} The fundamental slide polynomials of Assaf and Searles~\cite{AS} lift fundamental quasisymmetric functions to the full polynomial ring and admit positive product rules. Compatible words and slide
expansions are therefore natural tools for turning polynomial identities into insertion and enumeration problems.

A second motivation comes from coinvariant algebras. Replacing symmetric polynomials by quasisymmetric polynomials gives the quotient
\[
\mathrm{QSCoinv}_n
=\Z[x_1,\ldots,x_n]\big/
\big\langle\mathrm{QSym}_n^+\big\rangle,
\]
where $\mathrm{QSym}_n^+$ denotes the positive-degree part of the ring of quasisymmetric polynomials in $n$ variables. Aval, Bergeron, and Bergeron~\cite{ABB} showed that the corresponding quotient over $\mathbb{Q}$ has Catalan dimension and a basis indexed by Dyck paths. This places the search for a quasisymmetric Schubert calculus in a Catalan setting: binary trees and noncrossing partitions become natural indexing objects for the algebra.

Forest polynomials, $\fP_F$, supply a basis adapted to this setting. Nadeau and Tewari~\cite{NT} introduced the polynomials $\fP_F(\bx)$, indexed by indexed forests, and showed that they form a basis of $\Z[x_1,x_2,\ldots]$. Nadeau, Spink, and Tewari~\cite{NST} then established the divided-difference formalism for forest polynomials. Their subsequent work~\cite{NSTExpansions} places Schubert, forest, and slide expansions in a common operator framework. This geometric origin also connects the theory to mixed Eulerian numbers and the polyhedral setting of Postnikov~\cite{Postnikov}.

The geometric theory develops this analogy at the level of cycles and torus actions. Nadeau, Spink, and Tewari~\cite{NSTGeometry} constructed quasisymmetric Schubert cycles using toric Richardson varieties. Bergeron, Gagnon, Nadeau, Spink, and Tewari~\cite{BGNSTQF} constructed the quasisymmetric flag variety, a toric complex whose ordinary cohomology is the quasisymmetric coinvariant ring and whose torus-fixed points are indexed by noncrossing partitions. Their double forest polynomials $\fP_F(\bx;\bt)$ and equivariant trimming operators~\cite{BGNST} provide the
corresponding polynomial calculus. The newer grove polynomials of Nadeau, Spink, and Tewari~\cite{NSTGrove} extend this program to $K$-theory through set-valued fillings of forests, illustrating the scope of the surrounding geometric theory.

Several existing multiplication results are particularly relevant to the present work. Nadeau, Spink, and Tewari~\cite[Theorem~7.4]{NST} proved that $(x_1+\cdots+x_i)\fP_F(\bx)$ is multiplicity-free and forest-positive. They observed that higher Pieri products can have multiplicities and asked for direct enumerative descriptions of coefficients \cite[Remark~7.5]{NST}. More recently, Bansal~\cite{Bansal} studied the forest order generated by ordinary Monk supports, including its chain structure and intervals, and gave a word-insertion description of the degree-one product.

%Binary trees also bring a distinct algebraic structure to these questions. Local rotations generate the Tamari order and describe the edges of the associahedron; an explicit polyhedral realization is due to Loday~\cite{Loday}. Planar binary trees index the Hopf algebra of Loday and Ronco~\cite{LR}, and the Sylvester monoid of Hivert, Novelli, and Thibon~\cite{HNT} organizes words through binary-search-tree insertion. These constructions provide a broader context for the two descriptions used here: constrained rotations act directly on forests, while insertion organizes compatible words into Sylvester classes. The rotations and the word conventions relevant to our formulas are specified in Sections \ref{sec:4} and \ref{sec:m1}.

%The relation between multiplication and labelled chains has a similarly broad setting. Bergeron, Mykytiuk, Sottile, and van Willigenburg~\cite{BMSW} developed noncommutative Pieri operators on graded posets, linking such operators to quasisymmetric generating functions and Hopf-algebra constructions. In our setting, the left and right endpoints of intervals attached to forest moves provide the labels for two descriptions of each Pieri coefficient. This makes it possible to pass from degree-one multiplication to higher-degree products through monotonicity conditions on chains.

The main contribution of this paper is to make the degree-one coefficients explicit through local tree operations, the Tamari moves, including their equivariant weights, and to express ordinary Pieri coefficients through interval-labelled chains of these operations. Our first main result is introducing the cover relations $\lessdot_i$ on indexed forests that models the Monk and Pieri rule. Their combinatorial definitions can be found in Section \ref{sec:4}.

\begin{theorem}[Monk rule, Theorem \ref{thm:m}]
	For indexed forest $F$, we have
	$$\fP_F(\bx)(x_1+\cdots+x_i)=\sum_{F\lessdot_i G}\fP_G(\bx).$$
\end{theorem}

Chains in the $\lessdot_i$ relation lead to the $h$- and $e$-Pieri rules as follows. In the formula below, $[a_k,b_k)$ are certain the interval of the cover relation $F_{k-1}\lessdot_i F_k$, $1\leq k\leq r$.

\begin{theorem}[Pieri rule, Theorem \ref{thm:pieri}]
	For indexed forest $F$, we have
	$$\fP_F(\bx)h_r(x_1,\dots,x_i)=\sum_{F=F_0\lessdot_i F_1\lessdot_i\cdots\lessdot_i F_r=G \atop a_1\leq\cdots\leq a_r}\fP_G(\bx)=\sum_{F=F_0\lessdot_i F_1\lessdot_i\cdots\lessdot_i F_r=G \atop b_1<\cdots< b_r}\fP_G(\bx),$$
	$$\fP_F(\bx)e_r(x_1,\dots,x_i)=\sum_{F=F_0\lessdot_i F_1\lessdot_i\cdots\lessdot_i F_r=G \atop a_1>\cdots> a_r}\fP_G(\bx)=\sum_{F=F_0\lessdot_i F_1\lessdot_i\cdots\lessdot_i F_r=G \atop b_1\geq\cdots\geq b_r}\fP_G(\bx).$$
\end{theorem}

In the other words, the Pieri rule counts the number of certain chains in the $\lessdot_i$ relation, and unlike the Schubert counterpart, the coefficients could be more than $1$.

In the equivariant setting, we give an answer to the equivariant Monk problem posed in \cite[Question~1.2(4)]{BGNST}: every nonzero off-diagonal coefficient is either one or a single positive root, with precise conditions for its occurrence. More precisely, we introduce a binary relation $\prec_i$ on indexed forests with the same number of internal nodes. These relation describe the support of the same-size, off-diagonal terms in the Monk rule. The coefficients $\mathrm{wt}_{F,G}$ can also be obtained directly from the $\prec_i$ relation.

\begin{theorem}[Equivariant Monk rule, Theorem \ref{thm:emonk}]
	For indexed forest $F$, we have the following formula. We emphasize that each of the coefficients $\mathrm{wt}_{F,G}$ is always a single positive root.
	$$\fP_F(\bx;\bt)\fS_{s_i}(\bx;\bt)=\left(\sum_{k=1}^i(t_{w_F(k)}-t_k)\fP_F(\bx;\bt)\right)+\left(\sum_{F\prec_i G}\mathrm{wt}_{F,G}\fP_G(\bx;\bt)\right)+\left(\sum_{F\lessdot_i G}\fP_G(\bx;\bt)\right).$$
\end{theorem}

Our starting point is the localization theorem of the equivariant cohomology of the quasisymmetric flag variety developed in \cite{BGNSTQF}. The equivariant Monk coefficients can be divided into same-size off-diagonal ones and ordinary Monk coefficients. We then apply the equivariant operators, and perform significant reductions using $\star$-Thompson monoid relations. Finally, both same-size off-diagonal and ordinary Monk coefficients can be indexed by admissible sequences, that will be introduced in section \ref{sec:3}. These results get re-interpreted using the Tamari moves on indexed forests in Section \ref{sec:4}, and the Monk and Pieri formulations above are obtained.

In Section \ref{sec:m1}, we describe our second approach to the Monk multiplication, using the fundamental slide polynomials $\mathfrak{F}_w$ of Assaf and Searles~\cite{AS}. We introduce an insertion algorithm for the Sylvester words of indexed forests, that leads to a suitable Monk rule for slide polynomials. Our construction ensures that the resulting slide polynomials can be grouped back into forest polynomials. This yields another formulation fo the Monk rule for forest polynomials. This approach is independent of the previous one.

\begin{theorem}[Monk rule, Theorem \ref{thm:monk}]
	We have $(x_1+\cdots+x_i)\fP_F=\sum_{G\in\Omega_i(F)}\fP_G$. In particular, all coefficients are multiplicity-free.
\end{theorem}

This paper is organized as follows. Section \ref{sec:2} reviews indexed forests, ordinary and double forest polynomials, and the operators used in their study. Section \ref{sec:3} determines the
equivariant Monk coefficients through localization and admissible sequences. Section \ref{sec:4} translates these formulas into Tamari moves and establishes the ordinary Monk and Pieri rules and the equivariant Monk rule. Section \ref{sec:m1} develops the insertion algorithm, proves the slide-polynomial multiplication formula, and derives another Monk rule for the ordinary forest polynomials.

\section{Preliminaries}\label{sec:2}

This section provides notation and terminology of indexed forests, forest polynomials, double forest polynomials and their established properties that will be used in this article. For more details, we refer the readers to \cite{BGNST,NT}.

\subsection{Indexed forests}

\begin{definition}
	A \textbf{(planar complete) binary tree} is a tree such that every vertex has either $0$ or $2$ children. Vertices with $2$ children are \textbf{internal nodes}, and other vertices are \textbf{leaves}. The \textbf{degree} of a binary tree is the number of internal nodes. The left and right children of $v$ are denoted by $v_L$ and $v_R$, respectively.
\end{definition}

\begin{definition}
	An \textbf{indexed forest} is an infinite sequence $(T_1,T_2,\dots)$ of binary trees such that all but finitely many of them are trivial trees. The leaves of the trees are indexed from left to right by $\mathbb{N}$. The set of indexed forests is denoted by $\For_\infty$. Each binary tree $T_i$ in the sequence of is called a \textbf{component} of the indexed forest. The \textbf{support} of a component $T$ is the set of indexes of its leaves, the \textbf{support} of $F$ is the union of supports of $T_1,\dots,T_\ell$ where $T_\ell$ is the last non-trivial component.
\end{definition}

By abuse of notation, we write $v\in F$ for $v$ is an internal node in $F$.

\begin{definition}
	Given an indexed forest $F$, we define the following labelings on internal nodes. For $v\in F$,
	\begin{enumerate}
		\item $\rho(v)$ is the index of the leftmost leaf of the subtree rooted at $v$.
		\item $\gamma(v)$ is the index of the rightmost leaf of the subtree rooted at $v_L$.
		\item $\eta(v)$ is the number of internal nodes in the subtree rooted at $v$.
		\item $(\rho+\eta)(v)$ is the index of the rightmost leaf of the subtree rooted at $v$.
	\end{enumerate}
	%The set $\{\gamma(v):v\in F\}$ is denoted by $L(F)$.
\end{definition}

\begin{definition}
	A node $v\in F$ is \textbf{terminal} if both children of $v$ are leaves. The set of terminal nodes is denoted by $\Ter(F)$. The \textbf{left terminal} of $F$ is given by $\LTer(F)=\{\rho(v):v\in \Ter(F)\}$.
\end{definition}

\begin{example}\label{ex:For}
	A binary tree of degree $4$ and an indexed forest with $\rho$-labeled internal nodes.
	\begin{center}
		\begin{tikzpicture}
			\node at (0,0) {\begin{tikzpicture}[scale=0.25,baseline=0pt]
					\draw[blue, thick] (0,4) -- (-4,0);
					\draw[blue, thick] (0,4) -- (4,0);
					\draw[blue, thick] (-3,1) -- (-2,0);
					\draw[blue, thick] (2,2) -- (0,0);
					\draw[blue, thick] (1,1) -- (2,0);
					\filldraw[black] (0,4) circle (5pt)  {};
					\filldraw[black] (2,2) circle (5pt)  {};
					\filldraw[black] (1,1) circle (5pt)  {};
					\filldraw[black] (-3,1) circle (5pt)  {};
					
					\node at (0,-1.2) {};
			\end{tikzpicture}};
			
			\node at (8,0) {\begin{tikzpicture}[scale=0.25,baseline=0pt]
					\draw[blue, thick] (-7,1) -- (-8,0);
					\draw[blue, thick] (-7,1) -- (-6,0);
					\filldraw[black] (-7,1) circle (5pt)  {};
					
					\node at (-7,2) {\tiny 2};
					
					\draw[blue, thick] (0,4) -- (-4,0);
					\draw[blue, thick] (0,4) -- (4,0);
					\draw[blue, thick] (-3,1) -- (-2,0);
					\draw[blue, thick] (2,2) -- (0,0);
					\draw[blue, thick] (1,1) -- (2,0);
					\filldraw[black] (0,4) circle (5pt)  {};
					\filldraw[black] (2,2) circle (5pt)  {};
					\filldraw[black] (1,1) circle (5pt)  {};
					\filldraw[black] (-3,1) circle (5pt)  {};
					
					\node at (-3,2) {\tiny 4};
					\node at (0,5) {\tiny 4};
					\node at (1,2) {\tiny 6};
					\node at (2,3) {\tiny 6};
					
					\filldraw[black] (-10,0) circle (5pt)  {};
					\filldraw[black] (-8,0) circle (5pt)  {};
					\filldraw[black] (-6,0) circle (5pt)  {};
					\filldraw[black] (-4,0) circle (5pt)  {};
					\filldraw[black] (-2,0) circle (5pt)  {};
					\filldraw[black] (0,0) circle (5pt)  {};
					\filldraw[black] (2,0) circle (5pt)  {};
					\filldraw[black] (4,0) circle (5pt)  {};
					\filldraw[black] (6,0) circle (5pt)  {};
					\filldraw[black] (8,0) circle (5pt)  {};
					\filldraw[black] (10,0) circle (5pt)  {};
					\filldraw[black] (12,0) circle (5pt)  {};
					\filldraw[black] (14,0) circle (5pt)  {};
					\filldraw[black] (16,0) circle (5pt)  {};
					\filldraw[black] (18,0) circle (5pt)  {};
					
					\draw[blue, thick] (10,2) -- (8,0);
					\draw[blue, thick] (10,2) -- (12,0);
					\draw[blue, thick] (11,1) -- (10,0);
					\filldraw[black] (10,2) circle (5pt)  {};
					\filldraw[black] (11,1) circle (5pt)  {};
					
					\node at (10,3) {\tiny 10};
					\node at (11.2,2) {\tiny 11};
					
					\node at (-10,-1.2) {1};
					\node at (-8,-1.2) {2};
					\node at (-6,-1.2) {3};
					\node at (-4,-1.2) {4};
					\node at (-2,-1.2) {5};
					\node at (0,-1.2) {6};
					\node at (2,-1.2) {7};
					\node at (4,-1.2) {8};
					\node at (6,-1.2) {9};
					\node at (8,-1.2) {10};
					\node at (10,-1.2) {11};
					\node at (12,-1.2) {12};
					\node at (14,-1.2) {13};
					\node at (16,-1.2) {14};
					\node at (18,-1.2) {15};
					
					\draw[black] (-12,0) -- (20,0);
			\end{tikzpicture}};
		\end{tikzpicture}
	\end{center}
	Its left terminal is $\LTer(F)=\{2,4,6,11\}$.
\end{example}

\begin{definition}
	Let $F$ be an indexed forest, and $i\in\mathbb{N}$, the \textbf{blossomed forest}, $F\cdot i$, is obtained from $F$ by giving the leaf at $i$ two children. If $i\in\LTer(F)$, the \textbf{trimmed forest}, $F/i$, is obtained from $F$ by removing the two children of the terminal node with label $i$. 
\end{definition}

\begin{example}
	Suppose $F$ is the indexed forest given in Example \ref{ex:For}, then, $F\cdot 5$ is
	\begin{center}
		\begin{tikzpicture}
			\node at (0,0) {\begin{tikzpicture}[scale=0.25,baseline=0pt]
					\draw[blue, thick] (-7,1) -- (-8,0);
					\draw[blue, thick] (-7,1) -- (-6,0);
					\filldraw[black] (-7,1) circle (5pt)  {};
					
					\draw[blue, thick] (1,5) -- (-4,0);
					\draw[blue, thick] (1,5) -- (6,0);
					\draw[blue, thick] (-2,2) -- (0,0);
					\draw[blue, thick] (-1,1) -- (-2,0);
					\draw[blue, thick] (4,2) -- (2,0);
					\draw[blue, thick] (3,1) -- (4,0);
					\filldraw[black] (1,5) circle (5pt)  {};
					\filldraw[black] (4,2) circle (5pt)  {};
					\filldraw[black] (3,1) circle (5pt)  {};
					\filldraw[black] (-2,2) circle (5pt)  {};
					\filldraw[black] (-1,1) circle (5pt)  {};
					
					\filldraw[black] (-10,0) circle (5pt)  {};
					\filldraw[black] (-8,0) circle (5pt)  {};
					\filldraw[black] (-6,0) circle (5pt)  {};
					\filldraw[black] (-4,0) circle (5pt)  {};
					\filldraw[black] (-2,0) circle (5pt)  {};
					\filldraw[black] (0,0) circle (5pt)  {};
					\filldraw[black] (2,0) circle (5pt)  {};
					\filldraw[black] (4,0) circle (5pt)  {};
					\filldraw[black] (6,0) circle (5pt)  {};
					\filldraw[black] (8,0) circle (5pt)  {};
					\filldraw[black] (10,0) circle (5pt)  {};
					\filldraw[black] (12,0) circle (5pt)  {};
					\filldraw[black] (14,0) circle (5pt)  {};
					\filldraw[black] (16,0) circle (5pt)  {};
					\filldraw[black] (18,0) circle (5pt)  {};
					
					\draw[blue, thick] (12,2) -- (10,0);
					\draw[blue, thick] (12,2) -- (14,0);
					\draw[blue, thick] (13,1) -- (12,0);
					\filldraw[black] (12,2) circle (5pt)  {};
					\filldraw[black] (13,1) circle (5pt)  {};
					
					\node at (-10,-1.2) {1};
					\node at (-8,-1.2) {2};
					\node at (-6,-1.2) {3};
					\node at (-4,-1.2) {4};
					\node at (-2,-1.2) {5};
					\node at (0,-1.2) {6};
					\node at (2,-1.2) {7};
					\node at (4,-1.2) {8};
					\node at (6,-1.2) {9};
					\node at (8,-1.2) {10};
					\node at (10,-1.2) {11};
					\node at (12,-1.2) {12};
					\node at (14,-1.2) {13};
					\node at (16,-1.2) {14};
					\node at (18,-1.2) {15};
					
					\draw[black] (-12,0) -- (20,0);
			\end{tikzpicture}};
		\end{tikzpicture}
	\end{center}and $F/6$ is shown below.
	\begin{center}
		\begin{tikzpicture}
			\node at (0,0) {\begin{tikzpicture}[scale=0.25,baseline=0pt]
					\draw[blue, thick] (-7,1) -- (-8,0);
					\draw[blue, thick] (-7,1) -- (-6,0);
					\filldraw[black] (-7,1) circle (5pt)  {};
					
					\draw[blue, thick] (-1,3) -- (-4,0);
					\draw[blue, thick] (-1,3) -- (2,0);
					\draw[blue, thick] (-3,1) -- (-2,0);
					\draw[blue, thick] (1,1) -- (0,0);
					\filldraw[black] (-1,3) circle (5pt)  {};
					\filldraw[black] (1,1) circle (5pt)  {};
					\filldraw[black] (-3,1) circle (5pt)  {};
					
					\filldraw[black] (-10,0) circle (5pt)  {};
					\filldraw[black] (-8,0) circle (5pt)  {};
					\filldraw[black] (-6,0) circle (5pt)  {};
					\filldraw[black] (-4,0) circle (5pt)  {};
					\filldraw[black] (-2,0) circle (5pt)  {};
					\filldraw[black] (0,0) circle (5pt)  {};
					\filldraw[black] (2,0) circle (5pt)  {};
					\filldraw[black] (4,0) circle (5pt)  {};
					\filldraw[black] (6,0) circle (5pt)  {};
					\filldraw[black] (8,0) circle (5pt)  {};
					\filldraw[black] (10,0) circle (5pt)  {};
					\filldraw[black] (12,0) circle (5pt)  {};
					\filldraw[black] (14,0) circle (5pt)  {};
					\filldraw[black] (16,0) circle (5pt)  {};
					\filldraw[black] (18,0) circle (5pt)  {};
					
					\draw[blue, thick] (8,2) -- (6,0);
					\draw[blue, thick] (8,2) -- (10,0);
					\draw[blue, thick] (9,1) -- (8,0);
					\filldraw[black] (8,2) circle (5pt)  {};
					\filldraw[black] (9,1) circle (5pt)  {};
					
					\node at (-10,-1.2) {1};
					\node at (-8,-1.2) {2};
					\node at (-6,-1.2) {3};
					\node at (-4,-1.2) {4};
					\node at (-2,-1.2) {5};
					\node at (0,-1.2) {6};
					\node at (2,-1.2) {7};
					\node at (4,-1.2) {8};
					\node at (6,-1.2) {9};
					\node at (8,-1.2) {10};
					\node at (10,-1.2) {11};
					\node at (12,-1.2) {12};
					\node at (14,-1.2) {13};
					\node at (16,-1.2) {14};
					\node at (18,-1.2) {15};
					
					\draw[black] (-12,0) -- (20,0);
			\end{tikzpicture}};
		\end{tikzpicture}
	\end{center}
\end{example}

\begin{definition}
	The \textbf{trimming sequence} of a degree $k$ indexed forest $F$ is defined recursively as
	$$\Trim(F)=\{(i_1,\dots,i_k):i_k\in\LTer(F),(i_1,\dots,i_{k-1})\in\Trim(F/i_k)\}$$
	with initial condition $\Trim(\varnothing)=\{\varnothing\}$.
\end{definition}

\begin{example}\label{ex:F2}
	Suppose $F$ is the indexed forest
	\begin{center}
		\begin{tikzpicture}
			\node at (0,0) {\begin{tikzpicture}[scale=0.25,baseline=0pt]
					\draw[blue, thick] (-7,1) -- (-8,0);
					\draw[blue, thick] (-7,1) -- (-6,0);
					\filldraw[black] (-7,1) circle (5pt)  {};
					
					\draw[blue, thick] (0,4) -- (-4,0);
					\draw[blue, thick] (0,4) -- (4,0);
					\draw[blue, thick] (-3,1) -- (-2,0);
					\draw[blue, thick] (2,2) -- (0,0);
					\draw[blue, thick] (1,1) -- (2,0);
					\filldraw[black] (0,4) circle (5pt)  {};
					\filldraw[black] (2,2) circle (5pt)  {};
					\filldraw[black] (1,1) circle (5pt)  {};
					\filldraw[black] (-3,1) circle (5pt)  {};
					
					\filldraw[black] (-10,0) circle (5pt)  {};
					\filldraw[black] (-8,0) circle (5pt)  {};
					\filldraw[black] (-6,0) circle (5pt)  {};
					\filldraw[black] (-4,0) circle (5pt)  {};
					\filldraw[black] (-2,0) circle (5pt)  {};
					\filldraw[black] (0,0) circle (5pt)  {};
					\filldraw[black] (2,0) circle (5pt)  {};
					\filldraw[black] (4,0) circle (5pt)  {};
					\filldraw[black] (6,0) circle (5pt)  {};
					\filldraw[black] (8,0) circle (5pt)  {};
					\filldraw[black] (10,0) circle (5pt)  {};
					
					\node at (-10,-1.2) {1};
					\node at (-8,-1.2) {2};
					\node at (-6,-1.2) {3};
					\node at (-4,-1.2) {4};
					\node at (-2,-1.2) {5};
					\node at (0,-1.2) {6};
					\node at (2,-1.2) {7};
					\node at (4,-1.2) {8};
					\node at (6,-1.2) {9};
					\node at (8,-1.2) {10};
					\node at (10,-1.2) {11};
					
					\draw[black] (-12,0) -- (12,0);
			\end{tikzpicture}};
		\end{tikzpicture}
	\end{center}
	Then, $\Trim(F)$ is
	\begin{center}
		\begin{tabular}{ccccc}
			$(2,4,4,6,6)$ & $(3,2,4,6,6)$ & $(2,4,5,4,6)$ & $(3,3,2,6,6)$ & $(3,2,5,4,6)$ \\
			$(2,4,5,5,4)$ & $(3,3,5,2,6)$ & $(3,4,2,4,6)$ & $(3,2,5,5,4)$ & $(3,4,3,2,6)$ \\
			$(3,3,5,5,2)$ & $(3,4,2,5,4)$ & $(3,4,3,5,2)$ & $(3,4,4,2,4)$ & $(3,4,4,3,2)$
		\end{tabular}
	\end{center}
\end{example}

\begin{definition}\label{def:trim}
	The \textbf{Bergeron-Sottile maps} are $R_i:\Z[x_1,x_2,\dots]\to\Z[x_1,x_2,\dots]$ with
	$$R_if(x_1,x_2,\dots)=f(x_1,\dots,x_{i-1},0,x_i,x_{i+1},\dots)$$
	and the \textbf{trimming operators} are $T_i:\Z[x_1,x_2,\dots]\to\Z[x_1,x_2,\dots]$ with
	$$T_if=\frac{R_{i+1}f-R_if}{x_i}.$$
	For an indexed forest $F$, we set $T_F=T_{i_1}\cdots T_{i_n}$ where $(i_1,\dots,i_n)\in\Trim(F)$. This is independent of the choice of the trimming sequence.
\end{definition}

\begin{definition}\label{thm:F}
	The polynomials $\{\fP_F:F\in\For_\infty\}$ are determined by
	\begin{itemize}
		\item $\fP_F$ are homogeneous;
		\item $\fP_\varnothing=1$;
		\item $T_i\fP_F=\begin{cases}
			\fP_{F/i} & i\in\LTer(F) \\
			0 & i\notin\LTer(F)
		\end{cases}$.
	\end{itemize}
	These polynomials $\fP_F$ are the \textbf{forest polynomial}.
\end{definition}

\begin{definition}\label{def:order}
	The set of internal nodes of $F$ is totally ordered by the lexicographic order on the pair $(\rho(v),-\eta(v))$. The resulting ordered list is given by $\{v\in F\}=\{v_1^F\leq v_2^F\leq\cdots\leq v_n^F\}$. In particular, $v_1^F$ is the leftmost root in $F$.
\end{definition}

\begin{definition}
	A \textbf{Sylvester word} of $F$ is an injective word $w_1\cdots w_{\ell}$ satisfying $\{w_k:1\leq k\leq\ell\}=\{\gamma(v):v\in F\}$, and $u$ is a node in the subtree of $F$ rooted at $v$ implies $i\leq j$ where $w_i=\gamma(v)$ and $w_j=\gamma(u)$. The set of Sylvester words of $F$ is denoted by $\mathrm{Syl}(F)$.
\end{definition}

\begin{definition}
	The set of indexed forests is partially ordered by $F\leq_R G$ if and only if there exists a trimming sequence $(j_1,\dots,j_n)$ of $G$ such that a suffix $(j_m,\dots,j_n)$ is a trimming sequence for $F$. In this case, we write $G/F$ to be the indexed forest with trimming sequence $(j_1,\dots,j_{m-1})$.
\end{definition}

\begin{definition}
	Given an indexed forest $F$, consider the subdiagram $\tilde{F}$ such that for all $v\in F$ such that $v=u_R$ for some $u\in F$, delete the edge $(v,v_L)$. We define a non-crossing partition $\mathrm{ForToNC}(F)$ in which $i,j$ are in the same block if and only if the two leaves $i,j$ are connected by the diagram $\tilde{F}$. We write $\pi_F=\mathrm{ForToNC}(F)$.
\end{definition}

\begin{definition}
	Given a non-crossing partition $\pi=B_1/\cdots/B_\ell$, we define a permutation $\mathrm{Perm}(\pi)=(B_1)\cdots (B_\ell)$ where $(B_i)$ is the cycle $(b_i^1\cdots b_i^k)$ and $B_i=\{b_i^1>\cdots>b_i^k\}$. We write $w_F=\mathrm{Perm}\circ\mathrm{ForToNC}(F)$. We denote $\mathrm{NC}_n=\{\pi_F:\mathrm{supp}(F)\subseteq [n]\}$ and $\NCPerm_n=\mathrm{Perm}(\mathrm{NC}_n)$. We write $F<G$ if $w_F<w_G$ in Bruhat order.
\end{definition}

\begin{definition}\label{def:h}
	For $F< G$, let $G$ have a trimming sequence $(j_1,\dots,j_n)$ with $j_1\leq\cdots\leq j_n$. Set $S_r$ be the forest with trimming sequence $(j_{r+1},j_{r+2},\dots,j_n)$. In particular, $S_0=G$ and $S_n=\varnothing$. Let $h(F,G)=\min\{r:S_r\leq_R F\}$, $h'(F,G)=\max\{r:j_r=j_{h(F,G)}\}$, and $v_F^G=v_{h(F,G)}^G$.
\end{definition}

\subsection{Equivariant operators}

\begin{definition}
	Given finite set $A\subseteq\mathbb{Z}_{>0}$, let $\bar{A}$ denote its complement $\mathbb{Z}_{>0}\setminus A$, and let $\bar{A}_i$ be the $i$-th smallest number in $\bar{A}$. Let $t_{i,A}=t_{\bar{A}_i}$ for $i\in\mathbb{Z}_{>0}$, and $\hat{\bt}_A=(t_{1,A},t_{2,A},\dots)$. Consider the following four operations.
	\begin{align*}
		\mathrm{ev}_Af(\bx;\bt)&=f(\hat{\bt}_A;\bt), \\
		R_{i,A}^-f(\bx;\bt)&=f(x_1,\dots,x_{i-1},t_{i,A},x_i,x_{i+1},\dots ;\bt), \\
		R_{i,A}^+f(\bx;\bt)&=f(x_1,\dots,x_{i-1},x_i,t_{i,A},x_{i+1},\dots ;\bt), \\
		E_{i,A}f(\bx;\bt)&=\frac{R_{i,A}^+f(\bx;\bt)-R_{i,A}^- f(\bx;\bt)}{x_i-t_{i,A}}.
	\end{align*}
\end{definition}

\begin{definition}
	\textbf{Double forest polynomials} $\fP_F(\bx;\bt)$ are the unique polynomials satisfying
	$$\fP_F(\bt;\bt)=\mathbf{1}_{F=\varnothing},$$
	$$E_{i,\varnothing}\fP_F(\bx;\bt)=\begin{cases}
		 \fP_{F/i}(\bx;\hat{\bt}_{\{i\}}) & i\in\LTer(F) \\ 
		 0 & i\notin \LTer(F)
	\end{cases}$$
	The double forest polynomials reduce to forest polynomials via $\fP_F(\bx;0)=\fP_F(\bx)$.
\end{definition}

\begin{definition}
	Given two finite sets $A,B\subseteq\mathbb{Z}_{>0}$, define $A\star B=\{\bar{B}_i:i\in A\}\cup B$. When $A=\{i\}$ is a singleton, we abbreviate $\{i\}\star B$ by $i\star B$.
\end{definition}

\begin{lemma}[\cite{BGNST} Lemma 9.2]\label{lem:BGNST}
	Let $F$ be an indexed forest with trimming sequence $(i_1,\dots,i_n)$. Then,
	$$\{\gamma(v):v\in F\}=i_1\star i_2\star\cdots\star i_n.$$
\end{lemma}

\begin{definition}
	Let $X_i$ be one of $\{R_i^-,R_i^+,E_i\}$, their \textbf{star product} is
	$$[X_{i_1}\star\cdots \star X_{i_n}]_A=X_{i_1,i_2\star\cdots\star i_n\star A}X_{i_2,i_3\star\cdots\star i_n\star A}\cdots X_{i_n, A}$$
	and 
	$$[\mathrm{ev}\star X_{i_1}\star\cdots \star X_{i_n}]_A=\mathrm{ev}_{i_1\star\cdots\star i_n\star A}[X_{i_1}\star\cdots\star X_{i_n}]_A.$$
	By convention, the subscript $A$ is omitted when $A=\varnothing$.
	
	Let $F$ be an indexed forest with trimming sequence $(i_1,\dots,i_n)$. For $A\subseteq\mathbb{Z}_>0$,
	\begin{align*}
		[E_F]_A&=E_{i_1,i_2\star\cdots\star i_n\star A}E_{i_2,i_3\star\cdots\star i_n\star A}\cdots E_{i_n,A},\\
		[\mathrm{ev}\star E_F]_A&=\mathrm{ev}_{i_1\star\cdots\star i_n\star A}[E_F]_A.\\
	\end{align*}
	$[E_F]_A$ and $[\mathrm{ev}\star E_F]_A$ are independent of the choices of trimming sequences of $F$.
\end{definition}

\begin{theorem}[\cite{BGNST} Theorem 10.10]\label{thm:BGNST2}
	For $f(\bx;\bt)\in\mathbb{Z}[\bx;\bt]$,
	$$f(\bx;\bt)=\sum_{F\in\For_\infty}\left([\mathrm{ev}\star E_F]f(\bx;\bt)\right)\fP_F(\bx;\bt).$$
\end{theorem}

\begin{theorem}[\cite{BGNST} Proposition 11.5]\label{thm:BGNST3}
	For $f(\bx;\bt),g(\bx;\bt)\in\mathbb{Z}[\bx;\bt]$ and $A\subseteq \mathbb{Z}_{>0}$,
	$$E_{i,A}(fg)=(E_{i,A}f)(R_{i,A}^-g)+(R_{i,A}^+f)(E_{i,A}g).$$
\end{theorem}

\begin{theorem}[\cite{BGNST} Proposition 10.9]\label{thm:BGNST4}
	For indexed forests $F,G$,
	$$[E_F]\fP_G(\bx;\bt)=\begin{cases}
		\fP_{G/F}(\bx;\hat{\bt}_{\{\gamma(v):v\in F\}}) & F\leq_R G\\
		0 & F\nleq_R G
	\end{cases}$$
\end{theorem}

\begin{theorem}[\cite{BGNST} Theorem 10.5]\label{thm:BGNST5}
	For $p<q$ and $A$, we have
	$$[E_q\star E_p]_A=[E_p\star E_{q+1}]_A,$$
	$$[R_q^+\star E_p]_A=[E_p\star R_{q+1}^+]_A,$$
	$$[E_q\star R_p^+]_A=[R_p^+\star E_{q+1}]_A.$$
\end{theorem}

\begin{theorem}[\cite{BGNST} Theorem 10.5 and Proposition 11.2]\label{thm:BGNST6}
	For $p$ and $A$, we have
	$$[E_p\star R_p^-]_A=[R_p^-\star E_{p+1}]_A,$$
	$$[E_p\star R_{p+1}^+]_A=[R_p^+\star E_p]_A,$$
	$$[E_p\star R_{p+1}^-]_A=[R_p^+\star E_p]_A+[R_p^-\star E_{p+1}]_A+(t_{p+1,A}-t_{p,A})[E_p\star E_{p+1}]_A,$$
	$$[E_p\star R_p^+]_A=[R_p^+\star E_p]_A+[R_p^-\star E_{p+1}]_A+(t_{p+1,A}-t_{p,A})[E_p\star E_p]_A.$$
\end{theorem}

\section{The admissible sequences}\label{sec:3}

This section provides an equivariant Monk rule for double forest polynomials using equivariant operators and admissible sequences.

\subsection{Localization}

This subsection determines the triangularity of the structure coefficients using the localization map.

\begin{theorem}[\cite{BGNST} Theorem 8.3]\label{thm:BGNST}
	The localization map $\Phi:H_T^*(\mathrm{QFL}_n)\to\bigoplus_{w\in\mathrm{NCPerm}_n}\mathbb{Z}[\mathbf{t}]/\sim$, where $\sim$ consists of the GKM divisibility relations, defined by
	$$\Phi=\bigoplus_{w\in\NCPerm_n}\Phi_w\text{ and }\Phi_w(f)=f(t_{w(1)},\dots,t_{w(n)};t_1,\dots,t_n)$$
	is a ring isomorphism. Moreover, $\{\Phi(\fP_F):\mathrm{supp}(F)\subseteq [n]\}$ forms the flow-up basis of $\bigoplus_{w\in\NCPerm_n}\mathbb{Z}[\mathbf{t}]/\sim$.
\end{theorem}

\begin{definition}
	For each $j\in\mathrm{supp}(F)$, define
	$$\mathrm{vert}_j(F)=\begin{cases}
		2|\{B\in\pi_F:\min B< j< \max B,j\notin B\}| & j\text{ is maximal in the block of }\pi_F\text{ containing }j \\
		2|\{B\in\pi_F:\min B< j< \max B,j\notin B\}|+1 & \text{otherwise}
	\end{cases}.$$
\end{definition}

\begin{theorem}[\cite{GW} Theorem 1.2]\label{thm:GW}
	We have $F\leq G$ if and only if $\mathrm{vert}_j(F)\leq\mathrm{vert}_j(G)$ for all $j$.
\end{theorem}

\begin{lemma}\label{lem:e1-1}
	We have $|F|=|\{j\in\mathrm{supp}(F):\mathrm{vert}_j(F)>0\}|$. In particular, $w_F\leq w_G$ implies $|F|\leq |G|$.
\end{lemma}

\begin{proof}
	For each $v\in F$, let $B$ be the block of $\pi_F$ containing the rightmost leaf of the sub-tree rooted at $v$. We have $\mathrm{vert}_{\gamma(v)}(F)>0$ since either $\gamma(v)\in B$ and $\gamma(v)<\max B$, or $\gamma(v)\notin B$ and $\min B<\gamma(v)<\max B$. Conversely, if $j\neq\gamma(v)$ for all $v\in F$, then $j$ is the rightmost leaf in the maximal sub-tree containing $j$. In both cases, $\mathrm{vert}_j(F)=0$. The rest follows from Theorem \ref{thm:GW}.
\end{proof}

\begin{corollary}\label{cor:em-1}
	We have
	$$\fP_F(\mathbf{x};\mathbf{t})\mathfrak{S}_{s_i}(\mathbf{x};\mathbf{t})=C_{F,F}^i(\mathbf{t})\fP_F(\mathbf{x};\mathbf{t})+\sum_{G> F,|G|=|F|}C_{F,G}^i(\mathbf{t})\fP_G(\mathbf{x};\mathbf{t})+\sum_{G>F,|G|=|F|+1}C_{F,G}^i\fP_G(\mathbf{x};\mathbf{t})$$
	for some $C_{F,F}^i(\mathbf{t}),C_{F,G}^i(\mathbf{t})\in\mathbb{Z}[\mathbf{t}]$ that are homogeneous of degree 1, and $C_{F,G}^i\in\mathbb{Z}$.
\end{corollary}

\begin{proof}
	It follows from Theorem \ref{thm:BGNST}, Lemma \ref{lem:e1-1} and comparison of degrees in $\mathbf{x}$ and $\mathbf{t}$.
\end{proof}

\begin{theorem}\label{thm:main1}
	We have
	$$C_{F,F}^i(\mathbf{t})=\Phi_{w_F}(\mathfrak{S}_{s_i})=\sum_{k=1}^i(t_{w_F(k)}-t_k).$$
\end{theorem}

\begin{proof}
	Localizing at $w_F$ on both sides of Corollary \ref{cor:em-1}, and using triangularity gives
	$$\Phi_{w_F}(\fP_F)\Phi_{w_F}(\fS_{s_i})=C_{F,F}^i(\bt)\Phi_{w_F}(\fP_F)+0+0.$$
\end{proof}

\subsection{The equivariant operators}

This subsection expresses $C_{F,G}^i$ using the equivariant operators, and performs certain reductions. In this subsection, we always assume $F<G$ and $G$ has trimming sequence $(j_1,\dots,j_n)$ with $j_1\leq\cdots\leq j_n$.

\begin{lemma}\label{lem:ev1}
	For any $f\in\mathbb{Z}[\bx;\bt]$, and any $j_1,\dots,j_n$, we have
	$$[\mathrm{ev}\star R_{j_1}^-\star\cdots\star R_{j_n}^-]f=f(\bt;\bt).$$
\end{lemma}

\begin{proof}
	For any $j$ and $A$, we have $j\star A=A\cup\{\bar{A}_j\}$ and $\hat{\bt}_{j\star A}=(t_{1,A},\dots,t_{j-1,A},t_{j+1,A},t_{j+2,A},\dots)$.
	\begin{align*}
		[\mathrm{ev}\star R_j^-]_Af(\bx;\bt)&=\mathrm{ev}_{j\star A}R_{j,A}^-f(\bx;\bt)\\
		&=\mathrm{ev}_{j\star A}f(x_1,\dots,x_{j-1},t_{j,A},x_j,x_{j+1},\dots;\bt)\\
		&=f(t_{1,A},\dots,t_{j-1,A},t_{j,A},t_{j+1,A},\dots;\bt)\\
		&=\mathrm{ev}_Af(\bx;\bt).
	\end{align*}
	Therefore,
	\begin{align*}
		[\mathrm{ev}\star R_{j_1}^-\star\cdots\star R_{j_n}^-]&=[\mathrm{ev}\star R_{j_1}^-]_{j_2\star\cdots\star j_n}[R_{j_2}^-\star\cdots\star R_{j_n}^-]=[\mathrm{ev}\star R_{j_2}^-\star\cdots\star R_{j_n}^-]=\mathrm{ev}.
	\end{align*}
\end{proof}

Recall from Definition \ref{def:order} that internal nodes of $G$ are ordered by $v_1^G<\cdots<v_n^G$. Note that the node corresponding to $j_r$ is $v_r^G$.

\begin{lemma}\label{lem:ev2}
	We have
	$$[\mathrm{ev}\star R_{j_1}^-\star\cdots\star R_{j_{r-1}}^-\star E_{j_r}\star R_{j_{r+1}}^-\star\cdots\star R_{j_n}^-]\fS_{s_i}(\bx;\bt)=\mathbf{1}_{\gamma(v_r^G)\leq i< (\rho+\eta)(v_r^G)}.$$
\end{lemma}

\begin{proof}
	We have
	$$[R_{j_{r+1}}^-\star\cdots\star R_{j_n}^-](x_k)=\begin{cases}
		t_k & k\in j_{r+1}\star\cdots\star j_n \\
		x_q & k\notin j_{r+1}\star\cdots\star j_n
	\end{cases}$$
	where $k=(\overline{j_{r+1}\star\cdots\star j_n})_q$. Recall that $\fS_{s_i}(\bx;\bt)=\sum_{k=1}^i \left(x_k- t_k\right)$. We have
	$$[R_{j_{r+1}}^-\star\cdots\star R_{j_n}^-]\fS_{s_i}(\bx;\bt)=\sum_{q\geq 1 \atop (\overline{j_{r+1}\star\cdots\star j_n})_q\leq i}\left(x_q-t_{q,j_{r+1}\star\cdots\star j_n}\right).$$
	Note that
	$$E_{j_r,A}(x_q)=\begin{cases}
		1 & q=j_r \\
		-1 & q=j_r+1 \\
		0 & q\notin\{j_r,j_r+1\}
	\end{cases}.$$
	Moreover, $(\overline{j_{r+1}\star\cdots\star j_n})_{j_r}=\gamma(v_r^G)$ and $(\overline{j_{r+1}\star\cdots\star j_n})_{j_r+1}=(\rho+\eta)(v_r^G)$. Therefore,
	$$[E_{j_r}\star R_{j_{r+1}}^-\star\cdots\star R_{j_n}^-]\fS_{s_i}(\bx;\bt)=\mathbf{1}_{(\overline{j_{r+1}\star\cdots\star j_n})_{j_r}\leq i}-\mathbf{1}_{(\overline{j_{r+1}\star\cdots\star j_n})_{j_r+1}\leq i}=\mathbf{1}_{\gamma(v_r^G)\leq i< (\rho+\eta)(v_r^G)}.$$
	The statement follows from the fact that $\mathrm{ev}_A$ and $R_{j,A}^-$ fix constants.
\end{proof}

\begin{lemma}\label{lem:ev3}
	We have
	$$C_{F,G}^i(\bt)=\sum_{r=1}^n\mathbf{1}_{\gamma(v_r^G)\leq i< (\rho+\eta)(v_r^G)}[\mathrm{ev}\star E_{j_1}\star\cdots \star E_{j_{r-1}}\star R_{j_r}^+ \star E_{j_{r+1}}\star\cdots \star E_{j_n}]\fP_F(\bx;\bt).$$
\end{lemma}

\begin{proof}
	By Theorem \ref{thm:BGNST2}, we have
	$$C_{F,G}^i(\bt)=[\mathrm{ev}\star E_G]\left(\fP_F(\bx;\bt)\fS_{s_i}(\bx;\bt)\right).$$
	Applying Theorem \ref{thm:BGNST3} repeatedly, we have
	$$C_{F,G}^i(\bt)=\sum_{(X_j,Y_j)\in\{(E_j,R_j^-),(R_j^+,E_j)\}}\left([\mathrm{ev}\star X_{j_1}\star\cdots\star X_{j_n}]\fP_F(\bx;\bt)\right)\left([\mathrm{ev}\star Y_{j_1}\star\cdots\star Y_{j_n}]\fS_{s_i}(\bx;\bt)\right).$$
	Note that $E_{j,A}\fS_{s_i}(\bx;\bt)\in\mathbb{Z}[\bt]$ for any $j,A$, and $E_{j,A}f(\bt)=0$. Hence, if $Y_{j_p}=E_{j_p}$ and $Y_{j_q}=E_{j_q}$ for some $p\neq q$, we have $\left([\mathrm{ev}\star Y_{j_1}\star\cdots\star Y_{j_n}]\fS_{s_i}(\bx;\bt)\right)=0$.
	
	Suppose all $Y_j=R_j^-$, by Lemma \ref{lem:ev1}, we have
	$$[\mathrm{ev}\star R_{j_1}^-\star\cdots\star R_{j_n}^-]\fS_{s_i}(\bx;\bt)=\fS_{s_i}(\bt;\bt)=0.$$
	Therefore, the statement follows from Lemma \ref{lem:ev2}.
\end{proof}

For simplicity, we set $\chi_r=[\mathrm{ev}\star E_{j_1}\star\cdots \star E_{j_{r-1}}\star R_{j_r}^+ \star E_{j_{r+1}}\star\cdots \star E_{j_n}]$. Recall $h(F,G)$, $h'(F,G)$ and $S_r$ from Definition \ref{def:h}.

\begin{lemma}\label{lem:ev4}
	 We have
	$$\chi_r\fP_F(\bx;\bt)=\begin{cases}
		\chi_{h(F,G)}\fP_F(\bx;\bt) & h(F,G)\leq r\leq h'(F,G)\\
		0 & r<h(F,G)\text{ or }r>h'(F,G)
	\end{cases}.$$
\end{lemma}

\begin{proof}
	(Case I) If $r<h(F,G)$, then $S_r\nleq_R F$. By Theorem \ref{thm:BGNST4}, we have
	$$[E_{j_{r+1}}\star\cdots\star E_{j_n}]\fP_F(\bx;\bt)=[E_{S_r}]\fP_F(\bx;\bt)=0.$$
	Consequently, $\chi_r\fP_F(\bx;\bt)=0$.
	
	(Case II) If $h(F,G)<r\leq h'(F,G)$, then let $a=j_r=j_{r-1}=\cdots=j_{h(F,G)}$. By Theorem \ref{thm:BGNST6}, we have
	\begin{align*}
		\chi_r\fP_F(\bx;\bt)&=[\mathrm{ev}\star E_{j_1}\star\cdots \star E_{j_{r-2}}\star E_{j_{r-1}}\star R_{j_r}^+ \star E_{j_{r+1}}\star\cdots \star E_{j_n}]\fP_F(\bx;\bt) \\
		&=[\mathrm{ev}\star E_{j_1}\star\cdots \star E_{j_{r-2}}\star R_a^+\star E_a \star E_{j_{r+1}}\star\cdots \star E_{j_n}]\fP_F(\bx;\bt)\\
		&+[\mathrm{ev}\star E_{j_1}\star\cdots \star E_{j_{r-2}}\star R_a^-\star E_{a+1} \star E_{j_{r+1}}\star\cdots \star E_{j_n}]\fP_F(\bx;\bt)\\
		&+(t_{a+1,j_{r+1}\star\cdots\star j_n}-t_{a,j_{r+1}\star\cdots\star j_n})[\mathrm{ev}\star E_{j_1}\star\cdots \star E_{j_{r-2}}\star E_a\star E_a \star E_{j_{r+1}}\star\cdots \star E_{j_n}]\fP_F(\bx;\bt).
	\end{align*}
	The last term vanishes by applying Theorem \ref{thm:BGNST2},
	$$[\mathrm{ev}\star E_{j_1}\star\cdots \star E_{j_{r-2}}\star E_a\star E_a \star E_{j_{r+1}}\star\cdots \star E_{j_n}]\fP_F(\bx;\bt)=[\mathrm{ev}\star E_G]\fP_F(\bx;\bt)=\mathbf{1}_{F=G}=0.$$
	By construction, $S_{r-1}\leq_R F$ and $S_{r-1}$ has trimming sequence $(a,j_{r+1},\dots,j_n)$. Let $H$ have a trimming sequence $(a+1,j_{r+1},\dots,j_n)$. Since left terminal nodes cannot be consecutive, we have $H\nleq_R F$ and the middle term vanishes by Theorem \ref{thm:BGNST2} as
	$$[E_{a+1}\star E_{j_{r+1}}\star\cdots\star E_{j_n}]\fP_F(\bx;\bt)=[E_H]\fP_F(\bx;\bt)=0.$$
	
	Therefore, only the first term survives, and we have
	$$\chi_r\fP_F(\bx;\bt)=[\mathrm{ev}\star E_{j_1}\star\cdots \star E_{j_{r-2}}\star R_{j_{r-1}}^+\star E_{j_r} \star \cdots \star E_{j_n}]\fP_F(\bx;\bt)=\chi_{r-1}\fP_F(\bx;\bt).$$
	
	(Case III) If $r>h'(F,G)$, set $b=j_r\neq a$ let $\ell=\min\{s:j_s=j_r\}$, in particular, $j_\ell>j_{\ell-1}$. By analogous arguments as in Case II, we have $\chi_r\fP_F(\bx;\bt)=\chi_\ell\fP_F(\bx;\bt)$.
	
	If $j_{\ell-1}=b-1$, applying the commutation relations $[E_p\star R_{p+1}^+]_A=[R_p^+\star E_p]_A$ in Theorem \ref{thm:BGNST6} gives
	$$\chi_\ell=[\mathrm{ev}\star E_{j_1}\star\cdots\star E_{j_{\ell-1}}\star R_b^+\star E_{j_{\ell+1}}\star\cdots\star E_{j_n}]=[\mathrm{ev}\star E_{j_1}\star\cdots\star E_{j_{\ell-2}}\star R_{b-1}^+\star E_{b-1}\star E_{j_{\ell+1}}\star\cdots\star E_{j_n}].$$
	Recall that $S_{\ell-1}$ has trimming sequence $(b,j_{\ell+1},\dots,j_n)$ and $S_{\ell-1}\leq_R F$. Let $H$ has trimming sequence $(b-1,j_{\ell+1},\dots,j_n)$. Since left terminal nodes cannot be consecutive, we have $H\nleq_R F$ and $\chi_\ell\fP_F(\bx;\bt)$ vanishes by Theorem \ref{thm:BGNST4} as
	$$[E_{b-1}\star E_{j_{\ell+1}}\star\cdots\star E_{j_n}]\fP_F(\bx;\bt)=[E_H]\fP_F(\bx;\bt)=0.$$
	
	If $j_{\ell-1}<b-1$, applying the commutation relations $[E_p\star R_{q+1}^+]_A=[R_q^+\star E_p]_A$ in Theorem \ref{thm:BGNST5} gives
	$$\chi_\ell=[\mathrm{ev}\star E_{j_1}\star\cdots\star E_{j_{\ell-1}}\star R_b^+\star E_{j_{\ell+1}}\star\cdots\star E_{j_n}]=[\mathrm{ev}\star E_{j_1}\star\cdots\star E_{j_{\ell-2}}\star R_{b-1}^+\star E_{j_{\ell-1}}\star E_{j_{\ell+1}}\star\cdots\star E_{j_n}].$$
	
	We keep applying these commutation relations. At each step, it gives
	$$\chi_\ell=[\mathrm{ev}\star E_{j_1}\star\cdots\star E_{j_{\ell-1}}\star R_b^+\star E_{j_{\ell+1}}\star\cdots\star E_{j_n}]=[\mathrm{ev}\star E_{j_1}\star\cdots \star E_{j_s}\star R_q^+\star E_{a_{s+2}}\star\cdots\star E_{a_n}]$$
	for some $(a_{s+2},\dots,a_n)$. Note that these commutation relations are the same as those for trimming sequences. 
	
	Suppose $\chi_\ell P_F(\bx;\bt)\neq 0$. We show by induction that the following properties hold, with (2) ommited when $s=0$.
	\begin{enumerate}
		\item $(j_1,\dots,j_s,q,a_{s+2},\dots,a_n)$ is a trimming sequence for $G$;
		\item $j_1\leq\cdots\leq j_s<q$;
		\item let $H$ has trimming sequence $(q,a_{s+2},\dots,a_n)$, then $H\leq_R F$.
	\end{enumerate}
	These properties hold before the first commutation, with $s=\ell-1$, $q=b$, and $H=S_{\ell-1}\leq_R F$. Suppose $s\geq 1$, and let $K$ have trimming sequence $(a_{s+2},\ldots,a_n)$. By (3), we have $q\in\operatorname{LTer}(F/K)$. Since $\chi_\ell P_F(x;t)\neq 0$, Theorems \ref{thm:BGNST5} and \ref{thm:BGNST6} gives
	\[
	[E_{j_s}\star E_{a_{s+2}}\star\cdots\star E_{a_n}]
	P_F(x;t) =\chi_\ell P_F(x;t)\neq 0.
	\]
	By Theorem \ref{thm:BGNST4}, this implies $j_s\in\operatorname{LTer}(F/K)$. If $j_s=q-1$, then $j_s$ and $q$ are consecutive left terminal labels, which is absurd. Therefore, $j_s<q-1$. Trimming the terminal node with label $j_s$ changes the terminal label $q$ to $q-1$. Consequently, the forest $H'$ with trimming sequence $(q-1,j_s,a_{s+2},\ldots,a_n)$
	satisfies $H'\leq_R F$. This proves (3) for the next step, where $s$ and $q$ are replaced by $s-1$ and $q-1$, respectively. Property (1) follows from the commutation relations for
	trimming sequences. For $s>1$, we have $j_{s-1}\leq j_s<q-1$, so (2) also holds.
	
	Thus, under the assumption $\chi_\ell P_F(x;t)\neq 0$, the process continues until $s=0$. In this case, (1) and (3) give $G\leq_R F$, which is absurd since $F<G$. Therefore, we have $\chi_\ell P_F(x;t)=0$, as desired.
\end{proof}

\begin{proposition}\label{prop:ev}
	We have
	$$C_{F,G}^i(\bt)=\mathbf{1}_{\rho (v_F^G)\leq i<(\rho+\eta)(v_F^G)}\chi_{h(F,G)}\fP_F(\bx;\bt).$$
\end{proposition}

\begin{proof}
	By Lemmas \ref{lem:ev3} and \ref{lem:ev4}, we have
	$$C_{F,G}^i(\bt)=\left(\sum_{r=h(F,G)}^{h'(F,G)}\mathbf{1}_{\gamma(v_r^G)\leq i<(\rho+\eta)(v_r^G)}\right)\chi_{h(F,G)}\fP_F(\bx;\bt).$$
	By construction $\gamma(v_r^G)=(\rho+\eta)(v_{r+1}^G)$ for all $h(F,G)\leq r< h'(F,G)$, and $\gamma(v_{h'(F,G)}^G)=\rho(v_{h'(F,G)}^G)$. Consequently,
	$$\sum_{r=h(F,G)}^{h'(F,G)}\mathbf{1}_{\gamma(v_r^G)\leq i<(\rho+\eta)(v_r^G)}=\mathbf{1}_{\rho(v_F^G)\leq i<(\rho+\eta)(v_F^G)}.$$
\end{proof}

\begin{lemma}\label{lem:ev5}
	Let $(j_1',\dots,j_n')$ be any trimming sequence of $G$, and suppose $j'_r$ corresponds to the vertex $v_F^G$. Then,
	$$\chi_{h(F,G)}=[\mathrm{ev}\star E_{j'_1}\star\cdots\star E_{j'_{r-1}}\star R_{j'_r}^+\star E_{j'_{r+1}}\star\cdots\star E_{j'_n}].$$
\end{lemma}

\begin{proof}
	This is a direct consequence of Theorem \ref{thm:BGNST5}, since the trimming sequences and the equivariant operators satisfy the same commutation relations.
\end{proof}

Consequently, we are safe to set $\chi_{F,G}=\chi_{h(F,G)}$. We will fix a trimming sequence, indexed for convenience, $(z_m,\dots,z_0,z_{m+1},\dots,z_{n-1})$ of $G$ such that $z_m,\dots,z_1$ correspond to all ancestors of $v_F^G$, and $z_0$ corresponds to $v_F^G$. By Lemma \ref{lem:ev5}, we have
$$\chi_{F,G}=[\mathrm{ev}\star E_{z_m}\star\cdots\star E_{z_1}\star R_{z_0}^+\star E_{z_{m+1}}\star\cdots\star E_{z_{n-1}}].$$
For $0\leq r\leq m$, let $u_r\in G$ be the vertex corresponding to $z_r$. In particular, $v_F^G=u_0$.

\subsection{Admissible sequences}\label{subsec:as}

This subsection expands $\chi_{F,G}$ using the commutation relations from Theorem \ref{thm:BGNST6}. The terms in the expansion will be indexed by admissible sequences. It is helpful to keep in mind that for $1\leq r\leq m$, $u_{r-1}=(u_r)_L$ implies $z_{r-1}=z_r$, and $u_{r-1}=(u_r)_R$ implies $z_{r-1}=z_r+1$.

\begin{definition}
	For $0\leq p\leq m$, a sequence $(\epsilon_0,\dots,\epsilon_p)\in\{+,-\}^{p+1}$ is \textbf{admissible} if
	\begin{enumerate}
		\item $\epsilon_0=+$,
		\item $(\epsilon_{r-1},\epsilon_r)=(+,-)$ implies $u_{r-1}=(u_r)_L$,
		\item $(\epsilon_{r-1},\epsilon_r)=(-,+)$ implies $u_{r-1}=(u_r)_R$.
	\end{enumerate}
\end{definition}

\begin{definition}
	Given an admissible sequence $\epsilon=(\epsilon_0,\dots,\epsilon_p)$, for $0\leq r\leq p$ we define $q_r^\epsilon=z_r+\mathbf{1}_{\epsilon_r=-}$, and
	$$D_r^\epsilon=q_r^\epsilon\star q_{r-1}^\epsilon\star\cdots\star q_1^\epsilon\star z_{m+1}\star\cdots\star z_{n-1}.$$
	In particular, $q_0^\epsilon=z_0$ and $D_0^\epsilon=z_{m+1}\star\cdots\star z_{n-1}$.
\end{definition}

\begin{lemma}\label{lem:ev6}
	We have
	\begin{align*}
		\chi_{F,G}=&\sum_{\epsilon=(\epsilon_0,\dots,\epsilon_m) \atop \text{admissible}}[\mathrm{ev}\star R_{z_m}^{\epsilon_m}\star E_{q_m^\epsilon}\star \cdots\star E_{q_1^\epsilon}\star E_{z_{m+1}}\star\cdots\star E_{z_{n-1}}]\\
		+&\sum_{p=1}^{m}\sum_{\epsilon=(\epsilon_0,\dots,\epsilon_{p-1}) \atop \text{admissible}}\left(\mathbf{1}_{u_{p-1}=(u_p)_L \atop \epsilon_{p-1}=+}+\mathbf{1}_{u_{p-1}=(u_p)_R \atop \epsilon_{p-1}=-}\right)\left(t_{z_p+1,D_{p-1}^\epsilon}-t_{z_p,D_{p-1}^\epsilon}\right)\\
		&\times[\mathrm{ev}\star E_{z_m}\star\cdots\star E_{z_{p-1}}\star E_{q_{p-1}^\epsilon}\star\cdots\star E_{q_1^\epsilon}\star E_{z_{m+1}}\star\cdots\star E_{z_{n-1}}].
	\end{align*}
\end{lemma}

\begin{proof}
	We apply relations in Theorem \ref{thm:BGNST6} to $[E_a\star R_b^\pm]$. Note that the terms in the expansion are indexed by admissible sequences $(\epsilon_0,\dots,\epsilon_p)$, namely,
	\begin{center}
		\begin{tabular}{cccccc}
			vertex relation & trimming relation & $\epsilon_{r-1}$ & local operator & $\epsilon_r$ & chosen term \\
			$u_{r-1}=(u_r)_L$ & $z_{r-1}=z_r$ & $+$ & $[E_{z_r}\star R_{z_{r-1}}^+]_A$ & $+$ & $[R_{z_r}^+\star E_{z_r}]_A$ \\
			$u_{r-1}=(u_r)_L$ & $z_{r-1}=z_r$ & $+$ & $[E_{z_r}\star R_{z_{r-1}}^+]_A$ & $-$ & $[R_{z_r}^-\star E_{z_r+1}]_A$ \\
			$u_{r-1}=(u_r)_L$ & $z_{r-1}=z_r$ & $+$ & $[E_{z_r}\star R_{z_{r-1}}^+]_A$ & $\varnothing$ & $(t_{z_r+1,A}-t_{z_r,A})[E_{z_r}\star E_{z_r}]_A$ \\
			$u_{r-1}=(u_r)_L$ & $z_{r-1}=z_r$ & $-$ & $[E_{z_r}\star R_{z_{r-1}}^-]_A$ & $-$ & $[R_{z_r}^-\star E_{z_r+1}]_A$ \\
			$u_{r-1}=(u_r)_R$ & $z_{r-1}=z_r+1$ & $+$ & $[E_{z_r}\star R_{z_{r-1}}^+]_A$ & $+$ & $[R_{z_r}^+\star E_{z_r}]_A$ \\
			$u_{r-1}=(u_r)_R$ & $z_{r-1}=z_r+1$ & $-$ & $[E_{z_r}\star R_{z_{r-1}}^-]_A$ & $+$ & $[R_{z_r}^+\star E_{z_r}]_A$ \\
			$u_{r-1}=(u_r)_R$ & $z_{r-1}=z_r+1$ & $-$ & $[E_{z_r}\star R_{z_{r-1}}^-]_A$ & $-$ & $[R_{z_r}^-\star E_{z_r+1}]_A$ \\
			$u_{r-1}=(u_r)_R$ & $z_{r-1}=z_r+1$ & $-$ & $[E_{z_r}\star R_{z_{r-1}}^-]_A$ & $\varnothing$ & $(t_{z_r+1,A}-t_{z_r,A})[E_{z_r}\star E_{z_r+1}]_A$
		\end{tabular}
	\end{center}
	Comparing with the definition of $q_r^\epsilon$ and $D_r^\epsilon$ yields the formula.
\end{proof}

\begin{lemma}\label{lem:ev7}
	Suppose $K$ has the trimming sequence $(j_2,\dots,j_p)$ and $H$ has the trimming sequence $(j_1,\dots,j_p)$. Then, we have
	$$[\mathrm{ev}\star R_{j_1}^-\star E_{j_2}\star\cdots\star E_{j_p}]\fP_F=\mathbf{1}_{F=K},$$
	and
	$$[\mathrm{ev}\star R_{j_1}^+\star E_{j_2}\star\cdots\star E_{j_p}]\fP_F=\mathbf{1}_{F=K}+(t_{j_1+1,j_2\star\cdots\star j_p}-t_{j_1,j_2\star\cdots\star j_p})\mathbf{1}_{F=H}.$$
\end{lemma}

\begin{proof}
	Lemma \ref{lem:ev1} gives $[\mathrm{ev}\star R_{j_1}^-]_A=\mathrm{ev}_A$. By Theorem \ref{thm:BGNST2}, we have
	$$[\mathrm{ev}\star R_{j_1}^-\star E_{j_2}\star\cdots\star E_{j_p}]\fP_F=[\mathrm{ev}\star E_K]\fP_F=\mathbf{1}_{F=K}.$$
	Note that we have $R_{j_1,A}^+-R_{j_1,A}^-=(x_{j_1}-t_{j_1,A})E_{j_1,A}$ with $A=j_2\star\cdots\star j_n$. Applying $\mathrm{ev}_{j_1\star A}$ gives
	$$[\mathrm{ev}\star R_{j_1}^+]_A=[\mathrm{ev}\star R_{j_1}^-]_A+(t_{j_1+1,A}-t_{j_1,A})[\mathrm{ev}\star E_{j_1}]_A$$
	as $\mathrm{ev}_{j_1\star A}(x_{j_1})=t_{j_1,j_1\star A}=t_{j_1+1,A}$. Applying $[\mathrm{ev}\star E_H]\fP_F=\mathbf{1}_{F=H}$ gives
	$$[\mathrm{ev}\star R_{j_1}^+\star E_{j_2}\star\cdots\star E_{j_p}]\fP_F=\mathbf{1}_{F=K}+(t_{j_1+1,j_2\star\cdots\star j_p}-t_{j_1,j_2\star\cdots\star j_p})\mathbf{1}_{F=H}.$$
\end{proof}

\begin{definition}
	Suppose $S$ and $T$ are two trees (possibly trivial trees), we set $S\vee T$ to be the tree obtained by attaching the root of $S$, $v$, and the root of $T$, $u$, by a new root $w$ such that $w_L=v$ and $w_R=u$.
\end{definition}

\begin{definition}\label{def:e1}
	Given an admissible sequence $(\epsilon_0,\dots,\epsilon_p)$, we construct a sequence of pairs $((A_0,B_0),\dots,(A_p,B_p))$ as follows.
	\begin{enumerate}
		\item $A_0$ and $B_0$ are the subtrees rooted at $(v_F^G)_L$ and $(v_F^G)_R$, respectively;
		\item if $u_{r-1}=(u_r)_L$ and $\epsilon_r=+$, then $(A_r,B_r)=(A_{r-1}\vee B_{r-1},F_r)$;
		\item if $u_{r-1}=(u_r)_L$ and $\epsilon_r=-$, then $(A_r,B_r)=(A_{r-1},B_{r-1}\vee F_r)$;
		\item if $u_{r-1}=(u_r)_R$ and $\epsilon_r=+$, then $(A_r,B_r)=(F_r\vee A_{r-1},B_{r-1})$;
		\item if $u_{r-1}=(u_r)_R$ and $\epsilon_r=-$, then $(A_r,B_r)=(F_r,A_{r-1}\vee B_{r-1})$,
	\end{enumerate}
	where $F_r$ is the subtree of $G$ rooted at $(u_r)_L$ (or $(u_r)_R$) if $u_{r-1}=(u_r)_R$ (or $(u_r)L$, respectively).
\end{definition}

\begin{lemma}\label{lem:ev8}
	The admissible sequence can be recovered from $(A_p,B_p)$.
\end{lemma}

\begin{proof}
	The sequence of pairs and admissible sequence can be recovered as follows. For $r$ goes from $p$ to $1$, we have
	\begin{enumerate}
		\item If $u_{r-1}=(u_r)_L$ and $B_r=F_r$, then $A_{r-1}\vee B_{r-1}=A_r$ and $\epsilon_r=+$;
		\item If $u_{r-1}=(u_r)_L$ and $B_r\neq F_r$, then $A_{r-1}=A_r$, $B_{r-1}\vee F_r=B_r$ and $\epsilon_r=-$;
		\item If $u_{r-1}=(u_r)_R$ and $A_r\neq F_r$, then $F_r\vee A_{r-1}=A_r$, $B_{r-1}=B_r$ and $\epsilon_r=+$;
		\item If $u_{r-1}=(u_r)_R$ and $A_r=F_r$, then $A_{r-1}\vee B_{r-1}=B_r$ and $\epsilon_r=-$.
	\end{enumerate}
\end{proof}

\begin{definition}\label{def:e2}
	Given an admissible sequence $\epsilon=(\epsilon_0,\dots,\epsilon_m)$ and its corresponding sequence of pairs $((A_0,B_0),\dots,(A_m,B_m))$, we define
	\begin{enumerate}
		\item if $p=m$, $K^\epsilon$ is the forest obtained from $G$ by replacing the component rooted at $u_m$ with two components, $A_m$ and $B_m$;
		\item if $p=m$, $H^\epsilon$ is the forest obtained from $G$ by replacing the component rooted at $u_m$ with $A_m\vee B_m$;
	\end{enumerate}
	For $1\leq p\leq m$, given an admissible sequence $\epsilon=(\epsilon_0,\dots,\epsilon_{p-1})$ and its corresponding sequence of pairs $((A_0,B_0),\dots,(A_{p-1},B_{p-1}))$, we define $H^\epsilon_p$ to be the forest obtained from $G$ by replacing the subtree rooted at $u_p$ with
		$$\begin{cases}
			(A_{p-1}\vee B_{p-1})\vee F_p & u_{p-1}=(u_p)_L \\
			F_p\vee (A_{p-1}\vee B_{p-1}) & u_{p-1}=(u_p)_R
		\end{cases}.$$
		Moreover, we set $w_p^\epsilon$ to be the root of $A_p\vee B_p$ in $H_{p+1}^\epsilon$, and by convention $H_{m+1}^\epsilon=H^\epsilon$.
\end{definition}

It is helpful to remember that $H_p^\epsilon$ is the forest obtained from $G$ by replacing the subtree rooted at $u_{p-1}$ by $A_{p-1}\vee B_{p-1}$.

\begin{lemma}\label{lem:ev9}
	Given an admissible sequence $\epsilon=(\epsilon_0,\dots,\epsilon_p)$. For $p=m$, we have
	\begin{enumerate}
		\item $(q_m^\epsilon,\dots,q_1^\epsilon,z_{m+1},\dots,z_{n-1})$ is a trimming sequence for $K^\epsilon$;
		\item $(z_m,q_m^\epsilon,\dots,q_1^\epsilon,z_{m+1},\dots,z_{n-1})$ is a trimming sequence for $H^\epsilon$;
		\item $t_{z_m+1,D_m^\epsilon}=t_{(\rho+\eta)(u_m)}$ and $t_{z_m,D_m^\epsilon}=t_{\gamma(w_m^\epsilon)}$.
	\end{enumerate}
	For $1\leq p\leq m$, we have
	\begin{enumerate}
		\item $(z_m,\dots,z_{p-1},q_{p-1}^\epsilon,\dots,q_1^\epsilon,z_{m+1},\dots,z_{n-1})$ is a trimming sequence for $H_p^\epsilon$;
		\item if $u_{p-1}=(u_p)_L$ then $t_{z_p,D_{p-1}^\epsilon}=t_{\gamma(w^\epsilon_{p-1})}$ and $t_{z_p+1,D_{p-1}^\epsilon}=t_{\gamma(u_p)}$;
		\item if $u_{p-1}=(u_p)_R$ then $t_{z_p,D_{p-1}^\epsilon}=t_{\gamma(u_p)}$ and $t_{z_p+1,D_{p-1}^\epsilon}=t_{\gamma(w^\epsilon_{p-1})}$.
	\end{enumerate}
\end{lemma}

\begin{proof}
	The trimming sequence can be determined by verifying the following cases. We label each vertex by its corresponding entry in the trimming sequence.
	
	When $u_{r-1}=(u_r)_L$, we have
	\begin{center}
		\begin{tikzpicture}
			\node at (0,0) {\scalebox{0.8}{\begin{tikzpicture}
					\node at (0,0) {$A_{r-1}$};
					\node at (2,0) {$B_{r-1}$};
					\node at (4,0) {$F_r$};
					
					\draw[dashed] (1,1.3) -- (2,0.3);
					\draw[dashed] (0,0.3) -- (2,2.3) -- (4,0.3);
					
					\node at (0.5,1.4) {$z_{r-1}$};
					\node at (2,2.6) {$z_r$};
					
					\node at (2,-0.6) {$G$};
				\end{tikzpicture}}};
			\node at (5.5,0) {\scalebox{0.8}{\begin{tikzpicture}
					\node at (0,0) {$A_{r-1}$};
					\node at (2,0) {$B_{r-1}$};
					\node at (4,0) {$F_r$};
					
					\draw[dashed] (2,0.3) -- (1,1.3);
					\draw[dashed] (0,0.3) -- (2,2.3) -- (4,0.3);
					
					\node at (0.6,1.4) {$q_r^\epsilon$};
					\node at (2,2.6) {$z_r$};
					
					\node at (2,-0.6) {$H^\epsilon_r$, $\epsilon_{r}=+$};
			\end{tikzpicture}}};
			\node at (11,0) {\scalebox{0.8}{\begin{tikzpicture}
					\node at (0,0) {$A_{r-1}$};
					\node at (2,0) {$B_{r-1}$};
					\node at (4,0) {$F_r$};
					
					\draw[dashed] (2,0.3) -- (3,1.3);
					\draw[dashed] (0,0.3) -- (2,2.3) -- (4,0.3);
					
					\node at (3.4,1.4) {$q_r^\epsilon$};
					\node at (2,2.6) {$z_r$};
					
					\node at (2,-0.6) {$H^\epsilon_r$, $\epsilon_{r}=-$};
				\end{tikzpicture}}};
		\end{tikzpicture}
	\end{center}
	When $u_{r-1}=(u_r)_R$, we have
	\begin{center}
		\begin{tikzpicture}
			\node at (0,0) {\scalebox{0.8}{\begin{tikzpicture}
						\node at (2,0) {$A_{r-1}$};
						\node at (4,0) {$B_{r-1}$};
						\node at (0,0) {$F_r$};
						
						\draw[dashed] (3,1.3) -- (2,0.3);
						\draw[dashed] (0,0.3) -- (2,2.3) -- (4,0.3);
						
						\node at (3.5,1.4) {$z_{r-1}$};
						\node at (2,2.6) {$z_r$};
						
						\node at (2,-0.6) {$G$};
			\end{tikzpicture}}};
			\node at (5.5,0) {\scalebox{0.8}{\begin{tikzpicture}
						\node at (2,0) {$A_{r-1}$};
						\node at (4,0) {$B_{r-1}$};
						\node at (0,0) {$F_r$};
						
						\draw[dashed] (2,0.3) -- (1,1.3);
						\draw[dashed] (0,0.3) -- (2,2.3) -- (4,0.3);
						
						\node at (0.6,1.4) {$q_r^\epsilon$};
						\node at (2,2.6) {$z_r$};
						
						\node at (2,-0.6) {$H^\epsilon_r$, $\epsilon_{r}=+$};
			\end{tikzpicture}}};
			\node at (11,0) {\scalebox{0.8}{\begin{tikzpicture}
						\node at (2,0) {$A_{r-1}$};
						\node at (4,0) {$B_{r-1}$};
						\node at (0,0) {$F_r$};
						
						\draw[dashed] (2,0.3) -- (3,1.3);
						\draw[dashed] (0,0.3) -- (2,2.3) -- (4,0.3);
						
						\node at (3.4,1.4) {$q_r^\epsilon$};
						\node at (2,2.6) {$z_r$};
						
						\node at (2,-0.6) {$H^\epsilon_r$, $\epsilon_{r}=-$};
			\end{tikzpicture}}};
		\end{tikzpicture}
	\end{center}
	In both cases, we have $q_r^\epsilon=z_r+\mathbf{1}_{\epsilon_r=-}$, as desired. This proves the statements for $K^\epsilon$, $H^\epsilon$ and $H_p^\epsilon$.
	
	Let $U_p^\epsilon$ be the forest with trimming sequence $(q_p^\epsilon,\dots,q_1^\epsilon,z_{m+1},\dots,z_{n-1})$. By Lemma \ref{lem:BGNST}, $D_p^\epsilon=\{\gamma(v):v\in U_p^\epsilon\}$. For $p=m$, by construction, we have $\overline{\{\gamma(v):v\in K^\epsilon\}}_{z_m}=\gamma(w_m^\epsilon)$ and $\overline{\{\gamma(v):v\in K^\epsilon\}}_{z_m+1}=(\rho+\eta)(u_m)$. This proves the statement for $t_{z_m+1,D_m^\epsilon}$ and $t_{z_m,D_m^\epsilon}$.
	
	For $1\leq p\leq m$, if $u_{p-1}=(u_p)_L$, by construction, we have $\overline{\{\gamma(v):v\in U_{p-1}^\epsilon\}}_{z_p}=\gamma(w_{p-1}^\epsilon)$ and $\overline{\{\gamma(v):v\in U_{p-1}^\epsilon\}}_{z_p+1}=\gamma(u_p)$. If $u_{p-1}=(u_p)_R$, we have $\overline{\{\gamma(v):v\in U_{p-1}^\epsilon\}}_{z_p}=\gamma(u_p)$ and $\overline{\{\gamma(v):v\in U_{p-1}^\epsilon\}}_{z_p+1}=\gamma(w_{p-1}^\epsilon)$. This settles the last two statements.
\end{proof}

\begin{definition}
	For $1\leq p\leq m$ and an admissible sequence $\epsilon=(\epsilon_0,\dots,\epsilon_{p-1})$, we set
	$$\psi_p^\epsilon=\begin{cases}
		t_{\gamma(u_p)}-t_{\gamma(w_{p-1}^\epsilon)} & u_{p-1}=(u_p)_L \text{ and }\epsilon_{p-1}=+ \\
		t_{\gamma(w_{p-1}^\epsilon)}-t_{\gamma(u_p)} & u_{p-1}=(u_p)_R \text{ and }\epsilon_{p-1}=- \\
		0 & \text{otherwise}
	\end{cases}.$$
\end{definition}

\subsection{The case $|G|=|F|+1$}

\begin{theorem}\label{thm:main2}
	Suppose $F<G$ and $|G|=|F|+1$. Then $C_{F,G}^i=1$ if the following hold.
	\begin{enumerate}
		\item $\rho(v_F^G)\leq i<(\rho+\eta)(v_F^G)$;
		\item Let $T$ be the full binary tree component in $G$ containing $v_F^G$ i.e. the subtree rooted at $u_m$, $F$ is obtained from $G$ by replacing $T$ with two components $X,Y$;
		\item the pair $(X,Y)$ recovers an admissible sequence.
	\end{enumerate}
	Otherwise, $C_{F,G}^i=0$.
\end{theorem}

\begin{proof}
	By Corollary \ref{cor:em-1}, $C_{F,G}^i\in\mathbb{Z}$. The rest follows from Proposition \ref{prop:ev}, Lemma \ref{lem:ev6}, \ref{lem:ev7}, \ref{lem:ev8} and \ref{lem:ev9}.
\end{proof}

\subsection{The case $|G|=|F|$}

In this section, we always assume $F<G$ and $|F|=|G|$.

\begin{lemma}\label{lem:tele1}
	We have
	$$\chi_{F,G}\fP_{F}(\bx;\bt)=\sum_{p=1}^m\sum_{\epsilon=(\epsilon_0,\dots,\epsilon_{p-1}) \atop \text{admissible}}\psi_p^\epsilon\mathbf{1}_{F=H_p^\epsilon}+\sum_{\epsilon=(\epsilon_0\dots,\epsilon_{m-1},+) \atop \text{admissible}}\left(t_{(\rho+\eta)(u_m)}-t_{\gamma(w_m^\epsilon)}\right)\mathbf{1}_{F=H^\epsilon}.$$
	Moreover, every $\psi_p^\epsilon$ and $t_{(\rho+\eta)(u_m)}-t_{\gamma(w_m^\epsilon)}$ is a non-negative root.
\end{lemma}

\begin{proof}
	By Corollary \ref{cor:em-1}, $\chi_{F,G}\fP_F(\bx;\bt)\in\mathbb{Z}[t]$ is homogeneous of degree $1$. Removing the constant terms from Lemma \ref{lem:ev6} and \ref{lem:ev7} yields the formula. Non-negativity of the coefficients follows by comparing $u_p$ and $w_{p-1}^\epsilon$ in all cases e.g. see the figure in Lemma \ref{lem:ev9}.
\end{proof}

\begin{lemma}\label{lem:tele2}
	For $1\leq p\leq m$ and an admissible sequence $\epsilon=(\epsilon_0,\dots,\epsilon_p)$, we have
	\begin{enumerate}
		\item If $u_{p-1}=(u_p)_L$ and $\epsilon_p=+$, then $\gamma(w_p^\epsilon)=\gamma(u_p)$;
		\item If $u_{p-1}=(u_p)_L$ and $\epsilon_p=-$, then $\gamma(w_p^\epsilon)<\gamma(u_p)$;
		\item If $u_{p-1}=(u_p)_R$ and $\epsilon_p=+$, then $\gamma(w_p^\epsilon)>\gamma(u_p)$;
		\item If $u_{p-1}=(u_p)_R$ and $\epsilon_p=-$, then $\gamma(w_p^\epsilon)=\gamma(u_p)$.
	\end{enumerate}
	Moreover, in cases (1) and (4), we have $H_{p+1}^\epsilon=H_p^{(\epsilon_0,\dots,\epsilon_{p-1})}$.
\end{lemma}

\begin{proof}
	The four cases can be verified directly as follows.
	
	When $u_{p-1}=(u_p)_L$, we have
	\begin{center}
		\begin{tikzpicture}
			\node at (0,0) {\scalebox{0.8}{\begin{tikzpicture}
						\node at (0,0) {$A_{p-1}$};
						\node at (2,0) {$B_{p-1}$};
						\node at (4,0) {$F_p$};
						
						\draw[dashed] (1,1.3) -- (2,0.3);
						\draw[dashed] (0,0.3) -- (2,2.3) -- (4,0.3);
						
						\node at (0.5,1.4) {$u_{p-1}$};
						\node at (2,2.6) {$u_p$};
						
						\node at (2,-0.6) {$G$};
			\end{tikzpicture}}};
			\node at (5.5,0) {\scalebox{0.8}{\begin{tikzpicture}
						\node at (0,0) {$A_{p-1}$};
						\node at (2,0) {$B_{p-1}$};
						\node at (4,0) {$B_p=F_p$};
						
						\draw[dashed] (2,0.3) -- (1,1.3);
						\draw[dashed] (0,0.3) -- (2,2.3) -- (4,0.3);
						
						\node at (0.6,1.4) {$A_p$};
						\node at (2,2.6) {$w_p^\epsilon$};
						
						\node at (2,-0.6) {$H^\epsilon_p$, $\epsilon_p=+$};
			\end{tikzpicture}}};
			\node at (11,0) {\scalebox{0.8}{\begin{tikzpicture}
						\node at (0,0) {$A_{p-1}=A_p$};
						\node at (2,0) {$B_{p-1}$};
						\node at (4,0) {$F_p$};
						
						\draw[dashed] (2,0.3) -- (3,1.3);
						\draw[dashed] (0,0.3) -- (2,2.3) -- (4,0.3);
						
						\node at (3.4,1.4) {$B_p$};
						\node at (2,2.6) {$w_p^\epsilon$};
						
						\node at (2,-0.6) {$H^\epsilon_p$, $\epsilon_{p}=-$};
			\end{tikzpicture}}};
		\end{tikzpicture}
	\end{center}
	When $u_{p-1}=(u_p)_R$, we have
	\begin{center}
		\begin{tikzpicture}
			\node at (0,0) {\scalebox{0.8}{\begin{tikzpicture}
						\node at (2,0) {$A_{p-1}$};
						\node at (4,0) {$B_{p-1}$};
						\node at (0,0) {$F_p$};
						
						\draw[dashed] (3,1.3) -- (2,0.3);
						\draw[dashed] (0,0.3) -- (2,2.3) -- (4,0.3);
						
						\node at (3.5,1.4) {$u_{p-1}$};
						\node at (2,2.6) {$u_p$};
						
						\node at (2,-0.6) {$G$};
			\end{tikzpicture}}};
			\node at (5.5,0) {\scalebox{0.8}{\begin{tikzpicture}
						\node at (2,0) {$A_{p-1}$};
						\node at (4,0) {$B_{p-1}=B_p$};
						\node at (0,0) {$F_p$};
						
						\draw[dashed] (2,0.3) -- (1,1.3);
						\draw[dashed] (0,0.3) -- (2,2.3) -- (4,0.3);
						
						\node at (0.6,1.4) {$A_p$};
						\node at (2,2.6) {$w_p^\epsilon$};
						
						\node at (2,-0.6) {$H^\epsilon_p$, $\epsilon_{p}=+$};
			\end{tikzpicture}}};
			\node at (11,0) {\scalebox{0.8}{\begin{tikzpicture}
						\node at (2,0) {$A_{p-1}$};
						\node at (4,0) {$B_{p-1}$};
						\node at (0,0) {$A_p=F_p$};
						
						\draw[dashed] (2,0.3) -- (3,1.3);
						\draw[dashed] (0,0.3) -- (2,2.3) -- (4,0.3);
						
						\node at (3.4,1.4) {$B_p$};
						\node at (2,2.6) {$w_p^\epsilon$};
						
						\node at (2,-0.6) {$H^\epsilon_p$, $\epsilon_{p}=-$};
			\end{tikzpicture}}};
		\end{tikzpicture}
	\end{center}
	To verify $H_{p+1}^\epsilon=H_p^{(\epsilon_0,\dots,\epsilon_{p-1})}$, both $(1)$ and $(4)$ can be combined into three cases.
	
	When $p<m$ and $u_p=(u_{p+1})_L$,
	\begin{center}
		\begin{tikzpicture}
			\node at (0,0) {\scalebox{0.8}{\begin{tikzpicture}
						\node at (0,0) {$A_{p}$};
						\node at (2,0) {$B_{p}$};
						\node at (4,0) {$F_{p+1}$};
						
						\draw[dashed] (1,1.3) -- (2,0.3);
						\draw[dashed] (0,0.3) -- (2,2.3) -- (4,0.3);
						
						\node at (0.5,1.4) {$u_p$};
						\node at (2,2.6) {$u_{p+1}$};
						
						\node at (2,-0.6) {$G$};
			\end{tikzpicture}}};
			\node at (5.5,0) {\scalebox{0.8}{\begin{tikzpicture}
						\node at (0,0) {$A_p$};
						\node at (2,0) {$B_p$};
						\node at (4,0) {$F_{p+1}$};
						
						\draw[dashed] (2,0.3) -- (1,1.3);
						\draw[dashed] (0,0.3) -- (2,2.3) -- (4,0.3);
						
						%\node at (0.6,1.4) {$A_{p+1}$};
						\node at (2,2.6) {$w_{p+1}^\epsilon$};
						
						\node at (2,-0.6) {$H^\epsilon_{p+1}$};
			\end{tikzpicture}}};
		\end{tikzpicture}
	\end{center}
	When $p<m$ and $u_p=(u_{p+1})_R$,
	\begin{center}
		\begin{tikzpicture}
			\node at (0,0) {\scalebox{0.8}{\begin{tikzpicture}
						\node at (2,0) {$A_{p}$};
						\node at (4,0) {$B_{p}$};
						\node at (0,0) {$F_{p+1}$};
						
						\draw[dashed] (3,1.3) -- (2,0.3);
						\draw[dashed] (0,0.3) -- (2,2.3) -- (4,0.3);
						
						\node at (3.5,1.4) {$u_{p}$};
						\node at (2,2.6) {$u_{p+1}$};
						
						\node at (2,-0.6) {$G$};
			\end{tikzpicture}}};
			\node at (5.5,0) {\scalebox{0.8}{\begin{tikzpicture}
						\node at (2,0) {$A_{p}$};
						\node at (4,0) {$B_{p}$};
						\node at (0,0) {$F_{p+1}$};
						
						\draw[dashed] (2,0.3) -- (3,1.3);
						\draw[dashed] (0,0.3) -- (2,2.3) -- (4,0.3);
						
						%\node at (3.4,1.4) {$B_p$};
						\node at (2,2.6) {$w_{p+1}^\epsilon$};
						
						\node at (2,-0.6) {$H^\epsilon_{p+1}$};
			\end{tikzpicture}}};
		\end{tikzpicture}
	\end{center}
	When $p=m$,
	\begin{center}
		\begin{tikzpicture}
			\node at (0,0) {\scalebox{0.8}{\begin{tikzpicture}
						\node at (2,0) {$A_{m}$};
						\node at (4,0) {$B_{m}$};
						
						\draw[dashed] (4,0.3) -- (3,1.3) -- (2,0.3);
						
						\node at (3,1.5) {$u_{m}$};
						
						\node at (3,-0.6) {$G$};
			\end{tikzpicture}}};
			\node at (5.5,0) {\scalebox{0.8}{\begin{tikzpicture}
						\node at (2,0) {$A_{m}$};
						\node at (4,0) {$B_{m}$};
						
						\draw[dashed] (4,0.3) -- (3,1.3) -- (2,0.3);
						
						\node at (3,1.5) {$w$};
						
						\node at (3,-0.6) {$H^\epsilon=H_{m+1}^\epsilon$};
			\end{tikzpicture}}};
		\end{tikzpicture}
	\end{center}
	Comparing the tree structure yields the result.
\end{proof}

\begin{definition}\label{def:tele}
	Suppose $F<G$ and $|G|=|F|$, we say $(F,G)$ is \textbf{compatible} if there exists a unique $k\in[m]$ satisfying the following.
	\begin{enumerate}
		\item $F$ is obtained from $G$ by replacing the subtree rooted at $u_k$ by $X\vee Y$ that has the same size;
		\item let $x$ be the root of $X\vee Y$ in $F$, then $\gamma(x)\neq\gamma(u_k)$;
		\item the pair $(X,Y)$ recovers an admissible sequence $(\epsilon_0,\dots,\epsilon_k)$ and $\epsilon_k=+$ if and only if $\gamma(x)>\gamma(u_k)$.
	\end{enumerate}
	When $(F,G)$ is compatible, we denote the root of $X\vee Y$ in $F$ by $x_{F,G}=x$, and $k_{F,G}=k$.
\end{definition}

\begin{lemma}\label{lem:tele3}
	$(F,G)$ is compatible if and only if for some $p$, and some admissible sequence $\epsilon=(\epsilon_0,\dots,\epsilon_{p-1})$, either $\psi_p^\epsilon\mathbf{1}_{F=H_p^\epsilon}\neq 0$ or $(t_{\rho+\eta}(u_m)-t_{\gamma(w_m^\epsilon)})\mathbf{1}_{F=H^\epsilon}\neq 0$.
\end{lemma}

\begin{proof}
	($\implies$) Let $k=k_{F,G}$ and let $\epsilon=(\epsilon_0,\dots,\epsilon_k)$ be the recovered admissible sequence. By Lemma \ref{lem:ev9}, $F=H_{k+1}^\epsilon$.
	
	($\impliedby$) Such a non-zero term gives $F=H_p^\epsilon$ for some $1\leq p\leq m+1$. Find the largest $p$. Suppose $\gamma(w_{p-1}^\epsilon)=\gamma(u_{p-1})$, Lemma \ref{lem:tele2} implies $F=H_p^\epsilon=H_{p-1}^{(\epsilon_0,\dots,\epsilon_{p-2})}$. Then we compare $\gamma(w_{p-1}^{(\epsilon_0,\dots,\epsilon_{p-2})})$ with $\gamma(u_{p-2})$, and so on.
	
	Let $k$ be the largest number such that $\gamma(w_k^{\epsilon_0,\dots,\epsilon_{k-1}})\neq\gamma(u_k)$. Such a term must exist, as we would have $F=G$ otherwise. Setting $x_{F,G}=w_k^{(\epsilon_0,\dots,\epsilon_{k-1})}$, and comparing Definition \ref{def:tele} with Lemma \ref{lem:tele2} yields the result.
\end{proof}

\begin{comment}

\begin{lemma}\label{lem:tele4}
	Suppose $(F,G)$ is compatible, $\chi_{F,G}\fP_F(\bx;\bt)\neq0$ if and only if one of the following holds:
	\begin{enumerate}
		\item $\gamma(x_{F,G})>\gamma(u_{k_{F,G}})$ and either $x_{F,G}$ is a root, or $x_{F,G}=v_L$ for some $v\in F$;
		\item $\gamma(x_{F,G})<\gamma(u_{k_{F,G}})$ and $x_{F,G}=v_R$ for some $v\in F$.
	\end{enumerate}
\end{lemma}

\begin{proof}
	Since all coefficients in the formula of Lemma \ref{lem:tele1} is a non-negative root, $\chi_{F,G}\fP_F(\bx;\bt)\neq 0$ if and only if at least one of the terms, $\psi_p^\epsilon\mathbf{1}_{F=H_p^\epsilon}$ or $(t_{(\rho+\eta)(u_m)}-t_{\gamma(w_m^\epsilon)})\mathbf{1}_{F=H^\epsilon}$, is non-zero.
	
	($\impliedby$) Let $k=k_{F,G}$ and suppose $\gamma(x_{F>G})>\gamma(u_{k_{F,G}})$, we have $\epsilon_{k}=+$ since $(F,G)$ is compatible. If $x_{F,G}=v_L$ for some $v\in F$, then $u_k=(u_{k+1})_L$
	
	($\implies$) Suppose $\gamma(x_{F,G})>\gamma(u_k)$, let $F=H_p^\epsilon$ for some $1\leq p\leq m$ and $\epsilon=(\epsilon_0,\dots,\epsilon_{p-1})$. 
\end{proof}

\end{comment}

\begin{lemma}\label{lem:tele5}
	Suppose $(F,G)$ is compatible, and $\gamma(x_{F,G})>\gamma(u_k)$. Let $x\in F$ be the vertex with maximal $\eta(x)$ among all vertices satisfying $\rho(x)=\rho(x_{F,G})$. Then, we have
	$$\chi_{F,G}\fP_F(\bx;\bt)=\begin{cases}
		t_{(\rho+\eta)(x)}-t_{\gamma(x_{F,G})} & x\text{ is a root of a component in }F,\\
		t_{\gamma(x)}-t_{\gamma(x_{F,G})} & {otherwise}.\\
	\end{cases}$$
\end{lemma}

\begin{proof}
	Let $k=k_{F,G}$ and let $(\epsilon_0,\dots,\epsilon_k)$ be the admissible sequence recovered from the pair in Definition \ref{def:tele}. Then $\epsilon_k=+$ and $F=H_{k+1}^{\epsilon}$, with $w_k^{\epsilon}=x_{F,G}$. Let $u_\ell$ be the vertex of $G$ corresponding to $x$. Thus, $u_{r-1}=(u_r)_L$ for $k<r\leq\ell$, and either $\ell=m$, or $u_\ell=(u_{\ell+1})_R$.
	
	Suppose $F=H_p^\epsilon$. We must have $p\geq k+1$, since replacing a subtree strictly below $u_k$ does not change the splitting label at $u_k$. For $k<r<p$, the splitting label at $u_r$ is unchanged. Applying Lemma \ref{lem:tele2} successively, starting with $r=p-1$, shows that this representation restricts to the recovered sequence $(\epsilon_0,\dots,\epsilon_k)$, and that
	$$\epsilon_r=\begin{cases}
		+ & u_{r-1}=(u_r)_L,\\
		- & u_{r-1}=(u_r)_R.
	\end{cases}$$
	Conversely, each admissible extension with these signs gives $F=H_p^\epsilon$ by Lemma \ref{lem:tele2}. Since $\epsilon_k=+$, these extensions are precisely those with $k+1\leq p\leq\ell+1$, and all their signs after $\epsilon_k$ are $+$. Indeed, if $\ell<m$, extending past $\ell$ would require $(\epsilon_\ell,\epsilon_{\ell+1})=(+,-)$ at a right child, which is not admissible. Uniqueness follows from Lemma \ref{lem:ev8}.
	
	For $k<r\leq\ell$, we have $\gamma(w_r^\epsilon)=\gamma(u_r)$, while $w_k^\epsilon=x_{F,G}$. If $\ell<m$, the term with $p=\ell+1$ vanishes, and Lemma \ref{lem:tele1} gives
	\begin{align*}
		\chi_{F,G}\fP_F(\bx;\bt)
		&=\sum_{p=k+1}^{\ell}\left(t_{\gamma(u_p)}-t_{\gamma(w_{p-1}^\epsilon)}\right)\\
		&=t_{\gamma(x)}-t_{\gamma(x_{F,G})}.
	\end{align*}
	If $\ell=m$, the last sum in Lemma \ref{lem:tele1} contributes one further term, and we obtain
	\begin{align*}
		\chi_{F,G}\fP_F(\bx;\bt)
		&=\sum_{p=k+1}^{m}\left(t_{\gamma(u_p)}-t_{\gamma(w_{p-1}^\epsilon)}\right)
		+t_{(\rho+\eta)(u_m)}-t_{\gamma(w_m^\epsilon)}\\
		&=t_{(\rho+\eta)(x)}-t_{\gamma(x_{F,G})}.
	\end{align*}
	This proves the statement.
\end{proof}

\begin{lemma}\label{lem:tele6}
	Suppose $(F,G)$ is compatible, and $\gamma(x_{F,G})<\gamma(u_k)$. Let $x\in F$ be the vertex with maximal $\eta(x)$ among all vertices satisfying $(\rho+\eta)(x)=(\rho+\eta)(x_{F,G})$. Then, we have
	$$\chi_{F,G}\fP_F(\bx;\bt)=t_{\gamma(x_{F,G})}-t_{\gamma(x)}.$$
\end{lemma}

\begin{proof}
	Let $\epsilon=(\epsilon_0,\dots,\epsilon_k)$ be the recovered admissible sequence. By Definition \ref{def:tele}, we have $\epsilon_k=-$ and $F=H_{k+1}^\epsilon$.
	
	The vertices with the same rightmost leaf as $x_{F,G}$ are obtained by following the ancestors of $x_{F,G}$ as long as the current vertex is a right child. Let $p$ be the index such that $x$ corresponds to $u_p$. Thus, $k\leq p\leq m$, we have $u_{r-1}=(u_r)_R$ for $k<r\leq p$, and either $p=m$ or $u_p=(u_{p+1})_L$.
	
	Extend $\epsilon$ by setting $\epsilon_r=-$ for $k<r\leq p$. These extensions are admissible, and Lemma \ref{lem:tele2} gives
	$$F=H_{k+1}^{(\epsilon_0,\dots,\epsilon_k)}=\cdots=H_{p+1}^{(\epsilon_0,\dots,\epsilon_p)},\qquad
	\gamma(w_r^\epsilon)=\gamma(u_r)\quad(k<r\leq p).$$
	Moreover, $w_k^\epsilon=x_{F,G}$ and $w_p^\epsilon=x$ under these identifications.
	
	As in the proof of Lemma \ref{lem:tele5}, every representation of $F$ must extend the recovered sequence, and preserving each subsequent splitting label forces the sign in Lemma \ref{lem:tele2}. If $p<m$, preserving the split at $u_{p+1}$ would require $\epsilon_{p+1}=+$, which is not admissible since $\epsilon_p=-$ and $u_p=(u_{p+1})_L$. Thus, no further $H_r$ can equal $F$, and $\psi_{p+1}^{(\epsilon_0,\dots,\epsilon_p)}=0$. If $p=m$, the final term in Lemma \ref{lem:tele1} is absent since $\epsilon_m=-$.
	
	Consequently, Lemma \ref{lem:tele1} gives
	\begin{align*}
		\chi_{F,G}\fP_F(\bx;\bt)
		&=\sum_{r=k+1}^{p}\left(t_{\gamma(w_{r-1}^\epsilon)}-t_{\gamma(u_r)}\right)\\
		&=t_{\gamma(x_{F,G})}-t_{\gamma(x)}.
	\end{align*}
	When $p=k$, the sum is empty and $x=x_{F,G}$, so the same formula holds.
\end{proof}

\begin{theorem}\label{thm:main3}
	Suppose $F<G$ and $|G|=|F|$. If all of the following holds:
	\begin{enumerate}
		\item $\rho(v_F^G)\leq i<(\rho+\eta)(v_F^G)$;
		\item $(F,G)$ is compatible;
		\item one of the following holds:
		\begin{enumerate}
			\item $\gamma(x_{F,G})>\gamma(u_{k_{F,G}})$ and either $x_{F,G}$ is a root, or $x_{F,G}=v_L$ for some $v\in F$;
			\item $\gamma(x_{F,G})<\gamma(u_{k_{F,G}})$ and $x_{F,G}=v_R$ for some $v\in F$.
		\end{enumerate}
	\end{enumerate}
	then the structure coefficient is a positive root given by
	$$C_{F,G}^i(\bt)=\begin{cases}
		t_{(\rho+\eta)(x)}-t_{\gamma(x_{F,G})} & \gamma(x_{F,G})>\gamma(u_{k_{F,G}})\text{ and }x\text{ is a root},\\
		t_{\gamma(x)}-t_{\gamma(x_{F,G})} & \gamma(x_{F,G})>\gamma(u_{k_{F,G}})\text{ and }x\text{ is not a root},\\
		t_{\gamma(x_{F,G})}-t_{\gamma(x)} & \gamma(x_{F,G})<\gamma(u_{k_{F,G}}),
	\end{cases}$$
	where $x$ is the vertex with maximal $\eta(x)$ satisfying $\rho(x)=\rho(x_{F,G})$ in case (3)(a), and $(\rho+\eta)(x)=(\rho+\eta)(x_{F,G})$ in case (3)(b). Otherwise, $C_{F,G}^i(\bt)=0$.
	\end{theorem}
	
\begin{proof}
	By Proposition \ref{prop:ev}, condition (1) is necessary, and when it holds we have $C_{F,G}^i(\bt)=\chi_{F,G}\fP_F(\bx;\bt)$. Lemmas \ref{lem:tele1} and \ref{lem:tele3} show that this vanishes unless $(F,G)$ is compatible. For a compatible pair, Lemmas \ref{lem:tele5} and \ref{lem:tele6} give the displayed formula.
	
	If $\gamma(x_{F,G})>\gamma(u_{k_{F,G}})$, the formula is positive when $x_{F,G}$ is a root or a left child. Otherwise, $x=x_{F,G}$ is not a root, and the formula is zero. If $\gamma(x_{F,G})<\gamma(u_{k_{F,G}})$, the formula is positive precisely when $x_{F,G}$ is a right child; otherwise, $x=x_{F,G}$ and the formula is zero. This proves the result.
\end{proof}

\begin{comment}

\begin{definition}\label{def:similar}
	Let $F$ and $G$ be indexed forests with $|F|=|G|$. We say $F$ and $G$ are \textbf{similar} if there exists $v\in F$ and $u\in G$ such that
	\begin{enumerate}
		\item $F_v$ and $G_u$ have the same support;
		\item replacing $F_v$ by $G_u$ changes $F$ into $G$;
		\item $\gamma(v)\neq\gamma(u)$,
	\end{enumerate}
	where $F_v$ and $G_u$ denote the subtree of $F$ and $G$ rooted at $v$ and $u$, respectively.
\end{definition}

\begin{definition}
	Suppose $F$ and $G$ are similar. Let $v\in F$ and $u\in G$ be maximal i.e. no children of $v$ and $u$ satisfy Definition \ref{def:similar}.
	\begin{enumerate}
		\item If $\gamma(v)>\gamma(u)$, let $v'\in F$ be the minimal such that $\rho(v)=\rho(v')$. Set $p(F,G)=\gamma(v)$,
		$$q(F,G)=\begin{cases}
			(\rho+\eta)(v') & v'\text{ is the root of a tree in }F \\
			\gamma(v') & \text{otherwise}
		\end{cases}.$$
		\item If $\gamma(v)<\gamma(u)$, let $v'\in F$ be the minimal such that $(\rho+\eta)(v)=(\rho+\eta)(v')$. Set $p(F,G)=\gamma(v')$ and $q(F,G)=\gamma(v)$.
	\end{enumerate}
\end{definition}

\begin{theorem}
	Let $|F|=|G|$. If $|F|$ and $|G|$ are similar, we have $C_{F,G}^i\in\{0,t_{q(F,G)}-t_{p(F,G)}\}$; otherwise, $C_{F,G}^i=0$.
\end{theorem}

\end{comment}

\section{Monk rule and Pieri rule}\label{sec:4}

This section introduces a local move on forests that describes the admissible sequences developed in Subsection \ref{subsec:as} and the results of Theorem \ref{thm:main2} and \ref{thm:main3}. A constructive combinatorial description of the Pieri rule and the equivariant Monk rule will be made explicit using these moves. These moves are Tamari moves on certain internal nodes.

\subsection{Tamari moves}

\begin{definition}
	Let $F$ be an indexed forest, and a subtree $(A\vee B)\vee C$ of $F$, the local move $(A\vee B)\vee C\mapsto A\vee (B\vee C)$ is called a \textbf{left Tamari move}, and the root of $(B\vee C)$ in the resulting forest is called an \textbf{active node}. Similarly, for subtree $A\vee (B\vee C)$ of $F$, the local move  $A\vee (B\vee C)\mapsto (A\vee B)\vee C$ is called a \textbf{right Tamari move}, and the root of $(A\vee B)$ in the resulting forest is called an \textbf{active node}.
\end{definition}

\begin{example}
	The following are left Tamari move and right Tamari move, $u^*$ is the active node.
	\begin{center}
		\begin{tikzpicture}
			\node at (5.5,0) {\scalebox{0.8}{\begin{tikzpicture}
						\node at (0,0) {$A$};
						\node at (2,0) {$B$};
						\node at (4,0) {$C$};
						
						\draw[blue, thick] (2,0.3) -- (1,1.3);
						\draw[blue, thick] (0,0.3) -- (2,2.3) -- (4,0.3);
						
						\node at (0.6,1.4) {$u$};
						\node at (2,2.6) {$v$};
						
						\filldraw[black] (2,2.3) circle (2pt)  {};
						\filldraw[black] (1,1.3) circle (2pt)  {};
			\end{tikzpicture}}};
			\node at (8,0) {$\mapsto$};
			\node at (11,0) {\scalebox{0.8}{\begin{tikzpicture}
						\node at (0,0) {$A$};
						\node at (2,0) {$B$};
						\node at (4,0) {$C$};
						
						\draw[blue, thick] (2,0.3) -- (3,1.3);
						\draw[blue, thick] (0,0.3) -- (2,2.3) -- (4,0.3);
						
						\node at (3.4,1.4) {$u^*$};
						\node at (2,2.6) {$v$};
						
						\filldraw[black] (2,2.3) circle (2pt)  {};
						\filldraw[black] (3,1.3) circle (2pt)  {};
			\end{tikzpicture}}};
		\end{tikzpicture}
	\end{center}
	\begin{center}
		\begin{tikzpicture}
			\node at (11,0) {\scalebox{0.8}{\begin{tikzpicture}
						\node at (0,0) {$A$};
						\node at (2,0) {$B$};
						\node at (4,0) {$C$};
						
						\draw[blue, thick] (2,0.3) -- (1,1.3);
						\draw[blue, thick] (0,0.3) -- (2,2.3) -- (4,0.3);
						
						\node at (0.6,1.4) {$u^*$};
						\node at (2,2.6) {$v$};
						
						\filldraw[black] (2,2.3) circle (2pt)  {};
						\filldraw[black] (1,1.3) circle (2pt)  {};
			\end{tikzpicture}}};
			\node at (8,0) {$\mapsto$};
			\node at (5.5,0) {\scalebox{0.8}{\begin{tikzpicture}
						\node at (0,0) {$A$};
						\node at (2,0) {$B$};
						\node at (4,0) {$C$};
						
						\draw[blue, thick] (2,0.3) -- (3,1.3);
						\draw[blue, thick] (0,0.3) -- (2,2.3) -- (4,0.3);
						
						\node at (3.4,1.4) {$u$};
						\node at (2,2.6) {$v$};
						
						\filldraw[black] (2,2.3) circle (2pt)  {};
						\filldraw[black] (3,1.3) circle (2pt)  {};
			\end{tikzpicture}}};
		\end{tikzpicture}
	\end{center}
\end{example}

\begin{definition}\label{def:tamari}
	Given a pair of subtrees $(A,B)$ in an indexed forest such that the leaves of $A$ are indexed by $[a,b]$ and those of $B$ are indexed by $[b+1,c]$. Let $C$ be a tree whose leaves are indexed by $[a,c]$. We say $(A,B)\lessdot_T C$ with \textbf{interval} $[r,s)$ if $C$ can be obtained from $(A,B)$ as follows.
	
	Let the root of $A\vee B$ be the active node. Denote the active node by $v$, and perform one of the following choices.
	\begin{enumerate}
		\item Stop the process, and record the interval $[\rho(v),(\rho+\eta)(v))$.
		\item Choose a node $u$ in the subtree rooted at $v$ such that $\rho(v)=\rho(u)$, and perform a left Tamari move on $(R\vee S)\vee T$ in which $u$ is the root. Let the active node from the Tamari move be the new active node $v$, and perform one of the choices again.
		\item Choose a node $u$ in the subtree rooted at $v$ such that $(\rho+\eta)(v)=(\rho+\eta)(u)$, and perform a right Tamari move on $R\vee(S\vee T)$ in which $u$ is the root. Let the active node from the Tamari move be the new active node, and perform one of the choices again.
	\end{enumerate}
\end{definition}

\begin{example}\label{ex:m}
	Let $A=\begin{tikzpicture}[scale=0.25,baseline=0pt]
				\draw[blue, thick] (-10,0) -- (-9,1) -- (-8,0);
				
				\filldraw[black] (-9,1) circle (5pt)  {};
				
				\filldraw[black] (-10,0) circle (5pt)  {};
				\filldraw[black] (-8,0) circle (5pt)  {};
				
				\node at (-10,-1.2) {1};
				\node at (-8,-1.2) {2};
				
				\draw[black] (-12,0) -- (-6,0);
			\end{tikzpicture}$, $B=\begin{tikzpicture}[scale=0.25,baseline=0pt]
			\draw[blue, thick] (-10,0) -- (-9,1) -- (-8,0);
			
			\filldraw[black] (-9,1) circle (5pt)  {};
			
			\filldraw[black] (-10,0) circle (5pt)  {};
			\filldraw[black] (-8,0) circle (5pt)  {};
			
			\node at (-10,-1.2) {3};
			\node at (-8,-1.2) {4};
			
			\draw[black] (-12,0) -- (-6,0);
			\end{tikzpicture}$. Then the following forests $C$ satisfy $(A,B)\lessdot C$, with final active node circled in red, and corresponding interval written below.
			\begin{center}
				\begin{tikzpicture}
					\node at (0,0) {\begin{tikzpicture}[scale=0.25,baseline=0pt]
							\draw[blue, thick] (-10,0) -- (-7,3) -- (-4,0);
							\draw[blue, thick] (-8,0) -- (-9,1);
							\draw[blue, thick] (-6,0) -- (-5,1);
							
							\draw[red] (-7,3) circle (10pt);
							\filldraw[red] (-7,3) circle (5pt)  {};
							\filldraw[black] (-9,1) circle (5pt)  {};
							\filldraw[black] (-5,1) circle (5pt)  {};
							
							\filldraw[black] (-10,0) circle (5pt)  {};
							\filldraw[black] (-8,0) circle (5pt)  {};
							\filldraw[black] (-6,0) circle (5pt)  {};
							\filldraw[black] (-4,0) circle (5pt)  {};
							
							\node at (-10,-1.2) {1};
							\node at (-8,-1.2) {2};
							\node at (-6,-1.2) {3};
							\node at (-4,-1.2) {4};
							
							\draw[black] (-12,0) -- (-2,0);
							
							\node at (-7,-3.5) {$[1,4)$};
							\node at (-7,-5.5) {$C_1$};
						\end{tikzpicture}};
					\node at (3,0) {\begin{tikzpicture}[scale=0.25,baseline=0pt]
							\draw[blue, thick] (-10,0) -- (-7,3) -- (-4,0);
							\draw[blue, thick] (-8,0) -- (-6,2);
							\draw[blue, thick] (-6,0) -- (-5,1);
							
							\filldraw[black] (-7,3) circle (5pt)  {};
							\draw[red] (-6,2) circle (10pt);
							\filldraw[red] (-6,2) circle (5pt)  {};
							\filldraw[black] (-5,1) circle (5pt)  {};
							
							\filldraw[black] (-10,0) circle (5pt)  {};
							\filldraw[black] (-8,0) circle (5pt)  {};
							\filldraw[black] (-6,0) circle (5pt)  {};
							\filldraw[black] (-4,0) circle (5pt)  {};
							
							\node at (-10,-1.2) {1};
							\node at (-8,-1.2) {2};
							\node at (-6,-1.2) {3};
							\node at (-4,-1.2) {4};
							
							\draw[black] (-12,0) -- (-2,0);
							
							\node at (-7,-3.5) {$[2,4)$};
							\node at (-7,-5.5) {$C_2$};
					\end{tikzpicture}};
					\node at (6,0) {\begin{tikzpicture}[scale=0.25,baseline=0pt]
							\draw[blue, thick] (-10,0) -- (-7,3) -- (-4,0);
							\draw[blue, thick] (-8,0) -- (-6,2);
							\draw[blue, thick] (-6,0) -- (-7,1);
							
							\filldraw[black] (-7,3) circle (5pt)  {};
							\draw[red] (-7,1) circle (10pt);
							\filldraw[black] (-6,2) circle (5pt)  {};
							\filldraw[red] (-7,1) circle (5pt)  {};
							
							\filldraw[black] (-10,0) circle (5pt)  {};
							\filldraw[black] (-8,0) circle (5pt)  {};
							\filldraw[black] (-6,0) circle (5pt)  {};
							\filldraw[black] (-4,0) circle (5pt)  {};
							
							\node at (-10,-1.2) {1};
							\node at (-8,-1.2) {2};
							\node at (-6,-1.2) {3};
							\node at (-4,-1.2) {4};
							
							\draw[black] (-12,0) -- (-2,0);
							
							\node at (-7,-3.5) {$[2,3)$};
							\node at (-7,-5.5) {$C_3$};
					\end{tikzpicture}};
					\node at (9,0) {\begin{tikzpicture}[scale=0.25,baseline=0pt]
							\draw[blue, thick] (-10,0) -- (-7,3) -- (-4,0);
							\draw[blue, thick] (-8,0) -- (-9,1);
							\draw[blue, thick] (-6,0) -- (-8,2);
							
							\filldraw[black] (-7,3) circle (5pt)  {};
							\filldraw[black] (-9,1) circle (5pt)  {};
							\draw[red] (-8,2) circle (10pt);
							\filldraw[red] (-8,2) circle (5pt)  {};
							
							\filldraw[black] (-10,0) circle (5pt)  {};
							\filldraw[black] (-8,0) circle (5pt)  {};
							\filldraw[black] (-6,0) circle (5pt)  {};
							\filldraw[black] (-4,0) circle (5pt)  {};
							
							\node at (-10,-1.2) {1};
							\node at (-8,-1.2) {2};
							\node at (-6,-1.2) {3};
							\node at (-4,-1.2) {4};
							
							\draw[black] (-12,0) -- (-2,0);
							
							\node at (-7,-3.5) {$[1,3)$};
							\node at (-7,-5.5) {$C_4$};
					\end{tikzpicture}};
					\node at (12,0) {\begin{tikzpicture}[scale=0.25,baseline=0pt]
							\draw[blue, thick] (-10,0) -- (-7,3) -- (-4,0);
							\draw[blue, thick] (-8,0) -- (-7,1);
							\draw[blue, thick] (-6,0) -- (-8,2);
							
							\filldraw[black] (-7,3) circle (5pt)  {};
							\draw[red] (-7,1) circle (10pt);
							\filldraw[red] (-7,1) circle (5pt)  {};
							\filldraw[black] (-8,2) circle (5pt)  {};
							
							\filldraw[black] (-10,0) circle (5pt)  {};
							\filldraw[black] (-8,0) circle (5pt)  {};
							\filldraw[black] (-6,0) circle (5pt)  {};
							\filldraw[black] (-4,0) circle (5pt)  {};
							
							\node at (-10,-1.2) {1};
							\node at (-8,-1.2) {2};
							\node at (-6,-1.2) {3};
							\node at (-4,-1.2) {4};
							
							\draw[black] (-12,0) -- (-2,0);
							
							\node at (-7,-3.5) {$[2,3)$};
							\node at (-7,-5.5) {$C_5$};
					\end{tikzpicture}};
				\end{tikzpicture}
			\end{center}
\end{example}

\subsection{The ordinary Pieri rule}

This subsection establish the combinatorial formula for Monk and Pieri rule of ordinary forest polynomials based on $\lessdot_i$. The technical part in the proof of the Monk rule will be deferred to subsection \ref{subsec:proof}

\begin{definition}
	Given indexed forests $F$ and $G$ with $|F|+1=|G|$, we say $F\lessdot G$ with \textbf{interval} $[a,b)$ if $G$ is obtained from $F$ by replacing two consecutive components $A,B$ of $F$ by $C$ with $(A,B)\lessdot_T C$ and corresponding interval is $[a,b)$.
\end{definition}

We remark that the transitive closure of $\lessdot$ gives a partial order on indexed forests.

\begin{definition}
	Given $i\geq 1$ and indexed forests $F$ and $G$, we say $F\lessdot_i G$ if $F\lessdot G$ with interval $[a,b)$, and additionally, $a\leq i<b$.
\end{definition}

\begin{theorem}\label{thm:m}
	The Monk rule for ordinary forest polynomials have the following formula.
	$$\fP_F(\bx)(x_1+\cdots+x_i)=\sum_{F\lessdot_i G}\fP_G(\bx).$$
\end{theorem}

\begin{proof}
	By Proposition~\ref{prop:tamari-admissible}(1), that will be proved later, and Theorem \ref{thm:main2}, the coefficient with $|G|=|F|+1$ is $1$ precisely when $F\lessdot_i G$. Setting $\bt=0$ in Corollary \ref{cor:em-1} removes all
	coefficients of degree one, and the result follows.
\end{proof}

\begin{example}
	We follow example \ref{ex:m} Let $F=A\uplus B$, then $\fP_F(\bx)\cdot x_1=\fP_{C_1}+\fP_{C_4}$ while $\fP_F(\bx)\cdot (x_1+x_2)=\sum_{i=1}^5\fP_{C_i}(\bx)$.
\end{example}

To obtain the Pieri rule, we need the following preparation.

\begin{lemma}\label{lem:p1}
	Let $F$ and $G$ be indexed forests.
	\begin{enumerate}
		\item Suppose $F\lessdot G$ with interval $[a,b)$, and $j\in\LTer(G)\setminus\LTer(F)$, then $[a,b)=[j,j+1)$.
		\item Suppose $F\lessdot G$ with interval $[a,b)$, and $j\in\LTer(F)\cap\LTer(G)$, then $F/j\lessdot G/j$. Moreover, if $F/j\lessdot G/j$ has interval $[a',b')$, then $a=a'+\mathbf{1}_{a'>j}$ and $b=b'+\mathbf{1}_{b'\geq j}$.
		\item Suppose $j\in\LTer(F)\cap\LTer(G)$ and $F/j\lessdot G/j$, then $F\lessdot G$.
	\end{enumerate}
\end{lemma}

\begin{proof}
	(1) If $j\in\LTer(G)\setminus\LTer(F)$, then the root of $j\vee (j+1)$ in $G$ is the final active node. By construction, $[a,b)=[j,j+1)$.
	
	(2) If $j\in\LTer(F)\cap\LTer(G)$, then the root of $j\vee (j+1)$ is not the active node. Therefore, $F/j\lessdot G/j$, and the sequence of choices in the construction as Definition \ref{def:tamari} is preserved. Trimming $j$ changes all indexes $r\geq j$ to $r-1$, consequently, the interval changes accordingly.
	
	(3) Conversely, if $F/j\lessdot G/j$, blossoming at $j$ preserves the sequence of choices in the construction as Definition \ref{def:tamari}, and hence, $F\lessdot G$.
\end{proof}

\begin{definition}
	Fix $i\geq 1$. For $1\leq a\leq i$, we define the linear operators on $\mathbb{Z}[\bx]$ by
	$$L_a^i\fP_F=\sum_{F\lessdot_i G \text{ with interval }[a,b) \text{ for some }b}\fP_G.$$
	For $b>i$, we define the linear operators on $\mathbb{Z}[\bx]$ by
	$$B_b^i\fP_F=\sum_{F\lessdot_i G\text{ with interval }[a,b) \text{ for some }a}\fP_G.$$
\end{definition}

It follows from Theorem \ref{thm:m} and linearity that for any $f\in\mathbb{Z}[\bx]$,
$$f(\bx)\fS_{s_i}(\bx)=\sum_{a=1}^i L_a^i(f)=\sum_{b>i}B_b^i(f).$$
Recall the trimming operators $R_i$ and $T_i$ from Definition \ref{def:trim}.

\begin{lemma}\label{lem:p2}
	Fix $1\leq a\leq i$, for $j>i$, we have $T_jL_a^i=L_a^iT_j$. For $j<i$, we have
	$$T_jL_a^i=\begin{cases}
		L_a^{i-1}T_j & a\leq j \\
		0 & a = j+1 \\
		L_{a-1}^{i-1}T_j & a > j+1
	\end{cases}.$$
	For $j=i$, we have
	$$T_iL_a^i=\begin{cases}
		L_a^{i-1}T_i & a<i \\
		R_{i+1} & a=i
	\end{cases}.$$
\end{lemma}

\begin{proof}
	Suppose $a=i=j$ does not hold. Note that we have
	$$T_jL_a^i\fP_F(\bx)=\sum_{F\lessdot_i G \text{ with interval }[a,b) \text{ for some }b \atop j\in\LTer(G)}\fP_{G/j}(\bx).$$
	If $j\notin \LTer(F)$, then $T_j\fP_F(\bx)=0$ and Lemma \ref{lem:p1} implies the interval of $F\lessdot G$ is $[j,j+1)$. Hence, $T_jL_a^i=0$. If $j\in\LTer(F)$, then Lemma \ref{lem:p1} implies $F\lessdot G$ if and only if $F/j\lessdot G/j$, and the interval changes accordingly.
	
	We now prove the case $a=i=j$. For any $f(\bx)$, note that we have
	\begin{align*}
		T_i(\fS_{s_i}(\bx)f(\bx))&=\frac{\fS_{s_i}(\bx)R_{i+1}f(\bx)-\fS_{s_{i-1}}(\bx)R_if(\bx)}{x_i}\\
		&=\fS_{s_{i-1}}(\bx)\frac{R_{i+1}f(\bx)-R_if(\bx)}{x_i}+R_{i+1}f(\bx)=\fS_{s_{i-1}}(\bx)T_if(\bx)+R_{i+1}f(\bx).
	\end{align*}
	Using the results from the previous cases, the statement follows as
	\begin{align*}
		T_iL_i^if(\bx)&=T_i(\fS_{s_i}(\bx)f(\bx))-\sum_{k=1}^{i-1}T_iL_k^if(\bx)=T_i(\fS_{s_i}(\bx)f(\bx))-\sum_{k=1}^{i-1}L_k^{i-1}T_if(\bx)\\
		&=T_i(\fS_{s_i}(\bx)f(\bx))-\fS_{s_{i-1}}(\bx)T_if(\bx)=R_{i+1}f(\bx).
	\end{align*}
\end{proof}

\begin{lemma}\label{lem:p3}
	Fix $i<b$, for $j<i$, we have $T_jB_b^i=B_{b-1}^{i-1}T_j$. For $j>i$, we have
	$$T_jB_b^i=\begin{cases}
		B_b^iT_j & b<j \\
		0 & b=j \\
		B_{b-1}^iT_j & b>j
	\end{cases}.$$
	For $j=i$, we have
	$$T_iB_b^i=\begin{cases}
		B_{b-1}^iT_i & b> i+1 \\
		R_i & b=i+1
	\end{cases}.$$
\end{lemma}

\begin{proof}
	The argument is parallel to Lemma \ref{lem:p2} so we omit it.
\end{proof}

\begin{definition}
	We define the generating functions
	$$\mathcal{H}_i(z)=\cdots (\id+zB_{i+2}^i)(\id+zB_{i+1}^i),$$
	$$\mathcal{E}_i(z)=(\id+zL_1^i)(\id+zL_2^i)\cdots (\id+zL_i^i),$$
	with the convention that $\mathcal{E}_0(z)=\mathcal{H}_0(z)=\id$.
\end{definition}

\begin{lemma}\label{lem:p4}
	For any polynomial $f(\bx)$ and $i\geq0$, we have
	$$\mathcal{H}_i(z)f(\bx)=\prod_{q=1}^i(1-zx_q)^{-1}f(\bx),$$
	$$\mathcal{E}_i(z)f(\bx)=\prod_{q=1}^i(1+zx_q)f(\bx).$$
\end{lemma}

\begin{proof}
	By Lemma \ref{lem:p2} and \ref{lem:p3}, we have
	$$T_j\mathcal{H}_i(z)=\begin{cases}
		\mathcal{H}_{i-1}(z)T_j & j<i \\
		\mathcal{H}_i(z)(T_i+zR_i) & j=i \\
		\mathcal{H}_i(z)T_j & j>i
	\end{cases},$$
	$$T_j\mathcal{E}_i(z)=\begin{cases}
		\mathcal{E}_{i-1}(z)T_j & j<i \\
		\mathcal{E}_{i-1}(z)(T_i+zR_{i+1}) & j=i \\
		\mathcal{E}_i(z)T_j & j>i
	\end{cases}.$$
	On the other hand, we also have
	$$T_j\left(\prod_{q=1}^i(1-zx_q)^{-1}f(\bx)\right)=\begin{cases}
		\prod_{1\leq q< i}(1-zx_q)^{-1}T_jf(\bx) & j<i \\
		\prod_{1\leq q\leq i}(1-zx_q)^{-1}\left(T_jf(\bx)+zR_if(\bx)\right) & j=i \\
		\prod_{1\leq q\leq i}(1-zx_q)^{-1}T_jf(\bx) & j>i
	\end{cases},$$
	$$T_j\left(\prod_{q=1}^i(1+zx_q)f(\bx)\right)=\begin{cases}
		\prod_{1\leq q <i}(1+zx_q)T_jf(\bx) & j<i \\
		\prod_{1\leq q<i}(1+zx_q)\left(T_jf(\bx)+zR_{i+1}f(\bx)\right) & j=i \\
		\prod_{1\leq q\leq i}(1+zx_q)T_jf(\bx) & j>i
	\end{cases}.$$
\end{proof}

\begin{theorem}\label{thm:pieri}
	The Pieri's rule for ordinary forest polynomials have the following formula.
	$$\fP_F(\bx)h_r(x_1,\dots,x_i)=\sum_{F=F_0\lessdot_i F_1\lessdot_i\cdots\lessdot_i F_r=G \atop a_1\leq\cdots\leq a_r}\fP_G(\bx)=\sum_{F=F_0\lessdot_i F_1\lessdot_i\cdots\lessdot_i F_r=G \atop b_1<\cdots< b_r}\fP_G(\bx)$$
	where $[a_k,b_k)$ is the interval of $F_{k-1}\lessdot_i F_k$, $1\leq k\leq r$. Moreover,
	$$\fP_F(\bx)e_r(x_1,\dots,x_i)=\sum_{F=F_0\lessdot_i F_1\lessdot_i\cdots\lessdot_i F_r=G \atop a_1>\cdots> a_r}\fP_G(\bx)=\sum_{F=F_0\lessdot_i F_1\lessdot_i\cdots\lessdot_i F_r=G \atop b_1\geq\cdots\geq b_r}\fP_G(\bx).$$
\end{theorem}

\begin{proof}
	The coefficients of $z^r$ in $\mathcal{H}_i(z)\fP_F(\bx)$ and $\mathcal{E}_i(z)\fP_F(\bx)$ are $\fP_F(\bx)h_r(x_1,\dots,x_i)$ and $\fP_F(\bx)e_r(x_1,\dots,x_i)$, respectively, by Lemma \ref{lem:p4}. On the other hand, we have
	$$\mathcal{H}_i(z)\fP_F(\bx)=\cdots (\id+zB_{i+2}^i)(\id+zB_{i+1}^i)\fP_F(\bx)=\sum_{r\geq 0}\left(\sum_{F=F_0\lessdot_i F_1\lessdot_i\cdots\lessdot_i F_r=G \atop b_1<\cdots< b_r}\fP_G(\bx)\right)z^r,$$
	$$\mathcal{E}_i(z)\fP_F(\bx)=(\id+zL_1^i)(\id+zL_2^i)\cdots (\id+zL_i^i)\fP_F(\bx)=\sum_{r\geq 0}\left(\sum_{F=F_0\lessdot_i F_1\lessdot_i\cdots\lessdot_i F_r=G \atop a_1>\cdots> a_r}\fP_G(\bx)\right)z^r.$$
	Consider $\widetilde{\mathcal{H}}_i(z)=\mathcal{H}_i(-z)^{-1}$ and $\widetilde{\mathcal{E}}_i(z)=\mathcal{E}_i(-z)^{-1}$. By Lemma \ref{lem:p4}, $\widetilde{\mathcal{H}}_i(z)\fP_F(\bx)=\mathcal{E}_i(z)\fP_F(\bx)$ and $\widetilde{\mathcal{E}}_i(z)\fP_F(\bx)=\mathcal{H}_i(z)\fP_F(\bx)$. Hence, the coefficients of $z^r$ in $\widetilde{\mathcal{H}}_i(z)\fP_F(\bx)$ and $\widetilde{\mathcal{E}}_i(z)\fP_F(\bx)$ are $\fP_F(\bx)e_r(x_1,\dots,x_i)$ and $\fP_F(\bx)h_r(x_1,\dots,x_i)$, respectively. On the other hand, we have
	$$\widetilde{\mathcal{H}}_i(z)\fP_F(\bx)=\left[(\id-zB_{i+1}^i)^{-1}(\id-zB_{i+2}^i)^{-1}\cdots\right]\fP_F(\bx)=\sum_{r\geq 0}\left(\sum_{F=F_0\lessdot_i F_1\lessdot_i\cdots\lessdot_i F_r=G \atop b_1\geq\cdots\geq b_r}\fP_G(\bx)\right)z^r,$$
	$$\widetilde{\mathcal{E}}_i(z)\fP_F(\bx)=(\id-zL_i^i)^{-1}\cdots (\id-zL_1^i)^{-1}\fP_F(\bx)=\sum_{r\geq 0}\left(\sum_{F=F_0\lessdot_i F_1\lessdot_i\cdots\lessdot_i F_r=G \atop a_1\leq\cdots\leq a_r}\fP_G(\bx)\right)z^r.$$
\end{proof}

\begin{example}\label{ex:pieri-illustrated}
	We retain the notation of Example~\ref{ex:m} and let $F=A\uplus B$. The forests $F,D_1,\ldots,D_7$ are shown below. In particular, $D_j=C_j\vee5$ for $1\leq j\leq5$, where $5$ denotes the leaf indexed by $5$.
	\begin{center}
		\begin{tabular}{@{}cccc@{}}
			\FPNonfreeForest{1}{F} &
			\FPNonfreeForest{2}{D_1} &
			\FPNonfreeForest{3}{D_2} &
			\FPNonfreeForest{4}{D_3} \\[0.9ex]
			\FPNonfreeForest{5}{D_4} &
			\FPNonfreeForest{6}{D_5} &
			\FPNonfreeForest{7}{D_6} &
			\FPNonfreeForest{8}{D_7}
		\end{tabular}
	\end{center}
	For $i=r=2$, the Pieri rule gives
	\[
	\begin{aligned}
		\fP_F(\bx)h_2(x_1,x_2)
		&={}\fP_{D_1}(\bx)+\fP_{D_2}(\bx)
		+2\fP_{D_3}(\bx)+\fP_{D_4}(\bx)\\
		&\quad+\fP_{D_5}(\bx)+\fP_{D_6}(\bx)+\fP_{D_7}(\bx),\\
		\fP_F(\bx)e_2(x_1,x_2)
		&={}\fP_{D_2}(\bx)+\fP_{D_3}(\bx)+\fP_{D_5}(\bx).
	\end{aligned}
	\]
	The coefficient $2$ of $\fP_{D_3}(\bx)$ in the first expansion
	comes from the two distinct chains
	\[
	\begin{array}{c@{\qquad}c@{\qquad}c}
		\text{chain} & [a_1,b_1) & [a_2,b_2) \\ \hline
		F\lessdot_2 C_1\lessdot_2 D_3 & [1,4) & [2,3) \\
		F\lessdot_2 C_5\lessdot_2 D_3 & [2,3) & [2,4)
	\end{array}
	\]
	Both satisfy $a_1\leq a_2$. Thus the Pieri expansion need not
	be multiplicity-free.
\end{example}

\subsection{The equivariant Monk rule}

This subsection establish the combinatorial formula for the equivariant Monk rule of double forest polynomials based on $\prec_i$. The technical part in the proof of the equivariant Monk rule will be deferred to subsection \ref{subsec:proof}

\begin{definition}
	Given indexed forests $F$ and $G$ with $|F|=|G|$, we say $F\prec G$ with \textbf{interval} $[a,b)$ and \textbf{weight} $\mathrm{wt}_{F,G}$ if $F\neq G$ and there exists a node $v\in F$ such that one of the following holds.
	\begin{enumerate}
		\item $v$ is the root of $(A\vee B)\vee C$ and either $v$ is the root of a component, or $v$ is a left child. $G$ is obtained from $F$ by replacing the subtree $(A\vee B)\vee C$ of $F$ by $A\vee U$ with $(B,C)\lessdot_T U$ and corresponding interval is $[a,b)$. Let $x\in F$ be the highest ancestor of $v$, possibly $x=v$, satisfying $\rho(x)=\rho(v)$, then the weight is
		$$\mathrm{wt}_{F,G}=\begin{cases}
			t_{(\rho+\eta)(x)}-t_{\gamma(v)} & x\text{ is a component root} \\
			t_{\gamma(x)} - t_{\gamma(v)} & x\text{ is a right child}
		\end{cases}.$$
		\item $v$ is the root of $A\vee (B\vee C)$ and $v$ is a right child. $G$ is obtained from $F$ by replacing the subtree $A\vee (B\vee C)$ of $F$ by $U\vee C$ with $(A,B)\lessdot_T U$ and corresponding interval is $[a,b)$. Let $x\in F$ be the highest ancestor of $v$, possibly $x=v$, satisfying $(\rho+\eta)(x)=(\rho+\eta)(v)$, then the weight is $\mathrm{wt}_{F,G}=t_{\gamma(v)}-t_{\gamma(x)}$.
	\end{enumerate}
\end{definition}

We emphasize that $\prec$ does not generate a partial order on indexed forests, since there possibly are cyclic relations. By abuse of notation, we sometimes say $G$ covers $F$ for $F\prec G$.

\begin{definition}
	Given $i\geq 1$ and indexed forests $F$ and $G$, we say $F\prec_i G$ with \textbf{weight} $\mathrm{wt}_{F,G}$ if $F\prec G$ with interval $[a,b)$ and weight $\mathrm{wt}_{F,G}$, and additionally, $a\leq i<b$.
\end{definition}

\begin{theorem}\label{thm:emonk}
	The equivariant Monk rule for double forest polynomials have the following formula. We emphasize that each of the coefficients $\mathrm{wt}_{F,G}$ is always a single positive root.
	$$\fP_F(\bx;\bt)\fS_{s_i}(\bx;\bt)=\left(\sum_{k=1}^i(t_{w_F(k)}-t_k)\fP_F(\bx;\bt)\right)+\left(\sum_{F\prec_i G}\mathrm{wt}_{F,G}\fP_G(\bx;\bt)\right)+\left(\sum_{F\lessdot_i G}\fP_G(\bx;\bt)\right).$$
\end{theorem}

\begin{proof}
	The diagonal coefficient is given by Theorem \ref{thm:main1}. By Proposition~\ref{prop:tamari-admissible}, that will be proved later, and Theorems \ref{thm:main2} and \ref{thm:main3}, the other nonzero coefficients are $1$ for $F\lessdot_i G$ and $\operatorname{wt}_{F,G}$ for $F\prec_i G$, respectively. Substituting these coefficients into Corollary \ref{cor:em-1} gives the formula.
\end{proof}

\begin{example}\label{ex:equivariant-monk-illustrated}
	We use the forests $T,U,C_4,C_5$ from
	Example~\ref{ex:pieri-illustrated}. The move $T\prec U$ has
	interval $[2,3)$ and weight $t_3-t_2$:
	\begin{center}
		\begin{tabular}{@{}ccc@{}}
			\FPExampleForest{2}{}{T} &
			$\xrightarrow[\mathrm{wt}_{T,U}=t_3-t_2]{[2,3)}$ &
			\FPExampleForest{3}{\FPExampleActive{3,1}}{U}
		\end{tabular}
	\end{center}
	The covers $T\lessdot C_4$ and $T\lessdot C_5$ have intervals $[1,4)$ and $[2,3)$, respectively. Their final active nodes are circled in red.
	\begin{center}
		\begin{tabular}{@{}c@{\qquad\qquad}c@{}}
			\FPExampleForest{5}{\FPExampleActive{3,3}}{C_4\quad[1,4)} &
			\FPExampleForest{6}{\FPExampleActive{3,1}}{C_5\quad[2,3)}
		\end{tabular}
	\end{center}
	Since $w_T(1)=3$ and $w_T(2)=1$, we obtain
	\[
	\fP_T(\bx;\bt)\fS_{s_1}(\bx;\bt)
	=(t_3-t_1)\fP_T(\bx;\bt)+\fP_{C_4}(\bx;\bt),
	\]
	while
	\[
	\begin{aligned}
		\fP_T(\bx;\bt)\fS_{s_2}(\bx;\bt)
		=(t_3-t_2)\fP_T(\bx;\bt)+(t_3-t_2)\fP_U(\bx;\bt)+\fP_{C_4}(\bx;\bt)+\fP_{C_5}(\bx;\bt).
	\end{aligned}
	\]
\end{example}

\subsection{The correspondence}\label{subsec:proof}

This subsection gives the correspondence between Tamari moves and admissible sequences.

\begin{proposition}\label{prop:tamari-admissible}
	Let $F,G$ be distinct indexed forests. When $F<G$, use the notation $u_0=v_F^G,u_1,\ldots,u_m$ and the constructions in Definitions \ref{def:e1} and \ref{def:e2}.
	\begin{enumerate}
		\item Suppose $|G|=|F|+1$. Then $F\lessdot G$ if and only if $F<G$ and $F=K^\epsilon$ for an admissible sequence $\epsilon=(\epsilon_0,\ldots,\epsilon_m)$. In this case, the two components joined in the definition of $F\lessdot G$ are $A_m$ and $B_m$.
		\item Suppose $|G|=|F|$. Then $F\prec G$ if and only if $F<G$
		and there exist $1\le k\le m$ and an admissible sequence
		$\epsilon=(\epsilon_0,\ldots,\epsilon_k)$ such that $F=H_{k+1}^\epsilon$, and one of the following holds:
		\begin{enumerate}
			\item $u_{k-1}=(u_k)_R$, $\epsilon_k=+$, and
			$w_k^\epsilon$ is a component root or a left child in $F$;
			\item $u_{k-1}=(u_k)_L$, $\epsilon_k=-$, and
			$w_k^\epsilon$ is a right child in $F$.
		\end{enumerate}
		Here $H_{m+1}^\epsilon=H^\epsilon$. In this case, $k$ and
		$\epsilon$ are unique, and the node $v$ in the definition of
		$F\prec G$ satisfies
		\[
		v=w_k^\epsilon=x_{F,G},\qquad k=k_{F,G}.
		\]
		Moreover, $\operatorname{wt}_{F,G}$ is the positive root displayed
		in Theorem \ref{thm:main3}.
	\end{enumerate}
	In both cases, the final active node is $v_F^G$, and the interval is
	\[
	[\rho(v_F^G),(\rho+\eta)(v_F^G)).
	\]
	Consequently, $F\lessdot_i G$ and $F\prec_i G$ correspond exactly
	to the nonvanishing conditions in Theorems \ref{thm:main2} and \ref{thm:main3}, respectively.
\end{proposition}

\begin{proof}
	Given $F\lessdot G$ or $F\prec G$, we first show the final active node is $v_F^G$. Recall from Definition \ref{def:h} that $v_F^G=v_{h(F,G)}^G$ and $h(F,G)=\min\{r:S_r\leq_R F:\Trim(S_r)=(j_{r+1},\dots,j_n)\}$, where $(j_1,\dots,j_n)$ is the unique non-increasing trimming sequence of $G$. In the other words, $h(F,G)$ corresponds to the highest node $v\in G$ such that $u\in G$ is not an ancestor of $v$ implies the subtree, with indexed leaves, of $G$ rooted at $u$ is identical to the subtree, with indexed leaves of $F$ rooted at some $u'\in F$. By construction of $\lessdot$ and $\prec$, $v_F^G$ is exactly the final active node. Consequently, $u_0=v_F^G$, which proves the assertion about the interval.
	
	Moreover, tracing the sequence of active nodes in the construction of $\lessdot_T$ backwards from $v_F^G$, it gives a sequence of nodes $u_0,u_1,\dots,u_k$ where $u_i=(u_{i+1})_L$ or $u_i=(u_{i+1})_R$, together with corresponding admissible sequence $\epsilon_0,\epsilon_1,\dots,\epsilon_k$. We have $\epsilon_i=-\epsilon_{i+1}$ if and only if a Tamari move is performed at $u_{i+1}$. Therefore, there is a one-to-one correspondence between admissible sequences and covers $\lessdot$ or $\prec$.
	
	It remains to verify $F<G$ for the pairs constructed by the moves. The coefficient-extraction proof of Lemma \ref{lem:ev3} and the expansions in Lemmas \ref{lem:ev6} and \ref{lem:ev7} apply to an arbitrary active node, without assuming $F<G$. Take $i=\gamma(u_0)$ and expand every active node term. In (1), the recovered sequence contributes $1$ to the term corresponding to $u_0$, and all other terms of degree zero are non-negative. In (2), the recovered sequence contributes the positive term
	\[
	\begin{cases}
		t_{(\rho+\eta)(v)}-t_{\gamma(v)},
		& \text{in (a), if }v\text{ is a component root},\\
		t_{\gamma(y)}-t_{\gamma(v)},
		& \text{in (a), if }v=y_L,\\
		t_{\gamma(v)}-t_{\gamma(y)},
		& \text{in (b), if }v=y_R.
	\end{cases}
	\]
	These are the final term in Lemma \ref{lem:ev7} or the term at $p=k+1$ in Lemma \ref{lem:ev6}. All other terms of degree one are non-negative roots. Thus the coefficient of $\mathfrak P_G(x;t)$ is nonzero in both cases, and Corollary \ref{cor:em-1} implies $F<G$. 
	
	Finally, in (2), the two cases are precisely the compatibility conditions in Definition \ref{def:tele} together with condition (3) of Theorem \ref{thm:main3}. Thus $v=x_{F,G}$ and $k=k_{F,G}$, as desired.
\end{proof}

\section{Another Monk rule via slide polynomial}\label{sec:m1}

%\begin{comment}

This section introduces an insertion algorithm that leads to a Monk rule for certain slide polynomials. We then give another Monk rule for forest polynomials by (1) expanding forest polynomials in slide polynomials via Nadeau-Tewari expansion (2) applying our Monk rule for slide polynomials, and (3) grouping the slide polynomial expansion back into forest polynomials.

\begin{definition}
	An \textbf{injective word} is a word $w=w_1\cdots w_\ell\in\mathbb{Z}^\ell$ such that all $w_k$ are distinct, $w$ is a \textbf{positive injective word} if all $w_k\in\mathbb{Z}_{>0}$. We also introduce the \textbf{decorated letters} $\mathbb{Z}\cup \{a^{[b]}:a\in\mathbb{Z},b\in\{-1,1\}\}$, totally ordered by $(a-1)^{[1]}<a^{[-1]}<a<a^{[1]}$. The \textbf{value} of $a^{[b]}$ is $\mathrm{val}(a^{[b]})=a$.
\end{definition}

\begin{lemma}\label{lem:m1-1}
	For every positive injective word $w$, there exists a unique indexed forest $F$ with $w\in\mathrm{Syl}(F)$.
\end{lemma}

\begin{proof}
	Decompose the support $\mathrm{supp}(w)$ into maximal consecutive intervals. Each consecutive interval $J$ gives a Sylvester word of the binary tree as follows. The smallest $j$ such that $w_j\in J$ corresponds to a root $v$, the subword of $w$ in $[1,w_j)\cap J$ is the left subtree of $v$ while the subword in $(w_j,\infty)\cap J$ is the right subtree of $v$, and proceed by induction.
\end{proof}

\begin{definition}
	Given a word $w=w_1\cdots w_\ell$, a word $b_1\cdots b_\ell\in\mathbb{Z}_{>0}^\ell$ is \textbf{$w$-compatible} if $b_1\leq\cdots\leq b_\ell$, $b_k\leq w_k$, and $w_k\leq w_{k+1}$ implies $b_k<b_{k+1}$. The set of $w$-compatible words is denoted by $\E(w)$. The \textbf{slide polynomial} indexed by $w$ is
	$$\mathfrak{F}_w=\sum_{b_1\cdots b_\ell\in\E(w)}x_{b_1}\cdots x_{b_\ell}.$$
\end{definition}

\begin{theorem}[\cite{NT} Proposition 3.8]
	We have
	$$\mathfrak{P}_F=\sum_{w\in\mathrm{Syl}(F)}\mathfrak{F}_w.$$
\end{theorem}

\begin{definition}
	A positive injective word $w$ is \textbf{legal} if every letter in the $k$-th decreasing run of $w$ is at least $k$. The set of legal Sylvester words of $F$ is denoted by $\mathrm{LSyl}(F)$.
\end{definition}

\begin{lemma}\label{lem:m1-2}
	A positive injective word $w$ is legal if and only if $\mathfrak{F}_w\neq 0$. Consequently,
	$$\mathfrak{P}_F=\sum_{w\in\mathrm{LSyl}(F)}\mathfrak{F}_w.$$
\end{lemma}

\begin{proof}
	Suppose decreasing run decomposition of $w$ is $D_1|\cdots|D_\ell$. If $w$ is legal, then $x_1^{|D_1|}\cdots x_\ell^{|D_\ell|}$ is a monomial in $\mathfrak{F}_w$. Otherwise, $\E(w)=\varnothing$ and $\mathfrak{F}_w=0$.
\end{proof}

\subsection{An insertion algorithm}

This subsection introduces an insertion algorithm, inspired by the insertion in \cite{NT}, to give a suitable Monk rule for slide polynomials.

\begin{definition}\label{def:ins}
	Given a positive injective word $w=w_1,\dots,w_\ell$,  $i\in\mathbb{Z}_{>0}$ and $1\leq p\leq \ell+1$. Set $A=\min\{k\geq i:k\notin\{w_p,\dots,w_\ell\}\}$, and define $\mathrm{Ins}_i(w,p)=w_1\cdots w_{p-1}A^{[1]}w_p\cdots w_\ell$. We then define $\Phi_i(w,p)$ as follows. Let $q$ be a value recorder and $j$ be a position recorder.
	\begin{enumerate}
		\item To initialize, set $\Phi_i(w,p)=w_1\cdots w_{p-1}A^{[1]}w_p\cdots w_\ell$, $q=A$ and $j=p$.
		\item If $w_r=q$ for some $1\leq r< j$, let $J\subset\{\mathrm{val}(\Phi_i(w,p)_{r+1}),\dots,\mathrm{val}(\Phi_i(w,p)_{\ell+1})\}$  be the maximal consecutive interval containing $q$, and let $t\geq r+1$ be the smallest number such that $\mathrm{val}(\Phi_i(w,p)_t)\in J$.
		\begin{enumerate}
			\item If $w_r>\Phi_i(w,p)_t$, then set $\Phi_i(w,p)_r=(\max J+1)^{[-1]}$.
			\item If $w_r<\Phi_i(w,p)_t$, then set $\Phi_i(w,p)_r=(\min J-1)^{[1]}$.
		\end{enumerate}
		\item Set $q$ to be the new $\mathrm{val}(\Phi_i(w,p)_r)$ and set $j=r$. Repeat step (2) and (3) until $w_r$ does not exist.
		\item Finally, remove all decorations.
	\end{enumerate}
	Note that $\Phi_i(w,p)$ is always injective, but not necessarily positive.
\end{definition}

We emphasize that there is a unique $q\notin\mathrm{supp}(w)$ such that $\mathrm{supp}(\Phi_i(w,p))=\mathrm{supp}(w)\cup\{q\}$, and $q$ is given by eventual value of the value recorder when the algorithm terminates.

We now define the (left) inverse algorithm.

\begin{definition}\label{def:inv}
	Given an indexed forest $F$, $i\in\mathbb{Z}_{>0}$. Let $u$ be a positive injective word $u=u_1\cdots u_{\ell+1}$ of the form $u=\Phi_i(w,p)$ for some $w\in\mathrm{Syl}(F)$ and $1\leq p\leq |F|+1$. We construct $w\in\mathrm{Syl}(F)$ and $p$ as follows. Let $q$ be the unique $q\notin\{\gamma(v):v\in F\}$ satisfying $\mathrm{supp}(u)=\{\gamma(v):v\in F\}\cup\{q\}$.
	\begin{enumerate}
		\item To initialize, set $w=u$, and let $j$ be the position recorder, initialized to the unique $u_j=q$.
		\item Let $F'$ be the indexed forest obtained from $F$ by removing all vertices with canonical label $\{w_1,\dots,w_{j-1}\}$.
		\begin{enumerate}
			\item  If $u_{j+1}\cdots u_{\ell+1}\in\mathrm{Syl}(F')$ and $u_j=\min\{k\geq i:k\notin\{u_{j+1},\dots,u_{\ell+1}\}\}$, then remove $u_j$ from $w$, set $p=j$, and terminate the algorithm.
			\item Otherwise, let $J\subseteq \{u_{j+1},\dots,u_{\ell+1}\}$ be the maximal consecutive interval containing $i$. Let $v$ be the root of the subtree in $F'$ with canonical labels $J$. Change $w_j$ to $\gamma(v)$, set $q=\gamma(v)$, and set $j$ to be the position of the unique $\gamma(v)$ in $\{w_{j+1},\dots,w_{\ell+1}\}$. Then, repeat step (2).
		\end{enumerate}
	\end{enumerate}
\end{definition}

\begin{proposition}\label{prop:m1}
	Given an indexed forest $F$ and $i\in\mathbb{Z}_{>0}$. If there are $w,w'\in\mathrm{Syl}(F)$ and $1\leq p,p'\leq \ell+1$ such that $\Phi_i(w,p)=\Phi_i(w',p')$ and is positive, then $w=w'$ and $p=p'$.
\end{proposition}

\begin{proof}
	It suffices to show that applying the inverse algorithms in Definition \ref{def:inv} on $\Phi_i(w,p)$ gives $(w,p)$. Let $u=\Phi_i(w,p)$, $\mathrm{supp}(u)\setminus\{\gamma(v):v\in F\}=\{q\}$, and $u_j=q$. Note that $u_k=w_k$ for $1\leq k<j$. In particular, $\{v\in F:\gamma(v)\in\{w_1,\dots,w_{j-1}\}\}$ is an upper ideal of $F$. Hence, $F'=F\setminus\{v\in F:\gamma(v)\in\{w_1,\dots,w_{j-1}\}\}$ remains to be an indexed forest. Now, there are two possibilities.
	\begin{enumerate}
		\item Suppose $j=p$ and the process of $\Phi_i(w,p)$ were to insert $i$ at position $j$ which results in $q$. In this case, to recover $w$, we simply need to remove $q$.
		\item Otherwise, $q$ comes from inserting $w_j$ to the suffix $u_{j+1}\cdots u_{\ell+1}$. In this case, $u_j$, $w_j$ and $i$ always stay in the same maximal consecutive interval in  $\{u_j,\dots,u_{\ell+1}\}$. Let $v$ be the root of the binary tree in $F'$ containing $i$, then we recover $w_j$ as $w_j$ is exactly $\gamma(v)$. The process continues until all $w_k$ are recovered.
	\end{enumerate}
	We just need to identify which possibility we are in. If $j=p$, then clearly $u_{j+1}\cdots u_{\ell+1}=w_j\cdots w_\ell\in\mathrm{Syl}(F')$, and we have $u_j=\min\{k\geq i:k\notin\{u_{j+1},\dots,u_{\ell+1}\}\}$. Conversely, if $j<p$, let $J\subseteq \{u_{j+1},\dots,u_{\ell+1}\}$ be the maximal consecutive interval containing $i$, then $u_j\in\{\min J - 1,\max J + 1\}$. If $u_j=\min J-1$, then $u_j<i$ and the condition $u_j=\min\{k\geq i:k\notin\{u_{j+1},\dots,u_{\ell+1}\}\}$ fails. If $u_j=\max J +1$, then the condition $u_{j+1}\cdots u_{\ell+1}\in\mathrm{Syl}(F')$ fails.
\end{proof}

\begin{example}
	If $w=9638245$, $i=4$ and $p=7$, then $\Phi_i(w,p)$ is obtained as follows.
	\begin{align*}
		& 9638245 \longrightarrow 963824{\color{red} 4^{[1]}}5 \xrightarrow{q=4<4^{[1]}} 96382{\color{red} 3^{[1]}}4^{[1]}5 \xrightarrow{q=3>2} \\
		& 96{\color{red} 6^{[-1]}}823^{[1]}4^{[1]}5 \xrightarrow{q=6>6^{[-1]}} 9{\color{red} 7^{[-1]}}6^{[-1]}823^{[1]}4^{[1]}5 \longrightarrow 97682345.
	\end{align*}
	The inverse construction is given as follows.
	\begin{center}
		\begin{tikzpicture}
			\node at (0,0) {\begin{tikzpicture}[scale=0.25,baseline=0pt]
					\draw[blue, thick] (-8,0) -- (-3,5) -- (2,0);
					\draw[blue, thick] (-4,4) -- (0,0);
					\draw[blue, thick] (-7,1) -- (-6,0);
					\draw[blue, thick] (-2,2) -- (-4,0);
					\draw[blue, thick] (-1,1) -- (-2,0);
					
					\filldraw[black] (-3,5) circle (5pt)  {};
					\filldraw[black] (-4,4) circle (5pt)  {};
					\filldraw[black] (-7,1) circle (5pt)  {};
					\filldraw[black] (-2,2) circle (5pt)  {};
					\filldraw[black] (-1,1) circle (5pt)  {};
					
					\draw[blue, thick] (8,0) -- (6,2) -- (4,0);
					\draw[blue, thick] (5,1) -- (6,0);
					\filldraw[black] (5,1) circle (5pt)  {};
					\filldraw[black] (6,2) circle (5pt)  {};
					
					\filldraw[black] (-10,0) circle (5pt)  {};
					\filldraw[black] (-8,0) circle (5pt)  {};
					\filldraw[black] (-6,0) circle (5pt)  {};
					\filldraw[black] (-4,0) circle (5pt)  {};
					\filldraw[black] (-2,0) circle (5pt)  {};
					\filldraw[black] (0,0) circle (5pt)  {};
					\filldraw[black] (2,0) circle (5pt)  {};
					\filldraw[black] (4,0) circle (5pt)  {};
					\filldraw[black] (6,0) circle (5pt)  {};
					\filldraw[black] (8,0) circle (5pt)  {};
					\filldraw[black] (10,0) circle (5pt)  {};
					
					\node at (-10,-1.2) {1};
					\node at (-8,-1.2) {2};
					\node at (-6,-1.2) {3};
					\node at (-4,-1.2) {4};
					\node at (-2,-1.2) {5};
					\node at (0,-1.2) {6};
					\node at (2,-1.2) {7};
					\node at (4,-1.2) {8};
					\node at (6,-1.2) {9};
					\node at (8,-1.2) {10};
					\node at (10,-1.2) {11};
					
					\draw[black] (-12,0) -- (12,0);
			\end{tikzpicture}};
			\node at (5,1) {$w$};
			\node at (8,1) {$u$};
			\node at (5,0) {$9$};
			\node at (8,0) {{\color{red} 7}682345};
			\node at (6.5,-1.5) {$\downarrow$};
			\node at (7.4,-1.5) {$\color{red} \gamma(v)=6$};
			\node at (0,-3) {\begin{tikzpicture}[scale=0.25,baseline=0pt]
					\draw[blue, thick] (-8,0) -- (-3,5) -- (2,0);
					\draw[blue, thick] (-4,4) -- (0,0);
					\draw[blue, thick] (-7,1) -- (-6,0);
					\draw[blue, thick] (-2,2) -- (-4,0);
					\draw[blue, thick] (-1,1) -- (-2,0);
					
					\draw[red] (-3,5) circle (10pt);
					\filldraw[red] (-3,5) circle (5pt)  {};
					\filldraw[black] (-4,4) circle (5pt)  {};
					\filldraw[black] (-7,1) circle (5pt)  {};
					\filldraw[black] (-2,2) circle (5pt)  {};
					\filldraw[black] (-1,1) circle (5pt)  {};
					
					\draw[blue, thick] (6,0) -- (5,1) -- (4,0);
					\filldraw[black] (5,1) circle (5pt)  {};
					
					\filldraw[black] (-10,0) circle (5pt)  {};
					\filldraw[black] (-8,0) circle (5pt)  {};
					\filldraw[black] (-6,0) circle (5pt)  {};
					\filldraw[black] (-4,0) circle (5pt)  {};
					\filldraw[black] (-2,0) circle (5pt)  {};
					\filldraw[black] (0,0) circle (5pt)  {};
					\filldraw[black] (2,0) circle (5pt)  {};
					\filldraw[black] (4,0) circle (5pt)  {};
					\filldraw[black] (6,0) circle (5pt)  {};
					\filldraw[black] (8,0) circle (5pt)  {};
					\filldraw[black] (10,0) circle (5pt)  {};
					
					\node at (-10,-1.2) {1};
					\node at (-8,-1.2) {2};
					\node at (-6,-1.2) {3};
					\node at (-4,-1.2) {4};
					\node at (-2,-1.2) {5};
					\node at (0,-1.2) {6};
					\node at (2,-1.2) {7};
					\node at (4,-1.2) {8};
					\node at (6,-1.2) {9};
					\node at (8,-1.2) {10};
					\node at (10,-1.2) {11};
					
					\draw[black] (-12,0) -- (12,0);
			\end{tikzpicture}};
			\node at (5,-3) {96};
			\node at (8,-3) {{\color{red} 6}82345};
			\node at (6.5,-4.5) {$\downarrow$};
			\node at (7.4,-4.5) {$\color{red} \gamma(v)=3$};
			\node at (0,-6) {\begin{tikzpicture}[scale=0.25,baseline=0pt]
					\draw[blue, thick] (-8,0) -- (-4,4);
					\draw[blue, thick] (-4,4) -- (0,0);
					\draw[blue, thick] (-7,1) -- (-6,0);
					\draw[blue, thick] (-2,2) -- (-4,0);
					\draw[blue, thick] (-1,1) -- (-2,0);
					
					\draw[red] (-4,4) circle (10pt);
					\filldraw[red] (-4,4) circle (5pt)  {};
					\filldraw[black] (-7,1) circle (5pt)  {};
					\filldraw[black] (-2,2) circle (5pt)  {};
					\filldraw[black] (-1,1) circle (5pt)  {};
					
					\draw[blue, thick] (6,0) -- (5,1) -- (4,0);
					\filldraw[black] (5,1) circle (5pt)  {};
					
					\filldraw[black] (-10,0) circle (5pt)  {};
					\filldraw[black] (-8,0) circle (5pt)  {};
					\filldraw[black] (-6,0) circle (5pt)  {};
					\filldraw[black] (-4,0) circle (5pt)  {};
					\filldraw[black] (-2,0) circle (5pt)  {};
					\filldraw[black] (0,0) circle (5pt)  {};
					\filldraw[black] (2,0) circle (5pt)  {};
					\filldraw[black] (4,0) circle (5pt)  {};
					\filldraw[black] (6,0) circle (5pt)  {};
					\filldraw[black] (8,0) circle (5pt)  {};
					\filldraw[black] (10,0) circle (5pt)  {};
					
					\node at (-10,-1.2) {1};
					\node at (-8,-1.2) {2};
					\node at (-6,-1.2) {3};
					\node at (-4,-1.2) {4};
					\node at (-2,-1.2) {5};
					\node at (0,-1.2) {6};
					\node at (2,-1.2) {7};
					\node at (4,-1.2) {8};
					\node at (6,-1.2) {9};
					\node at (8,-1.2) {10};
					\node at (10,-1.2) {11};
					
					\draw[black] (-12,0) -- (12,0);
			\end{tikzpicture}};
			\node at (5,-6) {96382};
			\node at (8,-6) {{\color{red} 3}45};
			\node at (6.5,-7.5) {$\downarrow$};
			\node at (7.4,-7.5) {$\color{red} \gamma(v)=4$};
			\node at (0,-9) {\begin{tikzpicture}[scale=0.25,baseline=0pt]
					
					\draw[blue, thick] (0,0) -- (-2,2) -- (-4,0);
					\draw[blue, thick] (-1,1) -- (-2,0);
					\draw[red] (-2,2) circle (10pt);
					\filldraw[red] (-2,2) circle (5pt)  {};
					\filldraw[black] (-1,1) circle (5pt)  {};
					
					\filldraw[black] (-10,0) circle (5pt)  {};
					\filldraw[black] (-8,0) circle (5pt)  {};
					\filldraw[black] (-6,0) circle (5pt)  {};
					\filldraw[black] (-4,0) circle (5pt)  {};
					\filldraw[black] (-2,0) circle (5pt)  {};
					\filldraw[black] (0,0) circle (5pt)  {};
					\filldraw[black] (2,0) circle (5pt)  {};
					\filldraw[black] (4,0) circle (5pt)  {};
					\filldraw[black] (6,0) circle (5pt)  {};
					\filldraw[black] (8,0) circle (5pt)  {};
					\filldraw[black] (10,0) circle (5pt)  {};
					
					\node at (-10,-1.2) {1};
					\node at (-8,-1.2) {2};
					\node at (-6,-1.2) {3};
					\node at (-4,-1.2) {4};
					\node at (-2,-1.2) {5};
					\node at (0,-1.2) {6};
					\node at (2,-1.2) {7};
					\node at (4,-1.2) {8};
					\node at (6,-1.2) {9};
					\node at (8,-1.2) {10};
					\node at (10,-1.2) {11};
					
					\draw[black] (-12,0) -- (12,0);
			\end{tikzpicture}};
			\node at (5,-9) {963824};
			\node at (8,-9) {{\color{red} 4}5};
			\node at (5,-9.6) {9638245};
		\end{tikzpicture}
	\end{center}
\end{example}

\subsection{A Monk rule for slide polynomials}

\begin{definition}
	Given an injective decorated word $w=w_1\cdots w_\ell$ and a word $\mathbf{a}=a_1\cdots a_\ell$, we define
	$$\CMcal{C}(w;\mathbf{a})=\{b_1\cdots b_\ell\in\mathbb{Z}_{>0}^\ell:b_1\leq\cdots\leq b_\ell, b_k\leq a_k, w_k<w_{k+1}\implies b_k<b_{k+1}\}.$$
	In particular, $\CMcal{C}(w,w)=\E(w)$ for injective word $w$.
	
	We define the \textbf{tight bounds} of $(w,\mathbf{a})$ to be the two words $\epsilon_1\cdots\epsilon_{\ell-1}$ and $\delta_1\cdots\delta_\ell$ where $\epsilon_k=\begin{cases}
		1 & w_k<w_{k+1} \\ 
		0 & w_k>w_{k+1}
	\end{cases}$ and $\delta_k=\min\{a_k,\delta_{k+1}-\epsilon_k\}$ with initial condition $\delta_\ell=a_\ell$.
\end{definition}

\begin{lemma}\label{lem:m1-3}
	The set $\CMcal{C}(w;\mathbf{a})$ is given by $\{b_1\cdots b_\ell\in\mathbb{Z}_{>0}^\ell:b_k+\epsilon_k\leq b_{k+1},b_k\leq\delta_k\}$.
\end{lemma}

\begin{proof}
	The conditions $b_k\leq b_{k+1}$ and $w_k<w_{k+1}\implies b_k<b_{k+1}$ is equivalent to $b_k+\epsilon_k\leq b_{k+1}$. Clearly $b_k\leq\delta_k\implies b_k\leq a_k$. Conversely, if $b_k\leq a_k$, then $b_\ell\leq a_\ell=\delta_\ell$, and for $k<\ell$, $b_k\leq b_{k+1}-\epsilon_k\leq\delta_{k+1}-\epsilon_k$ by induction.
\end{proof}

\begin{lemma}\label{lem:m1-4}
	Given a positive injective word $w=w_1\cdots w_\ell$, $i\in\mathbb{Z}_{>0}$ and $1\leq p\leq \ell+1$, let $\mathbf{a}=(w_1,\dots,w_{p-1},i,w_p,\dots,w_\ell)$. Then, we have
	$$\E(\Phi_i(w,p))=\CMcal{C}(\mathrm{Ins}_i(w,p);\mathbf{a}).$$
\end{lemma}

\begin{proof}
	Let $(\epsilon^L,\delta^L)$ and $(\epsilon^R,\delta^R)$ be the tight bounds for the left and right hand sides, respectively. By Lemma \ref{lem:m1-3}, it suffices to show that $(\epsilon^L,\delta^L)=(\epsilon^R,\delta^R)$.
	
	We first prove $\epsilon^L_k=\epsilon^R_k$ i.e. $\mathrm{Ins}_i(w,p)_k<\mathrm{Ins}_i(w,p)_{k+1}$ if and only if $\Phi_i(w,p)_k<\Phi_i(w,p)_{k+1}$. Note that $\mathrm{Ins}_i(w,p)=\cdots A^{[1]}w_p\cdots w_\ell$ and $\Phi_i(w,p)=\cdots Aw_p\cdots w_\ell$ with $A\notin\{w_p,\dots,w_\ell\}$. The statement holds for $k\geq p$.
	
	For $k<p$, we have $\mathrm{Ins}_i(w,p)_k=w_k$. Let $S_k=\{\Phi_i(w,p)_{k+1},\dots,\Phi_i(w,p)_{\ell+1}\}$ be the suffix support. We decompose $S_k$ into maximal consecutive intervals, and assume $\Phi_i(w,p)_{k+1}$ belongs to an interval $J$. By construction, $\mathrm{val}(\mathrm{Ins}_i(w,p)_{k+1})\in J$. 
	\begin{enumerate}
		\item Suppose $w_k\notin S_k$. In this case, $\Phi_i(w,p)_k=w_k$, and $\Phi_i(w,p)_k<\Phi_i(w,p)_{k+1}$ if and only if $w_k<\min J$ if and only if $w_k<\mathrm{Ins}_i(w,p)_{k+1}$.
		\item Suppose $w_k$ belongs to an interval $I$ of $S_k$ with $I\neq J$. In this case, $\Phi_i(w,p)_k\in\{\min I-1,\max I+1\}$, and $\Phi_i(w,p)_k<\Phi_i(w,p)_{k+1}$ if and only if $\max I<\min J$ if and only if $w_k<\mathrm{Ins}_i(w,p)_{k+1}$.
		\item Suppose $w_k\in J$. If $w_k<\mathrm{Ins}_i(w,p)_{k+1}$, then either $\Phi_i(w,p)_k=\min J-1<\Phi_i(w,p)_{k+1}$, or $w_k=\Phi_i(w,p)_{k+1}=\min J$ while in the intermediate step, $w_k<(\min J)^{[1]}$, and hence $\Phi_i(w,p)_k=\min J-1$. The case $w_k>\mathrm{Ins}_i(w,p)_{k+1}$ is analogous.
	\end{enumerate}
	
	We now prove $\delta^L_k=\delta^R_k$ i.e. $\min\left\{\Phi_i(w,p)_k,\delta^L_{k+1}-\epsilon^L_k\right\}=\min\left\{a_k,\delta^R_{k+1}-\epsilon^R_k\right\}$. We prove by induction on suffix. The base case is $\delta^L_{\ell+1}=\Phi_i(w,p)_{\ell+1}=a_{\ell+1}=\delta^R_{\ell+1}$. By induction, we have $\delta_{k+1}^L-\epsilon_k^L=\delta_{k+1}^R-\epsilon_k^R$.
	\begin{enumerate}
		\item When $p<k\leq\ell$, we have $\Phi_i(w,p)_k=w_{k-1}=a_k$. Hence, $\delta_k^L=\delta_k^R$.
		\item When $k=p$, let $S$ be the suffix support $\{w_p,\dots,w_\ell\}$. If $i\notin S$, then $\Phi_i(w,p)_p=i=a_p$ and we are done. If $i\in S$, we can find $j\geq p$ with $w_j=i$. Then, $\delta^L_{p+1}\leq\cdots\leq \delta_{j+1}^L\leq a_{j+1}=w_j=i$. In particular, $\delta_{p+1}^L-\epsilon_p^L\leq i=a_p<\Phi_i(w,p)_p$. Therefore, $\delta^L_p=\delta^R_p=\delta^L_{p+1}-\epsilon^L_p$.
		\item When $k<p$, let $S=\{\Phi_i(w,p)_{k+1},\dots,\Phi_i(w,p)_{\ell+1}\}$ be the suffix support. If $w_k\notin S$, then $\Phi_i(w,p)_k=w_k=a_k$, and we are done. If $w_k\in S$, let $J$ be the maximal consecutive interval in $S$ containing $w_k$. There are two sub-cases.
		\begin{enumerate}
			\item If $\Phi_i(w,p)_k=\max J+1>a_k$, we can find $\Phi_i(w,p)_j=w_k$ with $j>k$. In particular, $\delta^L_{k+1}\leq\cdots\leq\delta^L_j\leq \Phi_i(w,p)_j=w_k$, and hence, $\delta^L_{k+1}-\epsilon^L_k\leq w_k=a_k<\Phi_i(w,p)_k$. Therefore, $\delta_k^L=\delta_k^R=\delta_{k+1}^L-\epsilon^L_k$.
			\item If $\Phi_i(w,p)_k=\min J-1<a_k$, we can find $\Phi_i(w,p)_j=\min J$ with $j>k$. Note that there must exists an ascent $\Phi_i(w,p)_{q}<\Phi_i(w,p)_{q+1}$ for some $k\leq q\leq j-1$ i.e. $\epsilon_q^L=1$. Then, either $q=k$, or $q>k$ and $\delta^L_{k+1}\leq\cdots\leq\delta_q<\delta_{q+1}\leq\cdots\leq\delta^L_j\leq \Phi_i(w,p)_j=\min J$. In both cases, we have
			$$\delta^L_{k+1}-\epsilon^L_k\leq \delta_j^L-\sum_{r=k}^{j-1}\epsilon_r^L \leq \min J-1=\Phi_i(w,p)_k<a_k.$$
			Therefore, $\delta_k^L=\delta_k^R=\delta_{k+1}^L-\epsilon^L_k$.
		\end{enumerate}
	\end{enumerate}
\end{proof}

\begin{corollary}\label{cor:m1}
	For positive injective word $w=w_1\cdots w_\ell$, $i\in\mathbb{Z}_{>0}$ and $1\leq p\leq \ell+1$, if $\Phi_i(w,p)$ is legal, then $w$ is legal.
\end{corollary}

\begin{proof}
	If $\Phi_i(w,p)$ is legal, then $\E(\Phi_i(w,p))\neq\varnothing$. By Lemma \ref{lem:m1-4}, we can find $b_1\cdots b_{\ell+1}\in \CMcal{C}(\mathrm{Ins}_i(w,p);\mathbf{a})$. Then, $b_1\cdots b_{p-1}b_{p+1}\cdots b_{\ell+1}\in\E(w)$ i.e. $w$ is legal.
\end{proof}

\begin{lemma}\label{lem:m1-5}
	Given a legal word $w=w_1\cdots w_\ell$ and $i\in\mathbb{Z}_{>0}$, fix $1\leq q\leq i$ and $\mathbf{b}=b_1\cdots b_\ell\in\E(w)$. There exists a unique $p$ such that $\mathbf{b}'=(b_1,\dots,b_{p-1},q,b_p,\dots,b_\ell)\in\CMcal{C}(\mathrm{Ins}_{i}(w,p),\mathbf{a})$ where $\mathbf{a}=(w_1,\dots,w_{p-1},i,w_p,\dots,w_\ell)$.
\end{lemma}

\begin{proof}
	The condition $b'_k\leq a_k$ is clearly satisfied for all $k$. If $q\notin\{b_1,\dots,b_\ell\}$, then the condition $b'_k\leq b'_{k+1}$ determines a unique position $p$ to put $q$, and we are done.
	
	We now assume $q\in\{b_1,\dots,b_\ell\}$. Let $L=\min\{j:b_j=q\}$ and $R=\max\{j:b_j=q\}$, then we have $w_L>w_{L+1}>\cdots>w_R$. To satisfy the condition $b'_k\leq b'_{k+1}$, we must have $L\leq p\leq R+1$. Note that $\mathrm{Ins}_i(w,p)=w_1\cdots w_{p-1}A_p^{[1]}w_p\cdots w_\ell$ where $A_p=\min\{t\geq i:t\notin\{w_p,\dots,w_\ell\}\}$. We need to show there exists a unique $p\in[L,R+1]$ such that $w_{p-1} > A_p$ if $L<p\leq R+1 $ and $A_p\geq w_p$ if $L\leq p\leq R$.
	
	Recall that the insertion algorithm guarantees $A_p\neq w_p$. It is straightforward to check that $A_{p+1}=\begin{cases}
		A_p & w_p<i\text{ or }w_p>A_p \\
		w_p & i\leq w_p < A_p
	\end{cases}$. We set $p\in [L,R]$ to be the smallest number such that $A_p> w_p$, or $p=R+1$ if no such number exists. Note that $w_k>A_k$ for $L\leq k< p$ implies $A_L=A_{L+1}=\cdots=A_p$. In particular, we get the desired condition $p>L\implies w_{p-1}>A_{p-1}=A_p$.
	
	By construction, the desired condition $A_k> w_k$ fails for any $L\leq k<p$. For any $p<k\leq R+1$, note that $A_{k-1}>w_{k-1}$ implies either $A_k=A_{k-1}$ (when $w_{k-1}<i$) or $A_k=w_{k-1}$ (when $w_{k-1}\geq i$). In both cases, we have $A_k\geq w_{k-1}>w_k$, which leads to $A_{k+1}\geq w_k>w_{k+1}$ etc. In particular, the desired condition $w_{k-1}>A_k$ fails. Therefore, $p$ is unique. 
\end{proof}

\begin{lemma}\label{lem:m1-6}
	Given a positive injective word $w=w_1\cdots w_\ell$ and $i\in\mathbb{Z}_{>0}$. The map
	$$[i]\times \E(w)\to\biguplus_{p=1}^{\ell+1}\E(\Phi_i(w,p))$$
	defined by $(q,\mathbf{b})\mapsto (p,(b_1,\dots,b_{p-1},q,b_p,\dots,b_\ell))$, with $p$ from Lemma \ref{lem:m1-5}, is a bijection.
\end{lemma}

\begin{proof}
	Given $(p,b_1\cdots b_{\ell+1})$ from the right hand side, Lemma \ref{lem:m1-4} and \ref{lem:m1-5} guarantee that taking $q=b_p$ and $b=b_1\cdots b_{p-1}b_{p+1}\cdots b_{\ell+1}\in\E(w)$ gives the inverse map.
\end{proof}

\begin{definition}
	For a positive injective word $w=w_1\cdots w_\ell$, define $\mathcal{A}_i(w)=\{1\leq p\leq \ell+1:\Phi_i(w,p)\text{ is legal}\}$.
\end{definition}

\begin{theorem}\label{thm:slide}
	For positive injective word $w$, we have $(x_1+\cdots+x_i)\mathfrak{F}_w=\sum_{p\in \mathcal{A}_i(w)}\mathfrak{F}_{\Phi_i(w,p)}$.
\end{theorem}

\begin{proof}
	It follows directly from Lemma \ref{lem:m1-6}.
\end{proof}

\subsection{Monk rule for forest polynomials}

\begin{lemma}\label{lem:m1-8}
	Given indexed forests $F,G$ and $i\in\mathbb{Z}_{>0}$. Suppose there exist $w'\in\mathrm{Syl}(F)$ and $p'\in\mathcal{A}_i(w')$ with $\Phi_i(w',p')\in\mathrm{Syl}(G)$. Then, for any $u\in\mathrm{Syl}(G)$, there exist $w\in\mathrm{Syl}(F)$ and $p\in\mathbb{Z}_{>0}$ such that $\Phi_i(w,p)=u$.
\end{lemma}

\begin{proof}
	Let $\Phi_i(w',p')=u_1\cdots u_{\ell+1}$. We first assume that
	$$u=s_k\cdot \Phi_i(w',p')=u_1\cdots u_{k-1}u_{k+1}u_ku_{k+2}\cdots u_{\ell+1}$$
	for some $k$ satisfying that the maximal consecutive interval in the suffix support $\{(u_k,\dots,u_{\ell+1}\})$ containing $u_k$ does not contain $u_{k+1}$. In the other words, the nodes $v,v'$ are incomparable in $G$ where $\gamma(v)=u_k$ and $\gamma(v')=u_{k+1}$.
	
	Note that the whole set $\mathrm{Syl}(G)$ can be obtained from $\Phi_i(w',p')$ via a sequence of such simple transpositions. This follows from the well-known fact that any two linear extensions are connected by simple transpositions on incomparable elements. Therefore, it suffices to check that there exists $w\in\mathrm{Syl}(F)$ and $p\in\mathbb{Z}_{\geq 0}$ with $\Phi_i(w,p)=u$. There are the four cases.
	\begin{enumerate}
		\item If $k<p'-1$, set $w=s_k\cdot w'$ and $p=p'$.
		\item If $k=p'-1$, set $w=w'$ and $p=p'-1$.
		\item If $k=p'$, set $w=w'$ and $p=p'+1$.
		\item If $k>p'$, set $w=s_{k-1}\cdot w'$ and $p=p'$.
	\end{enumerate}
	In all these cases, the incomparability condition on $k$ implies $\Phi_i(w,p)=s_k\cdot\Phi_i(w',p')=u$.
\end{proof}

\begin{definition}
	For an indexed forest $F$ and $i\in\mathbb{Z}_{>0}$, define $\Omega_i(F)=\bigcup_{w\in\mathrm{LSyl}(F)}\bigcup_{p\in \mathcal{A}_i(w)}\{G:\Phi_i(w,p)\in\mathrm{LSyl}(G)\}$.
\end{definition}

\begin{lemma}\label{lem:m1-9}
	Given indexed forests $F,G$ and $i\in\mathbb{Z}_{>0}$. Suppose there exist $w'\in\mathrm{LSyl}(F)$ and $p'\in\mathcal{A}_i(w')$ with $\Phi_i(w',p')\in\mathrm{LSyl}(G)$. Then, the map
	$$\biguplus_{w\in\mathrm{LSyl}(F)}\biguplus_{p\in\mathcal{A}_i(w)}\left(\{\Phi_i(w,p)\}\cap\mathrm{LSyl}(G)\right)\to \mathrm{LSyl}(G)$$
	given by $(w,p,u)\mapsto u$ is a bijection.
	
	Consequently, given $F$ and $i$, we have the disjoint union decomposition
	$$\biguplus_{G\in\Omega_i(F)}\mathrm{LSyl}(G)=\biguplus_{w\in\mathrm{LSyl}(F)}\biguplus_{p\in\mathcal{A}_i(w)}\{\Phi_i(w,p)\}.$$
\end{lemma}

\begin{proof}
	Fix $u\in\mathrm{LSyl}(G)$. Lemma \ref{lem:m1-8} gives some $w\in\mathrm{Syl}(F)$ and $p\in\mathbb{Z}_{>0}$ such that $\Phi_i(w,p)=u$. Since $u$ is legal, Corollary \ref{cor:m1} implies $w$ is legal, and $p\in\mathcal{A}_i(w)$. This proves surjectivity. Injectivity follows from Proposition \ref{prop:m1}.
\end{proof}

\begin{theorem}[Monk rule]\label{thm:monk}
	We have $(x_1+\cdots+x_i)\fP_F=\sum_{G\in\Omega_i(F)}\fP_G$. In particular, all coefficients are multiplicity-free.
\end{theorem}

\begin{proof}
	By Lemma \ref{lem:m1-2}, Theorem \ref{thm:slide} and Lemma \ref{lem:m1-9}, we have
	\begin{align*}
		(x_1+\cdots+x_i)\fP_F&=\sum_{w\in\mathrm{LSyl}(F)}(x_1+\cdots+x_i)\mathfrak{F}_w\\
		&=\sum_{w\in\mathrm{LSyl}(F)}\sum_{p\in\mathcal{A}_i(w)}\mathfrak{F}_{\Phi_i(w,p)}\\
		&=\sum_{G\in\Omega_i(F)}\sum_{w\in\mathrm{LSyl}(G)}\mathfrak{F}_w\\
		&=\sum_{G\in\Omega_i(F)}\fP_G.
	\end{align*}
\end{proof}

\medskip
\noindent
{\bf Acknowledgement.} During the research and preparation of this work, the authors used OpenAI model Astra 6.1 to assist with literature searches, mathematical exploration, the development
and testing of computational code, and critical review of the manuscript.

The authors directed the research, critically reviewed, rewritten and revised
the resulting material, and take full responsibility for the final manuscript,
including its mathematical arguments, computational results, and references.


\begin{thebibliography}{99}
	\bibitem{AS}
	S.~Assaf and D.~Searles,
	\emph{Schubert polynomials, slide polynomials, Stanley symmetric functions and quasi-Yamanouchi pipe dreams},
	Adv. Math. \textbf{306} (2017), 89--122.
	\href{https://doi.org/10.1016/j.aim.2016.10.015}
	{doi:10.1016/j.aim.2016.10.015}.
	
	\bibitem{ABB}
	J.-C.~Aval, F.~Bergeron, and N.~Bergeron,
	\emph{Ideals of quasi-symmetric functions and super-covariant polynomials for $S_n$},
	Adv. Math. \textbf{181} (2004), no.~2, 353--367.
	\href{https://doi.org/10.1016/S0001-8708(03)00068-9}
	{doi:10.1016/S0001-8708(03)00068-9}.
	
	\bibitem{Bansal}
	R.~Bansal,
	\emph{A quasisymmetric analogue of the Bruhat order},
	arXiv:2611.02210 (2026).
	\url{https://arxiv.org/abs/2611.02210}.
	
	\bibitem{BGNST}
	N.~Bergeron, L.~Gagnon, P.~Nadeau, H.~Spink, and V.~Tewari,
	\emph{Equivariant quasisymmetry and noncrossing partitions},
	arXiv:2504.15234 (2025).
	\url{https://arxiv.org/abs/2504.15234}.
	
	\bibitem{BGNSTQF}
	N.~Bergeron, L.~Gagnon, P.~Nadeau, H.~Spink, and V.~Tewari,
	\emph{The quasisymmetric flag variety: A toric complex on noncrossing partitions},
	Forum Math. Pi \textbf{14} (2026), e18.
	\href{https://doi.org/10.1017/fmp.2026.10034}
	{doi:10.1017/fmp.2026.10034}.
	
	\bibitem{BMSW}
	N.~Bergeron, S.~Mykytiuk, F.~Sottile, and S.~van Willigenburg,
	\emph{Non-commutative Pieri operators on posets},
	J. Combin. Theory Ser. A \textbf{91} (2000), nos.~1--2, 84--110.
	\href{https://doi.org/10.1006/jcta.2000.3090}
	{doi:10.1006/jcta.2000.3090}.
	
	\bibitem{BS}
	N.~Bergeron and F.~Sottile,
	\emph{Schubert polynomials, the Bruhat order, and the geometry of flag manifolds},
	Duke Math. J. \textbf{95} (1998), no.~2, 373--423.
	\href{https://doi.org/10.1215/S0012-7094-98-09511-4}
	{doi:10.1215/S0012-7094-98-09511-4}.
	
	\bibitem{Billey}
	S.~C.~Billey,
	\emph{Kostant polynomials and the cohomology ring for $G/B$},
	Duke Math. J. \textbf{96} (1999), no.~1, 205--224.
	\href{https://doi.org/10.1215/S0012-7094-99-09606-0}
	{doi:10.1215/S0012-7094-99-09606-0}.
	
	\bibitem{Gessel}
	I.~M.~Gessel,
	\emph{Multipartite $P$-partitions and inner products of skew Schur functions},
	in \emph{Combinatorics and Algebra} (Boulder, CO, 1983),
	Contemp. Math., vol.~34, Amer. Math. Soc., Providence, RI, 1984,
	pp.~289--317.
	\href{https://doi.org/10.1090/conm/034/777705}
	{doi:10.1090/conm/034/777705}.
	
	\bibitem{GW}
	T.~Gobet and N.~Williams,
	\emph{Noncrossing partitions and Bruhat order},
	European J. Combin. \textbf{53} (2016), 8--34.
	\href{https://doi.org/10.1016/j.ejc.2015.10.007}
	{doi:10.1016/j.ejc.2015.10.007}.
	
	\bibitem{GKM}
	M.~Goresky, R.~Kottwitz, and R.~MacPherson,
	\emph{Equivariant cohomology, Koszul duality, and the localization theorem},
	Invent. Math. \textbf{131} (1998), 25--83.
	\href{https://doi.org/10.1007/s002220050197}
	{doi:10.1007/s002220050197}.
	
	\bibitem{Graham}
	W.~Graham,
	\emph{Positivity in equivariant Schubert calculus},
	Duke Math. J. \textbf{109} (2001), no.~3, 599--614.
	\href{https://doi.org/10.1215/S0012-7094-01-10935-6}
	{doi:10.1215/S0012-7094-01-10935-6}.
	
	\bibitem{HNT}
	F.~Hivert, J.-C.~Novelli, and J.-Y.~Thibon,
	\emph{The algebra of binary search trees},
	Theoret. Comput. Sci. \textbf{339} (2005), no.~1, 129--165.
	\href{https://doi.org/10.1016/j.tcs.2005.01.012}
	{doi:10.1016/j.tcs.2005.01.012}.
	
	\bibitem{LS}
	A.~Lascoux and M.-P.~Sch\"utzenberger,
	\emph{Polyn\^omes de Schubert},
	C. R. Acad. Sci. Paris S\'er. I Math. \textbf{294} (1982), 447--450.
	\url{https://www-igm.univ-mlv.fr/~berstel/Mps/Travaux/A/1982-1SchubertCras.pdf}.
	
	\bibitem{Loday}
	J.-L.~Loday,
	\emph{Realization of the Stasheff polytope},
	Arch. Math. (Basel) \textbf{83} (2004), 267--278.
	\href{https://doi.org/10.1007/s00013-004-1026-y}
	{doi:10.1007/s00013-004-1026-y}.
	
	\bibitem{LR}
	J.-L.~Loday and M.~O.~Ronco,
	\emph{Hopf algebra of the planar binary trees},
	Adv. Math. \textbf{139} (1998), no.~2, 293--309.
	\href{https://doi.org/10.1006/aima.1998.1759}
	{doi:10.1006/aima.1998.1759}.
	
	\bibitem{Monk}
	D.~Monk,
	\emph{The geometry of flag manifolds},
	Proc. London Math. Soc. (3) \textbf{9} (1959), no.~2, 253--286.
	\href{https://doi.org/10.1112/plms/s3-9.2.253}
	{doi:10.1112/plms/s3-9.2.253}.
	
	\bibitem{NST}
	P.~Nadeau, H.~Spink, and V.~Tewari,
	\emph{Quasisymmetric divided differences},
	arXiv:2406.01510 (2024), revised 2026.
	\url{https://arxiv.org/abs/2406.01510}.
	
	\bibitem{NSTGeometry}
	P.~Nadeau, H.~Spink, and V.~Tewari,
	\emph{The geometry of quasisymmetric coinvariants},
	arXiv:2410.12643 (2024).
	\url{https://arxiv.org/abs/2410.12643}.
	
	\bibitem{NSTExpansions}
	P.~Nadeau, H.~Spink, and V.~Tewari,
	\emph{Schubert polynomial expansions revisited},
	Forum Math. Sigma \textbf{13} (2025), e106.
	\href{https://doi.org/10.1017/fms.2025.10068}
	{doi:10.1017/fms.2025.10068}.
	
	\bibitem{NSTGrove}
	P.~Nadeau, H.~Spink, and V.~Tewari,
	\emph{Grove polynomials and $K$-theoretic quasisymmetry},
	arXiv:2605.22750 (2026).
	\url{https://arxiv.org/abs/2605.22750}.
	
	\bibitem{NT}
	P.~Nadeau and V.~Tewari,
	\emph{Forest polynomials and the class of the permutahedral variety},
	Adv. Math. \textbf{453} (2024), 109834.
	\href{https://doi.org/10.1016/j.aim.2024.109834}
	{doi:10.1016/j.aim.2024.109834}.
	
	\bibitem{NTFlags}
	P.~Nadeau and V.~Tewari,
	\emph{$P$-partitions with flags and back stable quasisymmetric functions},
	Combin. Theory \textbf{4} (2024), no.~2, Paper No.~4.
	\href{https://doi.org/10.5070/C64264232}
	{doi:10.5070/C64264232}.
	
	\bibitem{Postnikov}
	A.~Postnikov,
	\emph{Permutohedra, associahedra, and beyond},
	Int. Math. Res. Not. IMRN \textbf{2009} (2009), no.~6, 1026--1106.
	\href{https://doi.org/10.1093/imrn/rnn153}
	{doi:10.1093/imrn/rnn153}.
	
	\bibitem{Sottile}
	F.~Sottile,
	\emph{Pieri's formula for flag manifolds and Schubert polynomials},
	Ann. Inst. Fourier (Grenoble) \textbf{46} (1996), no.~1, 89--110.
	\href{https://doi.org/10.5802/aif.1508}{doi:10.5802/aif.1508}.
\end{thebibliography}
\end{document}